\RequirePackage{fix-cm}
\documentclass[final]{svjour3}       
\usepackage{springerpaper}
\newcommand{\cahiernumber}{00}  
\smartqed  
\usepackage{graphicx}
\usepackage{comment}
\usepackage{subcaption}
\usepackage{tcolorbox}
\usepackage{multirow}
\usepackage{pgfplots}
\usepackage{cases}
\usepackage{tikz}                  
\usetikzlibrary{arrows,positioning,shapes}
\newcommand{\skcp}{s_{k, \textup{cp}}}
\newcommand{\kdc}{\kappa_{\textup{mdc}}}
\newcommand{\ueps}{u_{\epsilon}}
\newcommand{\fpsib}{\psi_\beta}
\newcommand{\psib}[1]{\psi_\beta\left(#1\right)}

\makeatletter
\newcommand*{\defeq}{\mathrel{\vcenter{\baselineskip0.5ex \lineskiplimit0pt
                     \hbox{\scriptsize.}\hbox{\scriptsize.}}}%
                     =}
\makeatother

\def\cahiernumber{43}

\newcommand{\papertitle}{%
  Complexity of trust-region methods in the presence of unbounded Hessian approximations
}

\hypersetup{
  pdftitle={\papertitle},
  pdfauthor={Youssef Diouane and Mohamed Laghdaf Habiboullah and Dominique Orban},
  pdfsubject={Complexity of trust-region methods in the presence of unbounded Hessian approximations},
  pdfkeywords={Complexity analysis, Quasi-Newton methods, Trust-region methods, Unconstrained optimization, Nonconvex optimization, Convex optimization},
}

\title{\papertitle}
\author{%
  Youssef Diouane\footnote{%
    GERAD and Department of Mathematics and Industrial Engineering, Polytechnique Montr\'eal. E-mail: \href{mailto:youssef.diouane@polymtl.ca}{youssef.diouane@polymtl.ca}.
  }
   \thanks{Research supported by an NSERC Discovery grant.}
  \and
  Mohamed Laghdaf Habiboullah\footnote{%
    GERAD and Department of Mathematics and Industrial Engineering, Polytechnique Montr\'eal. E-mail: \href{mailto:mohamed.habiboullah@polymtl.ca}{mohamed.habiboullah@polymtl.ca}.
  }
  \and
  Dominique Orban\footnote{%
    GERAD and Department of Mathematics and Industrial Engineering, Polytechnique Montr\'eal. E-mail: \href{mailto:dominique.orban@gerad.ca}{dominique.orban@gerad.ca}.
  }
  \thanks{Research supported by an NSERC Discovery grant.}
}

\spdefaulttheorem{problemassumption}{Problem Assumption}{\bf}{\it}
\spdefaulttheorem{modelassumption}{Model Assumption}{\bf}{\it}
\spdefaulttheorem{stepassumption}{Step Assumption}{\bf}{\it}
\crefname{problemassumption}{Problem Assumption}{Problem Assumptions}
\Crefname{problemassumption}{Problem Assumption}{Problem Assumptions}
\crefname{modelassumption}{Model Assumption}{Model Assumptions}
\Crefname{modelassumption}{Model Assumption}{Model Assumptions}
\crefname{stepassumption}{Step Assumption}{Step Assumptions}
\Crefname{stepassumption}{Step Assumption}{Step Assumptions}
\Crefname{subsection}{Section}{Sections}
\crefname{subsection}{section}{sections}
\journalname{Mathematical Programming}
\begin{document}

\title{%
  \papertitle
  \thanks{Research supported by an NSERC Discovery grant and Fonds de recherche du Québec grant}
}

\titlerunning{Complexity of TR Methods with Unbounded Model Hessians}

\author{%
  Youssef Diouane \and
  Mohamed Laghdaf Habiboullah \and
  Dominique Orban
}

\authorrunning{Diouane, Habiboullah, Orban} 

\institute{%
  Y.~Diouane \at
  GERAD and Mathematics and Industrial Engineering, Polytechnique Montr\'eal.\\
  \email{youssef.diouane@polymtl.ca}
  \and
  M.L.~Habiboullah \at
  GERAD and Mathematics and Industrial Engineering, Polytechnique Montr\'eal.\\
  \email{mohamed-laghdaf-2.habiboullah@polymtl.ca}
  \and
  D.~Orban \at
  GERAD and Mathematics and Industrial Engineering, Polytechnique Montr\'eal.\\
  \email{dominique.orban@gerad.ca}
}

\date{Received: date / Accepted: date}

\pagestyle{myheadings}

\maketitle
\thispagestyle{mytitlepage}

\begin{abstract}
  We extend traditional complexity analyses of trust-region methods for unconstrained, possibly nonconvex, optimization.
  Whereas most complexity analyses assume uniform boundedness of the model Hessians, we work with potentially unbounded model Hessians.
  Boundedness is not guaranteed in practical implementations, in particular ones based on quasi-Newton updates such as PSB, BFGS and SR1.
  Our analysis is conducted for a family of trust-region methods that includes most known methods as special cases.
  We examine two regimes of Hessian growth: one bounded by a power of the number of successful iterations, and one bounded by a power of the number of iterations.
  This allows us to formalize and address the intuition of \citet{powell-2010}, who studied convergence under a special case of our assumptions, but whose proof contained complexity arguments.
  Specifically, for \(0 \leq p < 1\), we establish sharp \(O([(1-p)\epsilon^{-2}]^{1/(1-p)})\) evaluation complexity to find an \(\epsilon\)-stationary point when model Hessians are \(O(|\mathcal{S}_{k-1}|^p)\), where \(|\mathcal{S}_{k-1}|\) is the number of iterations where the step was accepted, up to iteration \(k-1\).
  For \(p = 1\), which is the case studied by Powell, we establish a sharp \(O(\exp(c_1\epsilon^{-2}))\) evaluation complexity for a certain constant \(c_1 > 0\).
  This is far better than the double exponential bound that \citet{powell-2010} suspected, and is far worse than other bounds surmised elsewhere in the literature.
  We establish similar sharp bounds when model Hessians are \(O(k^p)\), where \(k\) is the iteration counter, for \(0 \leq p < 1\).
  When \(p = 1\), the complexity bound depends on the parameters of the family, but reduces to \(O((1 - \log(\epsilon))\exp(c_2\epsilon^{-2}))\) for a certain constant \(c_2 > 0\) for the special case of the standard trust-region method.
  As special cases, we derive novel complexity bounds for (strongly) convex objectives under the same growth assumptions.
\end{abstract}

\bigskip

\noindent

\section{Introduction}%

We consider the solution of the nonconvex unconstrained problem
\begin{equation}%
  \label{eq:nlp}
  \minimize{x \in \R^n} \ f(x),
\end{equation}
where \(f: \R^n \rightarrow \R\) is continuously differentiable, by means of trust-region methods~\citep{conn-gould-toint-2000}.
At each iteration of a trust-region method, a model of the objective is used to compute a step.
In most cases, said model is a quadratic, and is required to match \(f\) up to first order at the current iterate.
Although there are variants allowing inexact first-order information \citep[\S\(8.4\)]{conn-gould-toint-2000}, there is far more freedom in the model Hessian \(B_k = B_k^T\) at iteration \(k\).
In particular, quasi-Newton approximations are a natural choice as they are known to yield fast local convergence under certain assumptions \citep{dennis-more-1977}.
Convergence of trust-region methods, both local and global, has been studied extensively for the past few decades---see \citep{conn-gould-toint-2000} and references therein---and their worst-case evaluation complexity has attracted much attention in the past decade and a half---see \citep{cartis-gould-toint-2022} and references therein.
Although convergence of many trust-region schemes has been established under general assumptions, which allow for unbounded model Hessians, complexity has focused on cases where the latter remain uniformly bounded \citep{cartis-gould-toint-2022,grapiglia-2014}, except for the work of \citet{leconte-orban-2023}, who study a trust-region method for nonsmooth optimization that does not reduce to a known classical method in the case of smooth optimization.

In practice, quasi-Newton updates, including PSB, BFGS, and SR1 are not guaranteed to result in uniformly bounded model Hessians.
Our main contributions are to address the complexity of a family of trust-region methods when the model Hessians may be unbounded, and to expand upon the results of~\citet{leconte-orban-2023}.
Our family of methods includes most, if not all, known classical variants, and is parameterized by two scalars \(-1 \leq \alpha \leq 1\) and \(0 \leq \beta \leq 1\) that control the presence of the gradient norm and the model Hessian norm in the trust-region radius.

Several authors show that usual quasi-Newton updates, including PSB, BFGS and SR1, grow at most linearly with \(k\), e.g., \citep{powell-2010} and \citep[\S\(8.4.1.2\)]{conn-gould-toint-2000}, although at the time of this writing, it is not known whether that bound is attained.
However, in a typical implementation, those approximations are not updated when a step is rejected, only when a step is accepted.
Thus, we also investigate the incidence on complexity of a weaker, but more realistic, assumption where the model Hessians grow at most linearly with the number of iterations in which a step is accepted, called \emph{successful} iterations.

Specifically, we establish that if \(\|B_k\| \leq \mu (1 + |\mathcal{S}_{k-1}|^p)\), where \(0 \leq p < 1\) and \(|\mathcal{S}_{k-1}|\) is the number of successful iterations up to iteration \(k-1\), our family of trust-region methods may require as many as \(O([(1-p) \epsilon^{-2}]^{1/(1 - p)})\) iterations to identify an \(\epsilon\)-stationary point.
Although \citet{leconte-orban-2023} presumably intended to use \(\mathcal{S}_{k-1}\) in their assumption, they formulated it in terms of \(\mathcal{S}_{k}\).
Our \Cref{asm:unb} makes the correction.
Even though we obtain a complexity bound of the same order as that of \citet{leconte-orban-2023}, our assumption combined with a finer analysis results in a more relevant constant hidden inside the \(O\) term: as \(p\) approaches \(1\), the analysis of \citet{leconte-orban-2023} offers no information, whereas we establish the bound \(O(\exp(c_1 \epsilon^{-2}))\) for some constant \(c_1 > 0\).
In both cases, we show that the bounds are sharp.
Although it is not currently known whether quasi-Newton updates may indeed grow exactly linearly with \(|\mathcal{S}_{k-1}|\), our analysis applies to them.

Attention then turns to the weaker assumption \(\|B_k\| \leq \mu (1 + k^p)\), where again, \(0 \leq p < 1\), which, among others, takes into account scenarios where \(B_k\) is also updated on unsuccessful iterations.
Here too, we establish a complexity bound \(O([(1-p) \epsilon^{-2}]^{1/(1 - p)})\).
However, when \(p = 1\), the complexity bound becomes \(O([1 + (1 - \alpha ) \log(\epsilon^{-1})] \exp(c_2 \epsilon^{-2}))\), for some constant \(c_2 > 0\).
The deterioration of the bound compared with the one obtained under the previous assumption is due to the update of \(B_k\) on unsuccessful iterations.
We demonstrate that our bound is sharp by constructing an objective function on which the trust-region algorithm performs both unsuccessful and successful iterations to reach the desired tolerance.
This type of example appears to be new; all existing minimization problems that we know of that are used to prove sharpness of a complexity bound are, by construction, designed to perform successful iterations only.

\citet{powell-1975} first investigated convergence of trust-region methods under the assumption \(\|B_k\| \leq \mu (1 + \sum_{j = 0}^{k-1} \|s_j\|)\), where \(\mu > 0\) is a constant, and \(s_j\) is the step at iteration \(j\).
Among others, the PSB quasi-Newton update \citep{powell-1970} satisfies that condition \citep{fletcher-1972}.
The last paragraph of his paper mentions another researcher investigating the weaker condition \(\|B_k\| \leq \mu (1 + k)\), and Powell conjectured that convergence continues to hold under the weaker assumption.
\citet{powell-1984} himself extended the convergence results of \citep{powell-1975} under the weaker assumption, and \citet{powell-2010} further extended the previous analysis to a wide family of trust-region methods.
In his conclusions, \citet{powell-2010} shared intuition on the number of iterations required to reduce the gradient norm, and described that number as ``monstrous''.
In \citep[\S 4]{powell-2010}, he hints at a complexity of order \(O(2^{j_0 \exp(c_0 \epsilon^{-2})})\) for some \(c_0 > 0\) and an integer \(j_0 > 0\), which is worse than the bound we establish.

Our complexity analysis also addresses the special cases of convex and strongly convex objectives in both scenarios of model Hessian growth.
We establish tighter bounds than those obtained in the nonconvex case.
Specifically, in order to reach the neighborhood \(f(x_k) - f_{\textup{low}} \leq \epsilon\) of a global minimum \(f_{\textup{low}}\) of convex and strongly convex functions under the assumption that \(\|B_k\| = O(|\mathcal{S}_{k-1}|^p)\), our family of trust-region algorithms requires at most \(O([(1-p) \epsilon^{-1}]^{1/(1 - p)})\) and \(O([(1-p)\log(\epsilon^{-1})]^{1/(1 - p)})\) function evaluations, respectively, when \(p < 1\).
When \(p = 1\), the bounds deteriorate to \(O(\exp(c_3 \epsilon^{-1}))\) and \(O(\epsilon^{-c_5})\) for convex and strongly convex functions, respectively, where \(c_3 > 0\) and \(c_5 > 0\) are constants.
If \(\|B_k\| = O(k)\), the previous bounds are multiplied by a factor of \((1 + (1 - \alpha ) \log(\epsilon^{-1}))\), as in the general case.
Finally, if the model Hessians are bounded, our bounds are similar to those established in \citep{cartis-gould-toint-2012,nesterov-2013,grapiglia-yuan-2016} for (strongly) convex functions.

\subsection{Contributions}

This research makes the following novel contributions:
\begin{itemize}
  \item A family of trust-region methods that includes most known methods as special cases.
    The family is parameterized by two scalars \(-1 \leq \alpha \leq 1\) and \(0 \leq \beta \leq 1\) that control the presence of the gradient norm and the model Hessian norm in the trust-region radius.
  \item A complexity analysis for this family when the model Hessians growth is bounded by the number of (successful) iterations.
    \citet{leconte-orban-2023} only consider the case where the bound on the model Hessians depends on the number of successful iterations.
    However, their analysis cannot be extended to the case where the bound depends on the total number of iterations (successful or not).
    Our analysis also covers that case.
  \item An exponential complexity bound for the family when the model Hessians grow linearly with the number of (successful) iterations.
    Although the resulting complexity bounds are, to borrow Powell's qualifier, ``monstrous'', they are shown to be sharp and significantly better than those conjectured by \citet{powell-2010}, who suggests a double-exponential complexity.
  \item A novel type of example to establish sharpness that requires the occurrence of unsuccessful iterations.
  \item A complexity analysis for the family applied to convex and strongly convex functions, covering cases with unbounded model Hessians.
\end{itemize}
\Cref{tab:complexity} summarizes the dominant terms of the complexity bounds depending on the convexity of \(f\) and the growth of the model Hessian \(\|B_k\|\).

\begin{table}[ht]%
  \centering
  \begin{tabular}{|l|l|l|l|}%
    \hline
    \(f\) \, \textbackslash \, \(\|B_k\|\) & \(O(|\mathcal{S}_{k-1}|^p)\) or \(O(k^p)\)
                                           & \(O(|\mathcal{S}_{k-1}|)\)
                                           & \(O(k)\)                                                       \\
    \hline
    \multirow{2}{*}{General}               & \([(1-p) \epsilon^{-2}]^{1/(1 - p)}\)
                                           & \(\exp(c_1 \epsilon^{-2})\)
                                           & \([1 + \log(\epsilon^{\alpha - 1})] \exp(c_2 \epsilon^{-2})\)  \\
                                           & \Cref{thm:complexity:S,thm:complexity-k}
                                           & \Cref{thm:complexity:S}
                                           & \Cref{thm:complexity-k}                                        \\
    \hline
    \multirow{2}{*}{Convex}
                                           & \([(1-p) \epsilon^{-1}]^{1/(1 - p)}\)
                                           & \(\exp(c_3 \epsilon^{-1})\)
                                           & \([1 + \log(\epsilon^{\alpha - 1})] \exp(c_4 \epsilon^{-1})\)  \\
                                           & \Cref{thm:complexity-convex:S,thm:complexity-convex:k}
                                           & \Cref{thm:complexity-convex:S}
                                           & \Cref{thm:complexity-convex:k}                                 \\
    \hline
    Strongly
                                           & \([(1-p)\log(\epsilon^{-1})]^{1/(1 - p)}\)
                                           & \(\epsilon^{-c_5}\)
                                           & \([1 + \log(\epsilon^{\alpha - 1})] \, \epsilon^{-c_6}\)       \\
    Convex                                 & \Cref{thm:complexity-convex-str:S,thm:complexity-convex-str:k}
                                           & \Cref{thm:complexity-convex-str:S}
                                           & \Cref{thm:complexity-convex-str:k}                             \\
    \hline
  \end{tabular}
  \caption{Summary of the dominant terms of the complexity bounds obtained under different assumptions on \(f\)  and the growth of the model Hessian \(\|B_k\|\).
    The bounds involving \(p\) are for \(0 \leq p < 1\), and the constants \(c_1, \ldots, c_6 > 0\).}%
  \label{tab:complexity}
\end{table}
\subsection{Organization of the paper}

The remainder of this paper is organized as follows.
In \Cref{sec:algorithm}, we introduce our family of trust-region methods along with key preliminary results.
In \Cref{sec:S,sec:k}, we derive complexity bounds for cases where the bound on model Hessians depends on successful iterations and the iteration counter, respectively.
In both sections, we establish that the bounds are sharp and derive specialized complexity bounds for convex and strongly convex functions.
In \Cref{sec:num}, we investigate the performance in practice of different members of our family of trust-region methods, including classical schemes.
Finally, we provide closing remarks in \Cref{sec:discussion}.

\subsection*{Notation}

For a finite set $\mathcal{A}$, $|\mathcal{A}|$ denotes its cardinality.
For a vector $x$ and matrix $B$, $\|x\|$ and $\|B\|$ denote their $\ell_2$ norm on $\R^n$ and the induced matrix norm, respectively.
We denote \(\N_0\) the set of positive integers and \(\R_+\) the set of nonnegative real numbers.
By convention, we set \(0^0 = 1\).

\section{A family of trust-region methods}%
\label{sec:algorithm}

In this section, we describe a parameterized family of trust-region methods for~\eqref{eq:nlp}, where the main novelty is the definition of a general trust-region radius in the subproblem.

\subsection{Algorithm}%

At each iteration \(k\), we compute a step \(s_k\) from current iterate \(x_k\) as an (inexact) solution of the subproblem
\begin{align}%
  \label{subprob:TR}
  \min_{s \in \R^{n}} m_k(s) \quad \text{s.t.} \ \|s\| \leq \frac{\|\nabla f(x_k)\|^\alpha}{(1 + \|B_k\|)^\beta} \Delta_k, \\
  m_k(s) \defeq f(x_k) + \nabla f(x_k)^{T} s + \tfrac{1}{2} s^{T} B_k s, \nonumber
\end{align}
where \(B_k = B_k^T\) and \(\Delta_k > 0\) is used to determine the trust-region radius.
Specific choices of \(\alpha \leq 1\) and \(\beta \leq 1\) reduce~\eqref{subprob:TR} to known formulations.
Namely, \(\alpha = \beta = 0\) leads to the classical method \citep{conn-gould-toint-2000}, while \(\alpha = 1\) and \(\beta = 0\) reduces to the choice of the trust-region radius used by \citep{fan-yuan-2001,curtis-lubberts-robinson-2018}---note however, that the update of \(\Delta_k\) differs in both \citep{fan-yuan-2001,curtis-lubberts-robinson-2018}.
When $\beta = 0$, our scaled trust-region radius can be seen as a particular case of more general nonlinear step-size control mechanisms \citep{toint-2011,grapiglia-2014}.
Choosing \(\beta \neq 0\) incorporates second-order information into the trust-region radius.
We are not aware of other trust-region methods doing so, except for \citep{leconte-orban-2023} in nonsmooth optimization.
The parameters $\alpha$ and $\beta$ are for generalization purposes, but also play a key role in our complexity analysis; a discussion on the values of these parameters, ensuring that the trust-region method enjoys favorable worst-case complexity bounds, is provided later.

The rest of the algorithm is standard.
Once a trial step $s_k$ has been determined, the decrease in \(f\) at $x_k + s_k$ is compared to the decrease predicted by the model.
If both are in sufficient agreement, $x_k + s_k$ becomes the new iterate, and \(\Delta_k\) is possibly increased.
If the model turns out to predict poorly the actual decrease, the trial point is rejected and \(\Delta_k\) is reduced.
\Cref{alg:TR-GEN} states the whole procedure.

\begin{algorithm}[ht]%
  \caption[caption]{%
    \label{alg:TR-GEN}
    A family of trust-region methods
  }
  \begin{algorithmic}[1]%
    \State Choose constants \(0 < \eta_1  \leq \eta_2 < 1\), \(0 < \gamma_1  \leq \gamma_2 < 1 \leq \gamma_3 \leq \gamma_4\) and \(0 < \kdc \le \frac{1}{2}\).
    \State Choose \(x_0 \in \R^n\) and \(\Delta_0 > 0\).
    \State Choose \(-1 \leq \alpha \leq 1\) and \(0 \leq \beta \leq 1\).
    \For{\(k = 0, 1, \dots\)}
    \State%
    Choose \(B_k = B_k^T  \in \R^{n \times n}\).
    \State%
    Compute \(\skcp\) as in~\eqref{eq:skcp} and an approximate minimizer \(s_k\) of~\eqref{subprob:TR} satisfying~\eqref{eq:suff-decrease}.
    \State%
    Compute the ratio
    \[
      \rho_k \defeq
      \frac{
        f(x_k) - f(x_k + s_k)
      }{
        m_k(0) - m_k(s_k)
      }.
    \]
    \State%
    If \(\rho_k \geq \eta_1\), set \(x_{k+1} \defeq x_k + s_k\).
    Otherwise, set \(x_{k+1} \defeq x_k\).
    \State%
    \label{alg:TR-GEN:step-update}%
    Update \(\Delta_k\) according to
    \[
      \Delta_{k+1} \in
      \begin{cases}
        \begin{aligned}
           & [\gamma_3 \Delta_k, \, \gamma_4 \Delta_k] &  & \text{ if } \rho_k \geq \eta_2,           &  & \quad \Comment{(very successful iteration)} \\
           & [\gamma_2 \Delta_k, \,  \Delta_k]         &  & \text{ if } \eta_1  \leq \rho_k < \eta_2, &  & \quad \Comment{(successful iteration)}      \\
           & [\gamma_1 \Delta_k, \, \gamma_2 \Delta_k] &  & \text{ if } \rho_k < \eta_1 .             &  & \quad \Comment{(unsuccessful iteration)}
        \end{aligned}
      \end{cases}
    \]
    \EndFor
  \end{algorithmic}
\end{algorithm}

Only an approximate solution of~\eqref{subprob:TR} is required; a step \(s_k\) should provide a decrease larger than or equal to a fraction of the decrease at the Cauchy point within the trust-region.
The Cauchy point \(\skcp\) is defined as
\begin{equation}%
  \label{eq:skcp}
  \skcp \defeq -t_{k} \nabla f(x_k),
\end{equation}
where,
\begin{equation*}%
  t_{k} \defeq \argmin{t \geq 0} \; m_{k}\left(-t \nabla f(x_k)\right) \text { s.t. }\left\|t \nabla f(x_k)\right\| \leq \frac{\|\nabla f(x_k)\|^\alpha}{(1 + \|B_k\|)^\beta}\Delta_k.
\end{equation*}
The decrease required of \(s_k\) is stated as
\begin{equation}%
  \label{eq:suff-decrease}
  m_k(0) - m_k(s_k) \ge 2\kappa_{\textup{mdc}} (m_k(0) - m_k(\skcp)),  
\end{equation}
where \(0 < \kdc \le \frac{1}{2}\).

We make the following assumption.

\begin{problemassumption}%
  \label{asm:problem-obj}
  The objective function \(f\) is continuously differentiable and is bounded below by a constant \(f_{\textup{low}}\).
\end{problemassumption}

Moreover, we replace the standard Lipschitz continuity assumption on the gradient with the following weaker assumption.

\begin{modelassumption}%
  \label{asm:model-error}
  There is \(L \ge 1\) such that, for all \(k \in \N\),
  \begin{equation*}%
    |f(x_k + s_k) - m_k(s_k)| \le \tfrac{1}{2} (L + \|B_k\|) \|s_k\|^2.
  \end{equation*}
\end{modelassumption}

This assumption is satisfied when \(\nabla f\) is Lipschitz continuous with constant \(L_f\), provided we set \(L = \max\{1, L_f\}\) \citep[Theorem 8.4.2]{conn-gould-toint-2000}.

Below, we allow \(\{B_k\}\) to be unbounded.
We study the effect of two different bounds on \(\{\|B_k\|\}\) on the worst-case complexity.
The next section provides results useful to both analyses.

\subsection{Preliminary results}

Our first result clarifies the Cauchy decrease~\eqref{eq:suff-decrease}.
The proof is similar to that of, e.g., \citep[Lemma 2.3.2]{cartis-gould-toint-2022} with the radius replaced with that of~\eqref{subprob:TR}.

\begin{lemma}%
  \label{lem:sufficient-cauchy-decrease}
  Let \Cref{asm:problem-obj,asm:model-error} hold.
  For all \(k \in \N\),
  \begin{equation*}
    m_{k}(0)-m_{k}\left(s_{k}\right) \geq  \kappa_{\textup{mdc}} \|\nabla f(x_k)\| \min \left\{\frac{\|\nabla f(x_k)\|}{1+\|B_k\|}, \, \frac{\|\nabla f(x_k)\|^\alpha}{(1 + \|B_k\|)^\beta} \Delta_k \right\}.
  \end{equation*}
\end{lemma}

The quantity
\begin{equation*}
  a_k \defeq \frac{\Delta_k \max_{j = 0, \ldots,k} [(L + \|B_j\|) (1 + \|B_j\|)^{- \beta}]}{\min_{j = 0, \ldots,k}\|\nabla f(x_j)\|^{1-\alpha}}
  \quad
  (k \in \N)
\end{equation*}
plays a central role in the analysis.
We first establish that \(\{a_k\}\) is uniformly bounded below.
The proof proceeds by induction as in \citep[Lemma 2.3.4]{cartis-gould-toint-2022}.

\begin{lemma}%
  \label{lem:gamma-min}
  Let \Cref{asm:problem-obj,asm:model-error} hold.
  For all \(k \in \N\),
  \begin{equation*}
    a_k \ge a_{\min} \defeq 2\gamma_1 \kdc (1 - \eta_2) \min\left\{ 1, \, \frac{\Delta_0}{\|\nabla f(x_0)\|^{1-\alpha}} \right\} > 0.
  \end{equation*}
  Equivalently, for all \(k \in \N\),
  \begin{equation*}
    \Delta_k \ge \frac{\min_{j = 0, \ldots,k}\|\nabla f(x_j)\|^{1-\alpha}}{\max_{j = 0, \ldots,k} [(L + \|B_j\|) (1 + \|B_j\|)^{- \beta}]} \ a_{\min}.
  \end{equation*}
\end{lemma}

\begin{proof}
  First, observe that for \(k = 0\),
  \begin{align*}
    a_0 = \frac{\Delta_0 (L + \|B_0\|) (1 + \|B_0\|)^{- \beta}}{\|\nabla f(x_0)\|^{1-\alpha}} & \ge \frac{\Delta_0  (1 + \|B_0\|)^{1 - \beta}}{\|\nabla f(x_0)\|^{1-\alpha}} \\
                                                                                              & \ge \frac{\Delta_0 }{\|\nabla f(x_0)\|^{1-\alpha}}
    \ge a_{\min},
  \end{align*}
  because \(L \ge 1\), \(1 - \beta \geq 0\) and \(0 < 2\gamma_1 \kdc (1 - \eta_2) \leq 1\).

  Assume that \(k \ge 0\) is the first iteration such that \(a_k \ge a_{\min}\) and \(a_{k+1} < a_{\min}\) in order to reach a contradiction.
  By \Cref{alg:TR-GEN:step-update}, \(\gamma_1 \Delta_{k} \le \Delta_{k+1}\), which implies that,
  \begin{align}
    \frac{\Delta_{k} \max_{j = 0, \ldots,k+1} (L + \|B_j\|) (1 + \|B_j\|)^{- \beta}}{\min_{j = 0, \ldots,k+1}\|\nabla f(x_j)\|^{1-\alpha}} & \le \frac{a_{k+1}}{\gamma_1} \nonumber \\
                                                                                                                                           & < \frac{ a_{\min}}{\gamma_1}
    \label{eq:contradiction}
    \le 2 \kdc (1 - \eta_2)\le 1.
  \end{align}

  On the other hand, it follows from the definition of \(\rho_k\), together with~\eqref{subprob:TR}, \Cref{lem:sufficient-cauchy-decrease}, \Cref{asm:model-error} and \(s_k\) satisfying~\eqref{eq:suff-decrease} that
  \begin{align*}
    |\rho_k -1| & = \frac{|f(x_k + s_k) - m_k(s_k)|}{m_k(0) - m_k(s_k)}                                                                                                                                                                                                                                             \\
                & \le \frac{ (L + \|B_k\|) \|s_k\|^2}{2\kappa_{\textup{mdc}} \|\nabla f(x_k)\| \min \left\{\frac{\|\nabla f(x_k)\|}{1+\|B_k\|}, \, \frac{\|\nabla f(x_k)\|^\alpha}{(1 + \|B_k\|)^\beta} \Delta_k \right\}}                                                                                          \\
                & \le \frac{ {(L + \|B_k\|) (1 + \|B_k\|)^{-2\beta}} \|\nabla f(x_k)\|^{2\alpha} \Delta_k^2}{2\kappa_{\textup{mdc}} \|\nabla f(x_k)\| \min \left\{\frac{\|\nabla f(x_k)\|}{1+\|B_k\|}, \, \frac{\|\nabla f(x_k)\|^\alpha}{(1 + \|B_k\|)^\beta} \Delta_k \right\}}                                   \\
                & = \frac{ {(L + \|B_k\|)(1 + \|B_k\|)^{-\beta}} \|\nabla f(x_k)\|^{\alpha-1} \Delta_k^2}{2\kappa_{\textup{mdc}} \frac{(1 + \|B_k\|)^{\beta}}{\|\nabla f(x_k)\|^\alpha} \min \left\{\frac{\|\nabla f(x_k)\|}{1+\|B_k\|}, \, \frac{\|\nabla f(x_k)\|^\alpha}{(1 + \|B_k\|)^\beta} \Delta_k \right\}} \\
                & =  \frac{ {(L + \|B_k\|)(1 + \|B_k\|)^{-\beta}} \|\nabla f(x_k)\|^{\alpha-1} \Delta_k^2}{2\kappa_{\textup{mdc}} \min \left\{\frac{\|\nabla f(x_k)\|^{1- \alpha}}{(1+\|B_k\|)^{1-\beta}}, \, \Delta_k \right\}}                                                                                    \\
                & \le \frac{ {(L + \|B_k\|)(1 + \|B_k\|)^{-\beta}} \|\nabla f(x_k)\|^{\alpha-1} \Delta_k^2}{2\kappa_{\textup{mdc}} \min \left\{\frac{\|\nabla f(x_k)\|^{1- \alpha}}{(L + \|B_k\|)(1 + \|B_k\|)^{-\beta}}, \, \Delta_k \right\}}                                                                     \\
                & = \frac{ {(L + \|B_k\|)(1 + \|B_k\|)^{-\beta}} \|\nabla f(x_k)\|^{\alpha-1} \Delta_k^2}{2\kappa_{\textup{mdc}} \Delta_k } \le (1 - \eta_2),
  \end{align*}
  where the final step follows from~\eqref{eq:contradiction}.
  Therefore, \(\rho_k \ge \eta_2\), which implies that iteration \(k\) is very successful, and \(\Delta_{k+1} \ge \Delta_k\).
  Thus,
  \begin{equation*}
    \frac{\Delta_{k} \displaystyle{\max_{0 \le j \le k}} (L + \|B_j\|) (1 + \|B_j\|)^{- \beta}}{\displaystyle{\min_{j = 0, \ldots,k}}\|\nabla f(x_j)\|^{1-\alpha}} \le \frac{\Delta_{k + 1} \displaystyle{\max_{0 \le j \le k+1}} (L + \|B_j\|) (1 + \|B_j\|)^{- \beta}}{\displaystyle{\min_{j = 0, \ldots,k+1}}\|\nabla f(x_j)\|^{1-\alpha}}.
  \end{equation*}
  The latter is strictly less than \(a_{\min}\), which contradicts our assumption.
  \qed
\end{proof}

Note that, if \(\Delta_0\) is selected large enough that
\begin{equation*}
  \Delta_{0} \geq \max(1, \, \|\nabla f(x_0)\|^2),
\end{equation*}
then, \(\Delta_{0} \geq \|\nabla f(x_0)\|^{1-\alpha}\) for all \(\alpha \in [-1, 1]\) and \(a_{\min}\) is independent of \(\alpha\) in \Cref{lem:gamma-min}.

We now state a technical lemma that simplifies the bound on the denominator of \(\Delta_k\) in \Cref{lem:gamma-min} depending on the values of \(\beta\).

\begin{lemma}%
  \label{lem:technical}
  Assume that \Cref{asm:problem-obj,asm:model-error} hold.
  Let \( k \in \N \), and consider any \( j = 0, \ldots, k \).
  Then,
  \begin{equation}
    \label{eq:maj-L-beta-k}
    (L + \|B_j\|)(1 + \|B_j\|)^{-\beta}%
    \le L (1 + \max_{j = 0, \ldots,k}\|B_j\|)^{1-\beta}.
  \end{equation}
  Moreover, if \(\beta \in [0, 1/L)\), then
  \begin{equation}
    \label{eq:maj-L-beta-k-2}
    (L + \|B_j\|)(1 + \|B_j\|)^{-\beta}%
    \le (L + \max_{j = 0, \ldots,k}\|B_j\|)(1 + \max_{j = 0, \ldots,k}\|B_j\|)^{-\beta}.
  \end{equation}
\end{lemma}

\begin{proof}
  Let \(k \in \N \).
  We have that, for all \(j = 0, \ldots, k\),
  \begin{align*}
    (L + \|B_j\|)(1 + \|B_j\|)^{-\beta} & = L \left(1 + \frac{\|B_j\|}{L}\right) (1 + \|B_j\|)^{- \beta} \\
                                        & \le L(1 + \|B_j\|)^{1-\beta}                                   \\
                                        & \le L (1 + \max_{j = 0, \ldots,k}\|B_j\|)^{1-\beta},
  \end{align*}
  since \( L \geq 1 \), which proves~\eqref{eq:maj-L-beta-k}.

  To prove~\eqref{eq:maj-L-beta-k-2}, define \(\varphi_{\beta}: \R_+ \rightarrow \R \) as \(\varphi_\beta(x) \defeq (L + x) (1 + x)^{- \beta}\).
  For \(x \geq 0\), the derivative of \(\varphi_\beta \) is
  \[
    \varphi'_\beta(x) = \frac{1 - \beta L + (1 - \beta) x}{(1 + x)^{1 + \beta}}.
  \]
  Thus, for any \(\beta \in [0, 1/L]\), \(\varphi_\beta \) is nondecreasing over \(\R_+\), which implies that \(\varphi_\beta(\|B_j\|) \le \varphi_\beta(\max_{j = 0, \ldots,k} \|B_j\|)\) and thus~\eqref{eq:maj-L-beta-k-2} holds.
  \qed
\end{proof}

\Cref{lem:technical} allows us to simplify the bound on \(\Delta_k\) in \Cref{lem:gamma-min} depending on the values of \(\beta \).
The following lemma provides a lower bound on \(\Delta_k\) that is useful in the complexity analysis of \Cref{alg:TR-GEN}.
The proof is a direct consequence of \Cref{lem:gamma-min,lem:technical}.

\begin{lemma}%
  \label{lem:Delta-min}
  Assume that \Cref{asm:problem-obj,asm:model-error} hold.
  Let \( k \in \N \), then,
  \begin{equation*}
    \Delta_k \ge \frac{\min_{j = 0, \ldots,k}\|\nabla f(x_j)\|^{1-\alpha}}{L_{\beta, k}} \ a_{\min},
  \end{equation*}
  where
  \begin{equation}
    \label{eq:L-beta-k}
    L_{\beta,k} \defeq
    \begin{cases}
      \begin{aligned}
         & (L + \max_{j = 0, \ldots,k}\|B_j\|)(1 + \max_{j = 0, \ldots,k}\|B_j\|)^{-\beta} &  & \text{ if } \beta \in [0, 1/L], \\
         & L (1 + \max_{j = 0, \ldots,k}\|B_j\|)^{1-\beta}                                 &  & \text{ if } \beta \in (1/L, 1].
      \end{aligned}
    \end{cases}
  \end{equation}
\end{lemma}

Note that the constant \(L_{\beta, k}\) can be simply set to \(L (1 + \max_{j = 0, \ldots, k} \|B_j\|)^{1-\beta}\).
However, setting \(L_{\beta, k}\) to \((L + \max_{j = 0, \ldots, k} \|B_j\|)(1 + \max_{j = 0, \ldots, k} \|B_j\|)^{-\beta}\) when \(\beta \in [0,1/L]\) yields a better complexity bound.
This choice also allows deriving complexity bounds that match existing ones, for instance, when \(\beta=0\)~\cite{cartis-gould-toint-2022}.

In the following, we repeatedly use
\begin{equation}
  \label{eq:L-beta}
  \psib{x} \defeq
  \begin{cases}
    \begin{aligned}
       & L + x     &  & \text{ if } \beta \in [0, 1/L], \\
       & L (1 + x) &  & \text{ if } \beta \in (1/L, 1].
    \end{aligned}
  \end{cases}
\end{equation}

The function \(\fpsib\) is related to \(L_{\beta,k}\) defined in~\eqref{eq:L-beta} via
\begin{equation}
  \label{eq:L-beta-psi}
  L_{\beta,k} = (1 + \max_{j = 0, \ldots,k}\|B_j\|)^{-\beta} \psib{\max_{j = 0, \ldots,k} \|B_j\|}.
\end{equation}
\Cref{lem:sufficient-cauchy-decrease,lem:Delta-min} allow us to quantify the model decrease in terms of \(a_{\min}\).

\begin{lemma}%
  \label{lem:sufficient-decrease}
  Let \Cref{asm:problem-obj,asm:model-error} hold.
  For all \(k \in \N\), \Cref{alg:TR-GEN} generates \(s_k\) such that
  \[
    m_{k}(0)-m_{k}\left(s_{k}\right) \geq   \kdc \frac{\|\nabla f(x_k)\|^{1 + \alpha} \min_{j = 0, \ldots,k} \|\nabla f(x_j)\|^{1 - \alpha}}{\psib{\max_{j = 0, \ldots,k} \|B_j\|}} a_{\min}.
  \]
\end{lemma}
\begin{proof}
  \Cref{lem:gamma-min,lem:sufficient-cauchy-decrease} and the fact that \(a_{\min} < 1\) combine to give
  \begin{align*}
     & m_{k}(0)-m_{k}\left(s_{k}\right)  \geq  \kappa_{\textup{mdc}} \|\nabla f(x_k)\| \min \left\{\frac{\|\nabla f(x_k)\|}{1+\|B_k\|}, \frac{\|\nabla f(x_k)\|^\alpha}{(1 + \|B_k\|)^\beta} \Delta_k \right\}                                                    \\
     & \quad \quad =  \kappa_{\textup{mdc}} \frac{\|\nabla f(x_k)\|^{1 + \alpha}}{(1 + \|B_k\|)^\beta} \min \left\{\frac{\|\nabla f(x_k)\|^{1 - \alpha}}{(1+\|B_k\|)^{1- \beta}}, \Delta_k \right\}                                                               \\
     & \quad \quad \geq  \kappa_{\textup{mdc}} \frac{\|\nabla f(x_k)\|^{1 + \alpha}}{(1 + \|B_k\|)^\beta} \min \left\{\frac{\|\nabla f(x_k)\|^{1 - \alpha}}{(L + \|B_k\|)(1+\|B_k\|)^{- \beta}}, \Delta_k \right\}                                                \\
     & \quad \quad \geq  \kappa_{\textup{mdc}} \frac{\|\nabla f(x_k)\|^{1 + \alpha}}{(1 + \|B_k\|)^\beta} \min \left\{\frac{\min_{j = 0, \ldots,k} \|\nabla f(x_j)\|^{1 - \alpha}}{\max_{j = 0, \ldots,k}[(L + \|B_j\|)(1+\|B_j\|)^{- \beta}]}, \Delta_k \right\} \\
     & \quad \quad =  \kappa_{\textup{mdc}} \frac{\|\nabla f(x_k)\|^{1 + \alpha} \min_{j = 0, \ldots,k} \|\nabla f(x_j)\|^{1 - \alpha}}{(1 + \|B_k\|)^\beta \max_{j = 0, \ldots,k}[(L + \|B_j\|)(1+\|B_j\|)^{- \beta}]} \min \left\{ 1, a_k \right\}.             \\
  \end{align*}
  By \Cref{lem:technical} we obtain \(\max_{j = 0, \ldots,k}[(L + \|B_j\|)(1+\|B_j\|)^{- \beta}] \le L_{\beta,k}\), where \(L_{\beta,k}\) is defined in~\eqref{eq:L-beta-k}.
  Thus,
  \begin{align*}
    m_{k}(0)-m_{k}\left(s_{k}\right) & \geq  \kappa_{\textup{mdc}} \frac{\|\nabla f(x_k)\|^{1 + \alpha} \min_{j = 0, \ldots,k} \|\nabla f(x_j)\|^{1 - \alpha}}{(1 + \max_{j = 0, \ldots,k} \|B_j\|)^\beta L_{\beta, k}} a_{\min} \\
                                     & =    \kappa_{\textup{mdc}} \frac{\|\nabla f(x_k)\|^{1 + \alpha} \min_{j = 0, \ldots,k} \|\nabla f(x_j)\|^{1 - \alpha}}{\psib{\max_{j = 0, \ldots,k} \|B_j\|}} a_{\min},
  \end{align*}
  where we used the relation~\eqref{eq:L-beta-psi} between \(L_{\beta,k}\) and \(\fpsib\).
  \qed
\end{proof}

We conclude this section with a technical result.
\begin{lemma}%
  \label{lem:sum-integral}
  Let $\mu > 0$ and $p > 0$.
  Then for any \(k_2 \in \N \),
  \begin{align*}
    \sum_{k = 0}^{ k_2} \frac{1}{\psib{\mu(1 + k^p)}} & \ge \frac{1}{\psib{\mu}} + \frac{1}{\psib{2\mu}}\int_{1}^{k_2+1} \frac{1}{t^p} \mathrm{d} t,
  \end{align*}
  and, for any \(k_1 \in \N\) and \(k_2 \in \N \) such that \(1 \le k_1 \le k_2\)
  \begin{align*}
    \sum_{k = k_1}^{ k_2} \frac{1}{\psib{\mu(1 + k^p)}} & \ge \frac{k_1^p}{\psib{\mu(1 + k_1^p)}}\int_{k_1}^{k_2+1} \frac{1}{t^p} \mathrm{d} t.
  \end{align*}
\end{lemma}
\begin{proof}%
  Define \(\phi_\beta: \R_+ \to \R\),
  \[
    \phi_\beta(x) \defeq x / \psib{\mu(1 + x)} =
    \begin{cases}
      \begin{aligned}
         & \frac{x}{L + \mu(1 + x)}     &  & \text{ if } \beta \in [0, 1/L], \\
         & \frac{x}{L (1 + \mu(1 + x))} &  & \text{ if } \beta \in (1/L, 1],
      \end{aligned}
    \end{cases}
  \]
  which is increasing.
  For \( k > 0 \),
  \begin{equation}
    \label{eq:integral}
    \int_k^{k+1} \frac{1}{t^p} \, \mathrm{d} t \leq \int_k^{k+1} \frac{1}{k^p} \, \mathrm{d} t = \frac{1}{k^p}.
  \end{equation}
  Then, for any \(k_2 \in \N \), we separate the first term in order to apply~\eqref{eq:integral},
  \begin{align*}
    \sum_{k = 0}^{ k_2} \frac{1}{\psib{\mu(1 + k^p)}} = \frac{1}{\psib{\mu}} + \sum_{k = 1}^{ k_2} \frac{\phi_\beta(k^p)}{k^p} & \geq \frac{1}{\psib{\mu}} + \phi_\beta(1) \sum_{k = 1}^{k_2} \frac{1}{k^p}                               \\
                                                                                                                               & \geq \frac{1}{\psib{\mu}} + \phi_\beta(1) \sum_{k = 1}^{k_2} \int_k^{k+1} \frac{1}{t^p} \, \mathrm{d} t.
  \end{align*}
  Otherwise, for any \(k_1 \in \N\) and \(k_2 \in \N \) such that \(1 \le k_1 \le k_2\),
  \begin{align*}
    \sum_{k = k_1}^{ k_2} \frac{1}{\psib{\mu(1 + k^p)}} = \sum_{k = k_1}^{ k_2} \frac{\phi_\beta(k^p)}{k^p} & \geq \phi_\beta(k_1^p) \sum_{k = k_1}^{k_2} \frac{1}{k^p}                               \\
                                                                                                            & \geq \phi_\beta(k_1^p) \sum_{k = k_1}^{k_2} \int_k^{k+1} \frac{1}{t^p} \, \mathrm{d} t.
    \tag*{\qed}
  \end{align*}
\end{proof}


In the next sections, we derive worst-case complexity analyses for \Cref{alg:TR-GEN} that allows for potentially unbounded model Hessians \(B_k\).
We will repeatedly use the notation
\begin{alignat*}{2}%
  \mathcal{S}   & \defeq \{ i \in \N \mid \rho_i \geq \eta_1 \}              &  & \qquad \text{(all successful iterations)}
  \\
  \mathcal{S}_k & \defeq \{ i \in \mathcal{S} \mid i \le k \}                &  & \qquad \text{(successful iterations until iteration \(k\))}
  \\
  \mathcal{U}_k & \defeq \{i \in \N \mid i\not \in \mathcal{S}, \ i \le k \} &  & \qquad \text{(unsuccessful iterations until iteration \(k\))}.
\end{alignat*}

\section{Complexity when the bound on \texorpdfstring{\boldmath{\(\|B_k\|\)}}{∥Bₖ∥} depends only on successful iterations}%
\label{sec:S}

Our assumption on the growth of the model Hessians is as follows.
\begin{modelassumption}%
  \label{asm:unb}
  There exist \(\mu > 0\) and \(0 \leq p \leq 1\) such that
  \begin{equation*}%
    \max_{j = 0, \ldots,k} \|B_j\| \le \mu (1 + |\mathcal{S}_{k-1}|^p) \qquad (k \in \N),
  \end{equation*}
  where \(|\mathcal{S}_{-1}| \defeq 0\) by convention.
\end{modelassumption}

Several comments on \Cref{asm:unb} are in order.

The special case \(p = 0\), which corresponds to bounded model Hessians, is the only one considered in the literature for complexity analyses (e.g., \citep{cartis-gould-toint-2022,curtis-lubberts-robinson-2018}).
\citet{leconte-orban-2023} provided the first complexity analysis allowing for potentially unbounded model Hessians.
They studied complexity when \(0 \le p < 1\) in the context of a trust-region method for nonsmooth optimization, although their results also apply to smooth optimization.
At iteration \(k\), \citet{leconte-orban-2023} impose a bound on \(\max_{j = 0, \ldots, k} \|B_j\|\) that depends on \(|\mathcal{S}_{k}|^p\), where \(0 \le p < 1\).

We improve upon the model Hessian bound assumption of \citet{leconte-orban-2023} in two ways.
First, in our model assumption, \(\max_{j = 0, \ldots, k} \|B_j\|\) depends on \(|\mathcal{S}_{k-1}|\) rather than \(|\mathcal{S}_{k}|\).
This choice is more natural since \(B_k\) depends on previous updates \(B_0, \dots, B_{k-1}\), and when updates occur only on successful iterations, \(B_k\) should be bounded by the number of successful iterations up to \(k-1\).
Second, unlike \citet{leconte-orban-2023}, \cref{asm:unb} also includes the case \(p = 1\).

We believe that the present work is the first to address the case where \(p = 1\); an important case as it encompasses known bounds on existing quasi-Newton approximations such as SR1.
Indeed, \citet[\S \(8.4.1.2\)]{conn-gould-toint-2000} show that the SR1 Hessian approximation satisfies
\begin{equation}%
  \label{eq:SR1}
  \|B_{k+1}\| \leq \|B_{k}\| + \mu_{f},
\end{equation}
where \(\mu_f > 0\) is a constant related to \(f\), and $k$ is the index of a successful iteration for which the SR1 update is well defined.
They provide a similar bound for the BFGS update when $f$ is convex.
\citet{powell-2010} establishes a similar bound for the PSB update, where the bound on \(\|B_{k+1}\|\) depends linearly on \(\|B_{k}\|\), as in~\eqref{eq:SR1}.
Even though it is not currently known whether those bounds on SR1, BFGS and PSB are tight, \cref{asm:unb} includes them when \(p = 1\).
Indeed by induction, and if we update the model Hessian only on successful iterations, we have that for all \(k \in \N\),
\begin{equation*}
  \|B_k\| \leq \|B_0\| + \mu _{f} |\mathcal{S}_{k - 1}| \leq \mu (1 + |\mathcal{S}_{k - 1}|),
\end{equation*}
with \(\mu = \max\{\|B_0\|, \mu_f\}\).
Note that, although uncommon, if we update the model Hessian \(B_k\) at each iteration \(k\), we obtain \( \|B_k\| = O(k) \).
This scenario will be considered in \Cref{sec:k}.

In practice, limited-memory quasi-Newton methods, such as LBFGS and LSR1, are commonly used.
For the LBFGS update, \citet[Lemma~\(3\)]{burdakov-gong-zirkin-yuan-2017} showed that if the sequence of initial matrices \(\{B_{k,0}\}\) at the current iteration is bounded, then \(\{B_k\}\) is also bounded provided that \(f\) is twice continuously differentiable with \(\nabla^2 f(x)\) uniformly bounded, and that the LBFGS update is used only when \(s_k^T y_k \ge \omega_1 \|s_k\| \|y_k\|\) for some \(0 < \omega_1 < 1\), where \(s_k \defeq x_{k+1} - x_k\) and \(y_k \defeq \nabla f(x_{k+1}) - \nabla f(x_k)\).
Regarding the LSR1 update, \citep[\S 4.2]{aravkin-baraldi-orban-2021} provides a procedure to maintain the boundedness of the extreme eigenvalues.
In fact, as long as the LSR1 updates are performed only for iterations such that \(\|s_k^T z_k\| \ge \omega_2 \|z_k\|^2\) for some \(\omega_2 > 0\), where \(z_k \defeq y_k - B_k s_k\), it follows that if \(\{B_{k,0}\}\) is bounded, then \(\{B_k\}\) is also bounded.
Hence, for both LBFGS and LSR1, the sequence of initial matrices \(B_{k,0}\) must remain bounded to ensure boundedness of \(B_k\).
As suggested in \citep{lu-1996}, when \(\nabla f\) is Lipschitz-continuous, the choice \(B_{k,0} \defeq \gamma_k I\) with \(\gamma_k \defeq s_k^T y_k / s_k^T s_k\) guarantees this boundedness.
Choosing \(\gamma_k \defeq y_k^T y_k / s_k^T y_k\) \citep{burdakov-gong-zirkin-yuan-2017,nocedal-wright-1999,kanzow-steck-2023,lu-1996} also ensures the boundedness when \(f\) is convex \citep{conn-gould-toint-2000}.
We note that, unlike LBFGS and LSR1, it is currently not known whether LPSB approximations \citep{kanzow-steck-2023} remain bounded.

In~\Cref{asm:unb}, \(p\) is not allowed to take a value larger than \(1\).
In fact, when \(p > 1\), global convergence is no longer guaranteed as the series \(\sum_{k=0}^{\infty} 1 / (1 + |\mathcal{S}_{k-1}|^p)\) may become convergent~\citep{toint-1988}.
Finally, due to the non-decreasing nature of \(\{|\mathcal{S}_{k-1}|\}_{k \in \mathbb{N}}\), \Cref{asm:unb} is equivalent to
\begin{equation*}
  \|B_k\| \le \mu(1 + |\mathcal{S}_{k - 1}|^p)
  \quad
  \textup{for all } k \in \N.
\end{equation*}

If \Cref{alg:TR-GEN} generates only a finite number of successful iterations, a first-order critical point is identified after a finite number of iterations.
The following result parallels \citep[Theorem~\(6.4.4\)]{conn-gould-toint-2000}.

\begin{theorem}%
  Let \Cref{asm:problem-obj,asm:model-error,asm:unb} be satisfied.
  If \Cref{alg:TR-GEN} generates finitely many successful iterations, then \(x_k = x^*\) for all sufficiently large \(k\) where \(\nabla f(x^*) = 0\).
\end{theorem}

\begin{proof}
  Assume by contradiction that there exists \(\nu > 0\) such that \(\|\nabla f(x_k)\| \geq \nu\) for all \(k \in \N\), and let \(k_f\) be the last successful iteration.
  Necessarily, \(x_k = x_{k_f}\) for all \(k \geq k_f\), and hence, \(\{x_k\} \to x^* \defeq x_{k_f}\).
  By \Cref{lem:gamma-min,lem:technical}, for \(k \geq k_f\),
  \begin{align*}
    \Delta_k & \geq \frac{(\min_{j = 0, \ldots, k} \|\nabla f(x_j)\|)^{1-\alpha} \, a_{\min}}{\max_{j = 0, \ldots,k} [(L + \|B_j\|) (1 + \|B_j\|)^{- \beta}] } \\
             & \geq \frac{\nu^{1-\alpha}\, a_{\min}}{ L(1 + \max_{j = 0, \ldots,k}\|B_j\|)^{1- \beta}}                                                         \\
             & \geq \frac{\nu^{1-\alpha} \, a_{\min}}{L(1 + \mu(1 + |\mathcal{S}_{k_f - 1}|^p))^{1- \beta}},
  \end{align*}
  where we used the fact that \(\mathcal{S}_{k - 1} = \mathcal{S}_{k_f - 1}\) for \(k \ge k_f\) and \(\alpha, \beta \le 1\).
  Thus, \(\Delta_k\) is bounded away from zero.
  However, the mechanism of \Cref{alg:TR-GEN} ensures that \(\Delta_k\) decreases on unsuccessful iterations and converges to \(0\), which is a contradiction, and shows that \(\liminf_{k \to \infty} \|\nabla f(x_k)\| = \|\nabla f(x^*)\| = 0\).
  \qed
\end{proof}

In the next section, we derive a bound on the number of successful iterations when \Cref{alg:TR-GEN} generates infinitely many successful iterations.
In \Cref{sec:convex-S}, we sharpen the analysis in the case where \(f\) is convex.

\subsection{Complexity on general objectives}%

Let \(\epsilon > 0\) and \(k_{\epsilon}\) be the first iteration of \Cref{alg:TR-GEN} such that \(\|\nabla f(x_{k_{\epsilon}})\| \le \epsilon\).
Define
\begin{align*}
  \mathcal{S}(\epsilon) & \defeq \mathcal{S}_{k_\epsilon-1} = \{ k \in \mathcal{S} \mid k < k_{\epsilon} \},
  \\
  \mathcal{U}(\epsilon) & \defeq \mathcal{U}_{k_\epsilon-1} = \{k \in \N \mid k \not \in \mathcal{S} \text{ and } k < k_{\epsilon} \}.
\end{align*}

The next theorem states a bound on the cardinality of $\mathcal{S}(\epsilon)$ when infinitely many successful iterations are generated.

\begin{theorem}%
  \label{thm:complexity:S}
  Let \Cref{asm:problem-obj,asm:model-error,asm:unb} hold.
  Assume that \Cref{alg:TR-GEN} generates infinitely many successful iterations and that \(f(x_k) \ge f_{\textup{low}}\) for all \(k \in \N \).
  Let  \(\mu \) and \(p\) be as in \Cref{asm:unb}, \(a_{\min}\) be as in \Cref{lem:gamma-min} and \(\fpsib \) is defined in~\eqref{eq:L-beta}.
  Define
  \begin{equation}
    \label{eq:kappa1-kappa2}
    \kappa_{0} \defeq\frac{\psib{2\mu}}{ \eta_1 \kdc a_{\min}} > 0,
    \mathhfill
    \kappa_{1} \defeq (f(x_0) - f_{\textup{low}}) \,  \kappa_0 \geq 0,
    \mathhfill
    \kappa_{2} \defeq\frac{\psib{2\mu}}{\psib{\mu}} > 0.
  \end{equation}
  Then, if \(p = 0\),
  \begin{equation}
    \label{eq:S-eps0}
    |\mathcal{S} (\epsilon)| \leq \kappa_1 \; \epsilon^{-2} = O(\epsilon^{-2}).
  \end{equation}
  If \(0 < p < 1\),
  \begin{equation}
    \label{eq:S-eps2}
    |\mathcal{S} (\epsilon)| \leq \left[ (1-p)\;  (\kappa_1 \; \epsilon^{-2}-\kappa_2) + 1\right]^{1/ (1-p)} = O\left([(1-p) \epsilon^{-2}]^{1/ (1-p)}\right).
  \end{equation}
  Otherwise, if \(p = 1\),
  \begin{equation}
    \label{eq:S-exp}
    |\mathcal{S} (\epsilon)| \leq \exp (\kappa_1\; \epsilon^{-2}-\kappa_2).
  \end{equation}
\end{theorem}

\begin{proof}
  If \(\mathcal{S}(\epsilon) = \varnothing \), then the theorem holds trivially as \(|\mathcal{S}(\epsilon)| = 0\).

  Otherwise, let \(\ell \in \mathcal{S}(\epsilon)\).
  By \Cref{lem:sufficient-decrease} and the fact that \(\alpha \in [-1,1]\),
  \begin{align*}
    f(x_{\ell}) - f(x_{\ell} + s_{\ell}) & \geq \eta_1 \left(m_{\ell}(0) - m_{\ell}(s_{\ell})\right)                                                                                                               \\
                                         & \geq \eta_1 \kdc a_{\min} \frac{ \|\nabla f(x_\ell)\|^{1 + \alpha} \min_{j = 0, \ldots, \ell} \|\nabla f(x_j)\|^{1 - \alpha}}{\psib{\max_{j = 0, \ldots,\ell} \|B_j\|}} \\
                                         & \geq \eta_1 \kdc a_{\min} \frac{ \epsilon^{1 + \alpha} \epsilon^{1 - \alpha}}{\psib{\max_{j = 0, \ldots,\ell} \|B_j\|}}.
  \end{align*}
  \Cref{asm:unb} and the increasing nature of \(\fpsib\) yield
  \begin{align*}
    f(x_{\ell}) - f(x_{\ell} + s_{\ell}) & \geq \eta_1 \kdc a_{\min}\epsilon^2 \frac{1}{\psib{\mu(1 + |\mathcal{S}_{\ell - 1}|^p)}}.
  \end{align*}
  We sum the above inequality over all \({\ell} \in \mathcal{S}(\epsilon)\), use a telescoping argument, and obtain
  \begin{align*}
    f(x_0) - f_{\textup{low}} & \geq \eta_1 \kdc a_{\min}\epsilon^2 \sum_{\ell \in \mathcal{S}(\epsilon) } \frac{1}{\psib{\mu(1 + |\mathcal{S}_{\ell - 1}|^p)}}           \\
                              & = \eta_1 \kdc a_{\min}\epsilon^2 \sum_{k = 0}^{ |\mathcal{S}(\epsilon)|-1} \frac{1}{\psib{\mu(1 + |\mathcal{S}_{\varphi(k) - 1}|^p)}}\; ,
  \end{align*}
  where \(\varphi\) is an increasing map from \(\{0,\ldots, |\mathcal{S}_\epsilon| -1\}\) to \(\mathcal{S}_\epsilon\).
  Hence, by definition of \(\varphi\) and \(\mathcal{S}_{\varphi(k)}\),
  \[
    \left|\mathcal{S}_{\varphi(k+1)}\right| =\left|\mathcal{S}_{\varphi(k)}\right| +1~~~\mbox{and}~~\left|\mathcal{S}_{\varphi(0)}\right| = 1.
  \]
  In other words, \(\left|\mathcal{S}_{\varphi(k)}\right| = k + 1\), and because \(\varphi(k - 1) \le \varphi(k) - 1\), we get \( k = \left|\mathcal{S}_{\varphi(k - 1)}\right| \le \left|\mathcal{S}_{\varphi(k) - 1}\right| < \left|\mathcal{S}_{\varphi(k)}\right| = k + 1\).
  Thus, \(|\mathcal{S}_{\varphi(k) - 1}| = k\), and
  \begin{equation}%
    \label{eq:sum-series}
    f(x_0) - {f}_{\textup{low}} \geq \eta_1 \kdc a_{\min}\epsilon^2 \sum_{k = 0}^{ |\mathcal{S}(\epsilon)|-1} \frac{1}{\psib{\mu(1 + k^p)}}.
  \end{equation}
  Hence, if \(p = 0\), we have
  \begin{equation*}
    \frac{f(x_0) - {f}_{\textup{low}}}{\eta_1 \kdc a_{\min}\epsilon^2} \geq \frac{1}{\psib{2\mu}} |\mathcal{S}(\epsilon)|,
  \end{equation*}
  which gives~\eqref{eq:S-eps0}.

  Otherwise, if \(0 < p \le 1\),~\eqref{eq:sum-series} and \Cref{lem:sum-integral} give
  \begin{equation*}
    \frac{f(x_0) - {f}_{\textup{low}}}{\eta_1 \kdc a_{\min}\epsilon^2} \geq \frac{1}{\psib{\mu}} + \frac{1}{\psib{2\mu}}\int_{1}^{|\mathcal{S}(\epsilon)|} \frac{1}{t^p} \mathrm{d} t.
  \end{equation*}
  Two cases follow.
  The first one, \(0 < p < 1\), gives
  \begin{equation*}
    \frac{f(x_0) - {f}_{\textup{low}}}{\eta_1 \kdc a_{\min}\epsilon^2} \geq   \frac{1}{\psib{\mu}} + \frac{1}{\psib{2\mu}} \left(\frac{{|\mathcal{S}(\epsilon)|}^{1-p} - 1}{1-p} \right),
  \end{equation*}
  which provides~\eqref{eq:S-eps2}.
  The second case, \(p = 1\), gives
  \begin{equation*}
    \frac{f(x_0) - {f}_{\textup{low}}}{\eta_1 \kdc a_{\min}\epsilon^2} \geq   \frac{1}{\psib{\mu}} + \frac{1}{\psib{2\mu}} \log\left(|\mathcal{S}(\epsilon)|  \right) ,
  \end{equation*}
  which establishes~\eqref{eq:S-exp}.
  \qed
\end{proof}

When \(0 \leq p < 1\) in \Cref{asm:unb}, \Cref{thm:complexity:S} improves the constant in the complexity bound of \citep[Lemma~\(2\)]{leconte-orban-2023}.
In particular, the bound given in \citep{leconte-orban-2023} diverges as $p \to 1$, whereas~\eqref{eq:S-exp} suggests that the complexity becomes exponential.

For \(p = 1\) in \Cref{asm:unb}, \Cref{thm:complexity:S} provides a complexity bound of order \( O(\exp(c_1 \epsilon^{-2})) \), where \(c_1 > 0\) is a constant.
This complexity bound is significantly better than Powell's intuition \citep[\S 4]{powell-2010}, who suggested a complexity of order \( O(2^{j_0 \exp(c_0 \epsilon^{-2})}) \) with \(c_0 > 0\) and \(j_0\) a positive integer.
Indeed, within the proof of the main theorem in \citep[\S 3]{powell-2010}, for some \(\epsilon > 0\) sufficiently small, the total number of iterations \(k_\epsilon\) must be less than \(2^{j_\epsilon + 1}\), where \(j_\epsilon\) satisfies \citep[Equation (3.33)]{powell-2010}, i.e.,
\begin{equation}%
  \label{eq:powell}
  f(x_0) - f_{\textup{low}} > b_1 \epsilon^2 \sum_{j=j_0}^{j_\epsilon} \frac{1}{b_2 \, j},
\end{equation}
where the constants \(b_1>0\), \(b_2>0\), and \(j_0\in \N\) are problem dependent.
Hence, by \Cref{lem:comparison-series},~\eqref{eq:powell} reduces to
\[
  f(x_0) - f_{\textup{low}} > \frac{b_1 \epsilon^{2}}{b_2 \, j_0} \int_{j_0}^{j_\epsilon+1} \frac{1}{t} \mathrm{d} t = \frac{b_1 \epsilon^{2}}{b_2 \, j_0} \left( \log(j_\epsilon+1) - \log(j_0) \right),
\]
which implies that
\[
  j_\epsilon + 1 < j_0 \exp\left(\frac{b_2 \, j_0 (f(x_0) - f_{\textup{low}})}{b_1 } \epsilon^{-2}\right).
\]
Consequently, we deduce that
\[
  k_\epsilon < 2^{j_\epsilon + 1} < 2^{j_0 \exp(c_0 \epsilon^{-2})},
\]
with \(c_0:= \frac{b_2 \, j_0 (f(x_0) - f_{\textup{low}})}{b_1 } > 0 \).

The next theorem is a classical result that bounds the number of unsuccessful iterations.

\begin{theorem}%
  \label{thm:complexity:U}
  Let \Cref{asm:problem-obj,asm:model-error,asm:unb} hold.
  Assume that \Cref{alg:TR-GEN} generates infinitely many successful iterations.
  Then,
  \begin{align*}
    |\mathcal{U}(\epsilon)| \leq |\log_{\gamma_2} (\gamma_4)| |\mathcal{S}(\epsilon)| & + (1-\alpha)\log_{\gamma_2}(\epsilon) -\log_{\gamma_2}(\psib{\mu(1 + |\mathcal{S}(\epsilon)|^p)})                       \\
                                                                                      & + \beta\log_{\gamma_2}(1 + \mu(1 + |\mathcal{S}(\epsilon)|^p)) + \log_{\gamma_2}\left(\frac{a_{\min}}{\Delta_0}\right),
  \end{align*}
  where \(L, \mu \) and \(p\) are as in \Cref{asm:unb,asm:model-error} and \(\fpsib \) is as in~\eqref{eq:L-beta}.
\end{theorem}
\begin{proof}
  As in \citep[Lemma~\(2.3.1\)]{cartis-gould-toint-2022}, the mechanism of \Cref{alg:TR-GEN} guarantees that for all \(k \in \N\),
  \[
    \left|\mathcal{U}_k\right| \leq |\log_{\gamma_2} (\gamma_4)| \left|\mathcal{S}_k\right| +  \log_{\gamma_2}\left(\frac{\displaystyle  \Delta_{k}}{\Delta_0}\right).
  \]
  Hence, \Cref{lem:gamma-min} yields
  \begin{align*}
    |\mathcal{U}(\epsilon)| & \leq |\log_{\gamma_2} (\gamma_4)| |\mathcal{S}(\epsilon)| + \log_{\gamma_2}\left(\frac{\displaystyle (\min_{j = 0, \ldots,k_\epsilon -1}\|\nabla f(x_j)\|)^{1-\alpha}a_{\min}}{L_{\beta,k} \;\Delta_0}\right) \\
                            & \leq |\log_{\gamma_2} (\gamma_4)| |\mathcal{S}(\epsilon)| + (1-\alpha)\log_{\gamma_2}(\epsilon) - \log_{\gamma_2}(L_{\beta,k}) + \log_{\gamma_2}\left(\frac{a_{\min}}{\Delta_0}\right)                        \\
                            & =    |\log_{\gamma_2} (\gamma_4)| |\mathcal{S}(\epsilon)| + (1-\alpha)\log_{\gamma_2}(\epsilon)     + \log_{\gamma_2}\left(\frac{a_{\min}}{\Delta_0}\right)                                                   \\
                            & \quad \quad \quad \quad - \log_{\gamma_2}\left(\psib{\max_{j = 0, \ldots,k} \|B_j\|} (1 + \max_{j = 0, \ldots,k} \|B_j\|)^{-\beta}\right)                                                                     \\
                            & \leq |\log_{\gamma_2} (\gamma_4)| |\mathcal{S}(\epsilon)| + (1-\alpha)\log_{\gamma_2}(\epsilon)                  + \log_{\gamma_2}\left(\frac{a_{\min}}{\Delta_0}\right)                                      \\
                            & \quad \quad \quad \quad - \log_{\gamma_2}\left(\psib{\mu(1 + |\mathcal{S}(\epsilon)|^p)} (1 + \mu (1 + |\mathcal{S}(\epsilon)|^p))^{-\beta}\right),
  \end{align*}
  where we recall \Cref{asm:unb} and the increasing nature of \(\fpsib\).
  \qed
\end{proof}

In the next theorem, we provide the optimal parameters \(\alpha\) and \(\beta\) that minimize the total number of iterations under mild assumptions on the initial trust-region radius.

\begin{theorem}%
  \label{thm:best-parameters}
  Under the assumptions of \Cref{thm:complexity:U}, assume also that \(\epsilon \in (0,1) \) and \(\Delta_0 \ge \max(1, \|\nabla f(x_0) \|^2)\).
  Then, setting \(\alpha = 1\) and \(\beta = 1/L\) in \Cref{alg:TR-GEN} ensures the best complexity bound.
\end{theorem}

\begin{proof}
  For \(|\mathcal{S}(\epsilon)|\), the leading term in the bound~\eqref{eq:S-eps2} is \(\kappa_1 \epsilon^{-2}\).

  First, note that \(\kappa_1\) can depend on \(\alpha \) only through \(a_{\min} \), see~\eqref{eq:kappa1-kappa2}.
  However, because \(\Delta_0 \ge \max(1, \|\nabla f(x_0) \|^2)\), for all \(\alpha \in [-1,1]\),
  \[
    \|\nabla f(x_0)\|^{1-\alpha} \leq \max(1, \|\nabla f(x_0) \|^2) \le \Delta_0.
  \]
  Hence,  \(a_{\min} = 2 \gamma_1 \kdc (1 - \eta_2)\) and thus \(\kappa_1\) is independent of \(\alpha \).

  Second, with respect to \(\beta \), \(\kappa_1\) depends on \(\beta \) only through \(\fpsib \), see~\eqref{eq:kappa1-kappa2}.
  Thus, as shown in~\eqref{eq:L-beta}, the constant \(\kappa_1\) is minimal for \(\beta \in [0,1/L]\).
  In conclusion, the complexity bound on \( |\mathcal{S}(\epsilon)| \) is optimal for any \( \alpha \in [-1,1] \) and \( \beta \in [0,1/L] \).

  Next, consider the bound \(|\mathcal{U}(\epsilon)|\) as given in \Cref{thm:complexity:U}, for any \( \alpha \in [-1,1] \) and \(\beta \in [0,1/L]\),
  \begin{align*}
    |\mathcal{U}(\epsilon)| & \leq |\log_{\gamma_2} (\gamma_4)| |\mathcal{S}(\epsilon)| + (1-\alpha)\log_{\gamma_2}(\epsilon) -\log_{\gamma_2}
    (L + \mu(1 + |\mathcal{S}(\epsilon)|^p))                                                                                                                                           \\
                            & \quad \quad + \beta \log_{\gamma_2}(1 + \mu(1 + |\mathcal{S}(\epsilon)|^p)) + \log_{\gamma_2}\left(\frac{2 \gamma_1 \kdc (1 - \eta_2)}{\Delta_0}\right).
  \end{align*}
  Since \(0 < \gamma_2 < 1\), the term \(\log_{\gamma_2}(1 + \mu(1 + |\mathcal{S}(\epsilon)|^p)) \le 0\) and hence, with respect to \(\beta \in [0,1/L]\), the term  \(\beta \log_{\gamma_2}(1 + \mu(1 + |\mathcal{S}(\epsilon)|^p))\) is minimized when \(\beta = 1/L\).
  With respect to \( \alpha \in [-1,1] \), it suffices to minimize the term \( (1-\alpha)\log_{\gamma_2}(\epsilon) \), which attains its minimum when \( \alpha = 1 \).
  \qed
\end{proof}

According to \Cref{thm:best-parameters}, the best iteration complexity bound is obtained by \(\alpha = 1\) and \(\beta = 1/L\) provided that \(\Delta_0 \ge \max(1, \|\nabla f(x_0) \|^2)\).
This setting is not practical, as the constant \(L\) is not available in practice—being generally related to the Lipschitz constant of \(\nabla f\).
However, an interesting alternative can be obtained by setting \(\alpha = 1\) and \(\beta = 0\).
In fact, in this case, the complexity bound in \Cref{thm:complexity:U} simplifies to
\begin{align*}
  |\mathcal{U}(\epsilon)| & \leq |\log_{\gamma_2} (\gamma_4)| |\mathcal{S}(\epsilon)| - \log_{\gamma_2}(L + \mu(1 + |\mathcal{S}(\epsilon)|^p)) \\
                          & \quad \quad + \log_{\gamma_2}\left(\frac{ 2 \gamma_1 \kdc (1 - \eta_2)}{\Delta_0}\right).
\end{align*}

In particular, when \(p=0\), i.e., the model Hessians in \Cref{asm:unb} are uniformly bounded by \(2 \mu \), the total number of iterations reduces to
\begin{align}
  \label{eq:complexity-4:alpha=1:beta=0:p=0}
  k_\epsilon & \leq \left(|\log_{\gamma_2} (\gamma_4)| + 1 \right) \frac{2(L + 2\mu )}{\gamma_1 \eta_1 (1 - \eta_2)} (f(x_0) - f_{\textup{low}}) \; \epsilon^{-2} \\
             & \quad \quad + \log_{\gamma_2}\left(\frac{\gamma_1 (1 - \eta_2)}{(L + 2 \mu) \;\Delta_0}\right), \nonumber
\end{align}
where we set \(\kdc = \tfrac{1}{2}\) to recover the classical trust-region Cauchy decrease condition~\citep{cartis-gould-toint-2022}.
This bound~\eqref{eq:complexity-4:alpha=1:beta=0:p=0} is of the same order as the one derived by \citet{curtis-lubberts-robinson-2018} in a similar setting, i.e., \(\alpha = 1\) and \(\beta = 0\).

For the classical trust-region method \citep{cartis-gould-toint-2022,conn-gould-toint-2000}, i.e., \(\alpha = \beta = 0\), our complexity bound reduces to
\begin{align*}
  k_\epsilon & \leq \left(|\log_{\gamma_2} (\gamma_4)| + 1 \right) \frac{2(L + 2 \mu)}{\gamma_1 \eta_1 (1-\eta_2)} (f(x_0) - f_{\textup{low}}) \; \epsilon^{-2} \\
             & \quad \quad + \log_{\gamma_2}(\epsilon) + \log_{\gamma_2} \left(\frac{\gamma_1 (1-\eta_2)}{(L + 2 \mu) \; \Delta_0}  \right),
\end{align*}
where we set \(\kdc = \tfrac{1}{2}\) and choose \(\Delta_0\) such that \(\Delta_0 \ge \|\nabla f(x_0)\|\). The obtained complexity bound is thus similar to the one derived in \citep[Theorem~\(2.3.7\)]{cartis-gould-toint-2022}.

For a better analysis on the impact of  selecting \(\alpha \) and \(\beta \), a numerical comparison of different choices of such parameters will be conducted in \Cref{sec:num}.

\subsubsection{Sharpness of the complexity bound}%
\label{sec:sharpness}

Let \(c > 0\) and \(\epsilon \in (0, \, 1)\).
Our goal is to construct smooth \(f: \R \to \R\) that satisfies \Cref{asm:problem-obj,asm:model-error,asm:unb} and for which \Cref{alg:TR-GEN} requires exactly
\begin{equation}%
  \label{eq:def-keps}
  k_{\epsilon} \defeq
  \begin{cases}
    \lfloor \epsilon^{-2/(1-p)} \rfloor              & \text{if } 0 \le p < 1, \\
    \left\lfloor \exp(c \epsilon^{-2}) \right\rfloor & \text{if } p = 1,
  \end{cases}
\end{equation}
function and gradient evaluations to produce \(x_{k_\epsilon}\) with  \(|f'(x_{k_\epsilon})| \leq \epsilon\).
The construction follows the guidelines of \citet{cartis-gould-toint-2022} and proceeds as \citet{leconte-orban-2023}.
We begin by recalling a key result on Hermite interpolation.

\begin{proposition}[{\protect \citealp[Theorem~\(A.9.2\)]{cartis-gould-toint-2022}; \citealp[Proposition~\(6\)]{leconte-orban-2023}}]%
  \label{prop:hermite-interpolation}
  Let \(k_\epsilon \in \N\).
  Consider real sequences \(\{f_k\}\), \(\{g_k\}\), and \(\{x_k\}\) for \(k \in \{0, \ldots, k_\epsilon\}\).
  For \(k = 0, \ldots, k_\epsilon\), let \(s_k \defeq x_{k+1} - x_k\), and assume that
  \begin{align*}
    |f_{k+1} - ( f_k + g_k s_k )| & \leq \kappa_f s_k^2, \\
    |g_{k+1} - g_k|               & \leq \kappa_f |s_k|,
  \end{align*}
  for \(k = 0, \ldots, k_\epsilon - 1\), where \(\kappa_f \geq 0\).
  Then, there exists continuously differentiable \(f: \R \to \R\) such that
  \[
    f(x_k) = f_k \quad \text{and} \quad f'(x_k) = g_k,
  \]
  for \(k = 0, \ldots, k_\epsilon\).
  Furthermore, if
  \[
    |f_k| \leq \kappa_f, \quad |g_k| \leq \kappa_f, \quad \text{and} \quad |s_k| \leq \kappa_f,
  \]
  for \(k = 0, \ldots, k_\epsilon\), then both \(|f|\) and \(|f'|\) are bounded by a constant depending only on \(\kappa_f\).
\end{proposition}

We proceed as \citet{leconte-orban-2023}, but consider a scenario where \Cref{asm:unb} holds for any \(p \in [0, 1]\).
For \(k \in \{0, \ldots, k_\epsilon\}\), we set
\begin{equation*}
  \omega_k \defeq \frac{k_\epsilon - k}{k_\epsilon} \quad \text{and} \quad
  g_k \defeq -\epsilon (1 + \omega_k).
\end{equation*}
By definition, \(|g_k| > \epsilon\) for all \(k \in \{0, \ldots, k_\epsilon - 1\}\) and \(|g_{k_\epsilon}| = \epsilon\).

Let $p \in [0,1]$.
Define
\begin{equation*}%
  B_0 \defeq 1 \quad \text{and} \quad B_k \defeq k^p \text{ for all } k = 1, \ldots, k_\epsilon,
\end{equation*}
and
\begin{equation*}%
  x_0 \defeq 0 \quad \text{and} \quad x_{k+1} \defeq x_k + s_k \text{ for all } k = 0, \ldots, k_\epsilon - 1,
\end{equation*}
where
\begin{equation}%
  \label{eq:s_k}
  s_k \defeq -B_k^{-1} g_k > 0 \text{ for all } k = 0, \ldots, k_\epsilon - 1.
\end{equation}
Finally, set \(f_0 \defeq
\begin{cases}
  8 \epsilon^2 + 4 c           & \text{if } p = 1       \\
  8 \epsilon^2 + \frac{4}{1-p} & \text{if } 0 \le p < 1
\end{cases}
\) and
\begin{equation}%
  \label{eq:f_k}
  f_{k+1} \defeq f_k + g_k s_k \text{ for } 0 \le k \le k_\epsilon - 1.
\end{equation}

The next lemma establishes properties of \(\{f_k\}\).
The proof is similar to that of \citep[Lemma~\(7\)]{leconte-orban-2023} and is omitted for brevity.
\begin{lemma}%
  \label{lem:fk}
  The sequence \(\{f_k\}\) defined in~\eqref{eq:f_k} is decreasing and
  \begin{equation*}%
    f_k \in [0, \, f_0] \quad \text{for all } k = 0, \ldots, k_\epsilon,
  \end{equation*}
  where \(k_\epsilon\) is defined in~\eqref{eq:def-keps}.
\end{lemma}

The next theorem establishes slow convergence of \Cref{alg:TR-GEN}.
The proof is omitted, as it is similar to that of \citep[Theorem~\(6\)]{leconte-orban-2023}.
\begin{theorem}%
  \label{thm:sharpness}
  Let \(c > 0\) and \(0 < \epsilon < 1\).
  \Cref{alg:TR-GEN} applied to minimize \(f: \R^n \to \R\) satisfying \Cref{asm:problem-obj,asm:model-error,asm:unb} with \(0\le p \le 1\) may require as many as
  \[
    k_{\epsilon} \defeq
    \begin{cases}
      \lfloor \epsilon^{-2/(1-p)} \rfloor  & \text{if } 0 \le p < 1, \\
      \lfloor \exp(c\epsilon^{-2}) \rfloor & \text{if } p = 1,
    \end{cases}
  \]
  iterations to produce \(x_{k_\epsilon}\) such that  \(\|\nabla f(x_{k_\epsilon})\| \leq \epsilon\).
\end{theorem}

\begin{figure}[ht]%
  \centering
  \begin{subfigure}{.32\textwidth}%
    \includegraphics[width=1\linewidth]{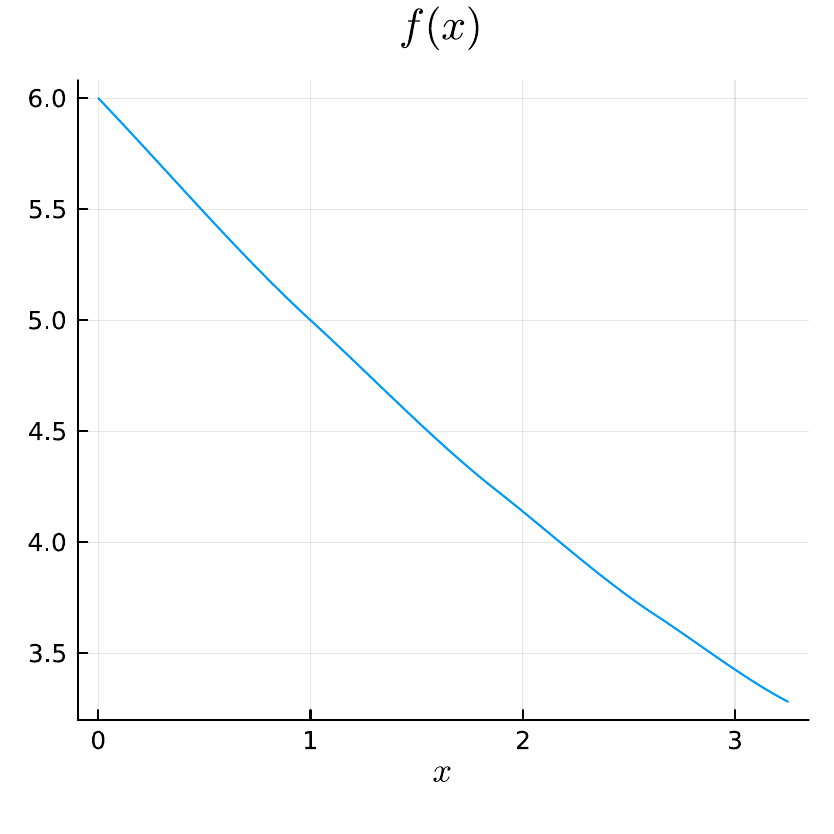}
    \\
    \includegraphics[width=1\linewidth]{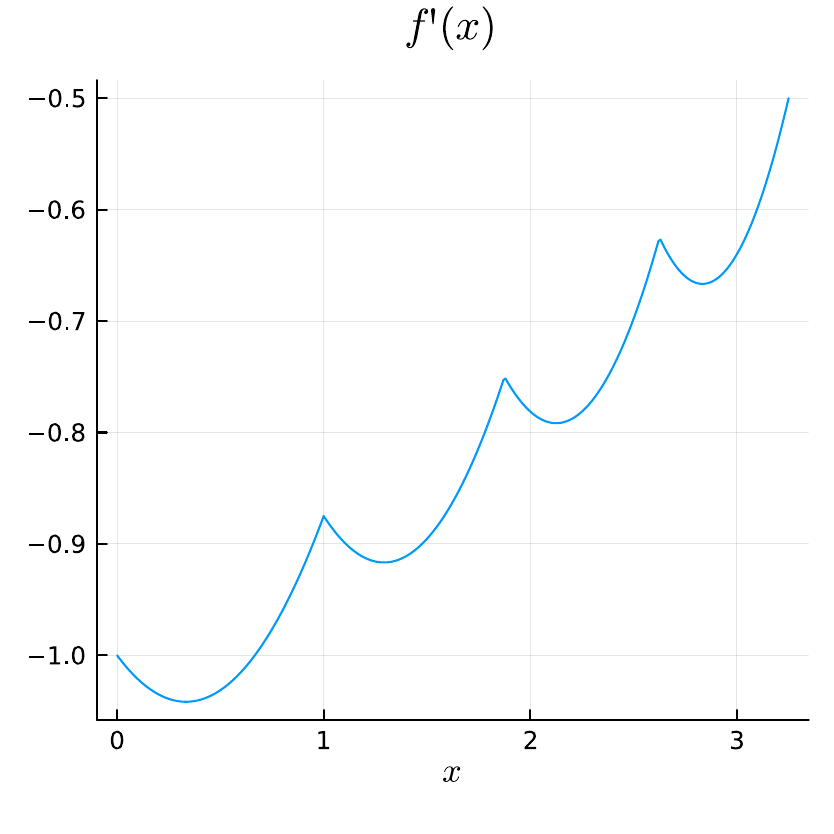}
    \caption{%
      \(p=0\).
    }
  \end{subfigure}
  \begin{subfigure}{.32\textwidth}%
    \includegraphics[width=1\linewidth]{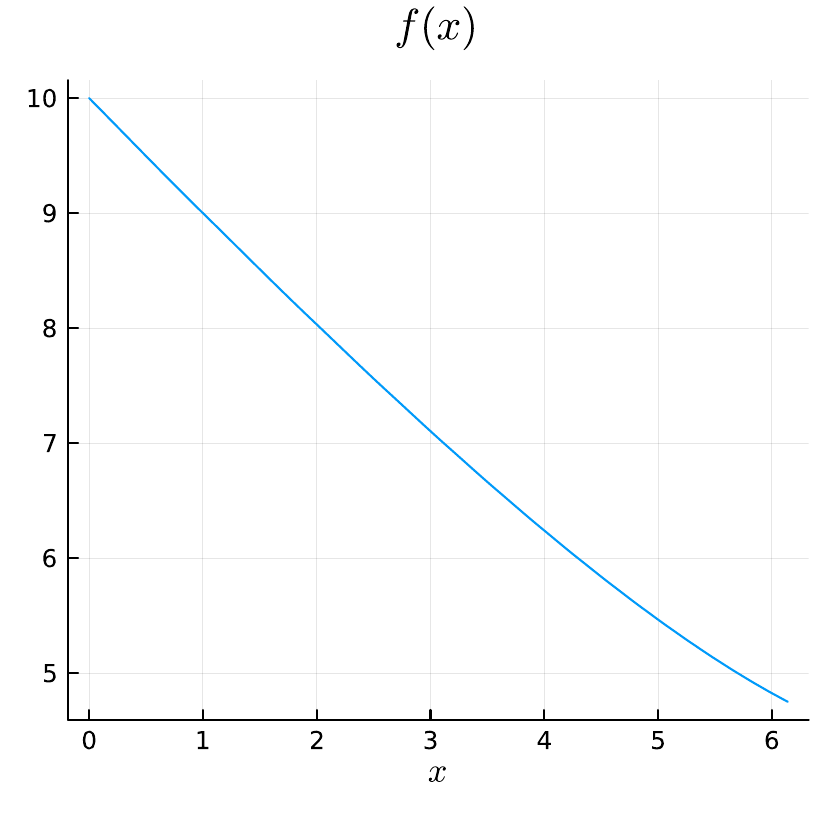}
    \\
    \includegraphics[width=1\linewidth]{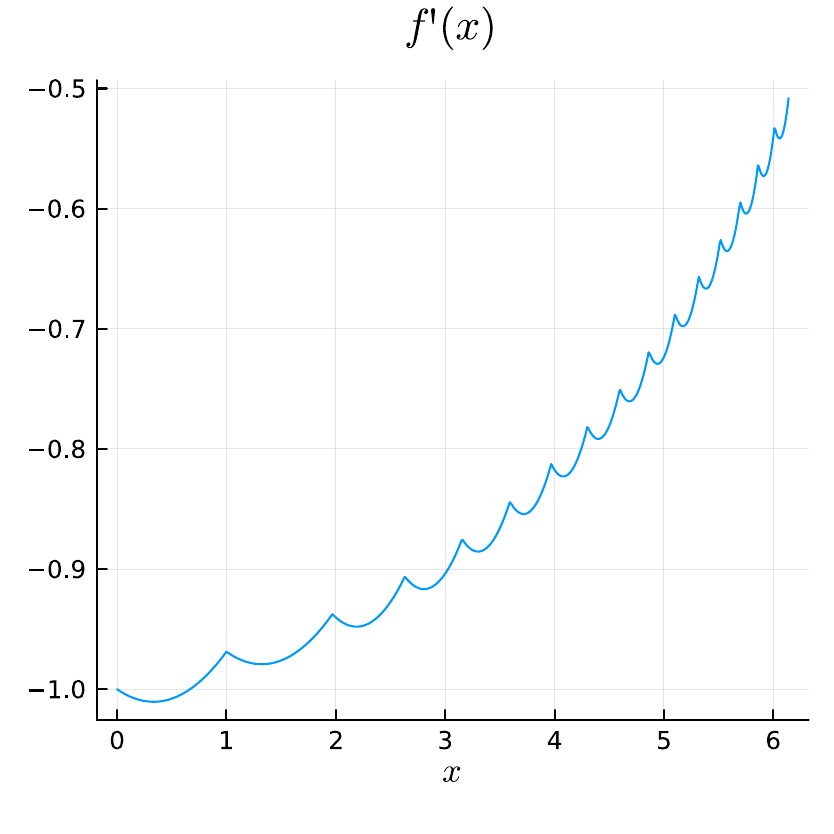}
    \caption{%
      \(p = 0.5\).
    }
  \end{subfigure}
  \begin{subfigure}{.32\textwidth}%
    \centering
    \includegraphics[width=1\linewidth]{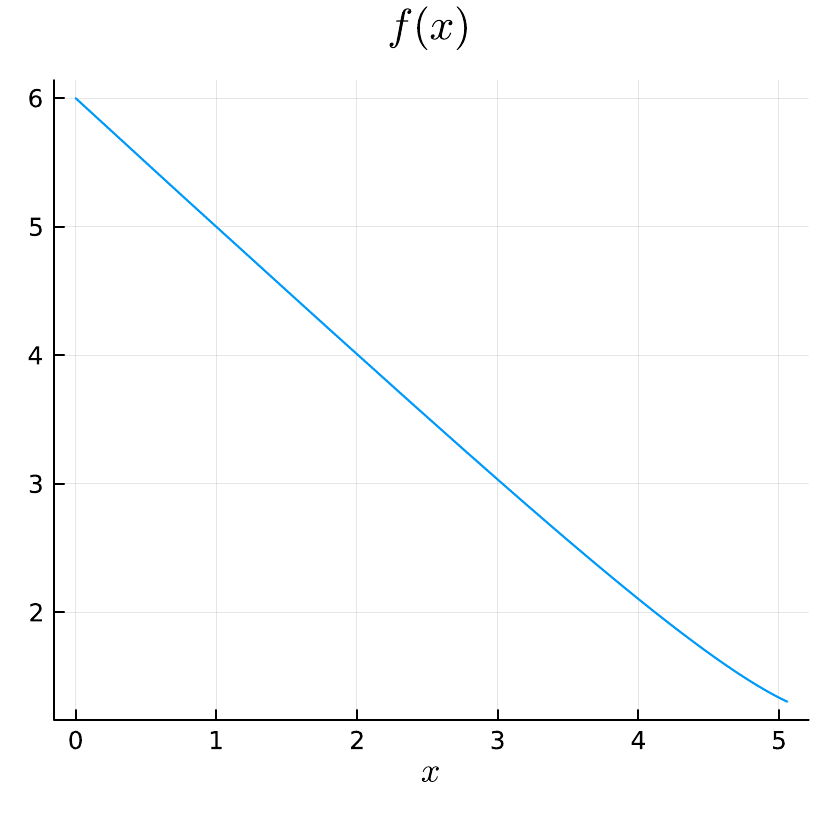}
    \centering
    \includegraphics[width=1\linewidth]{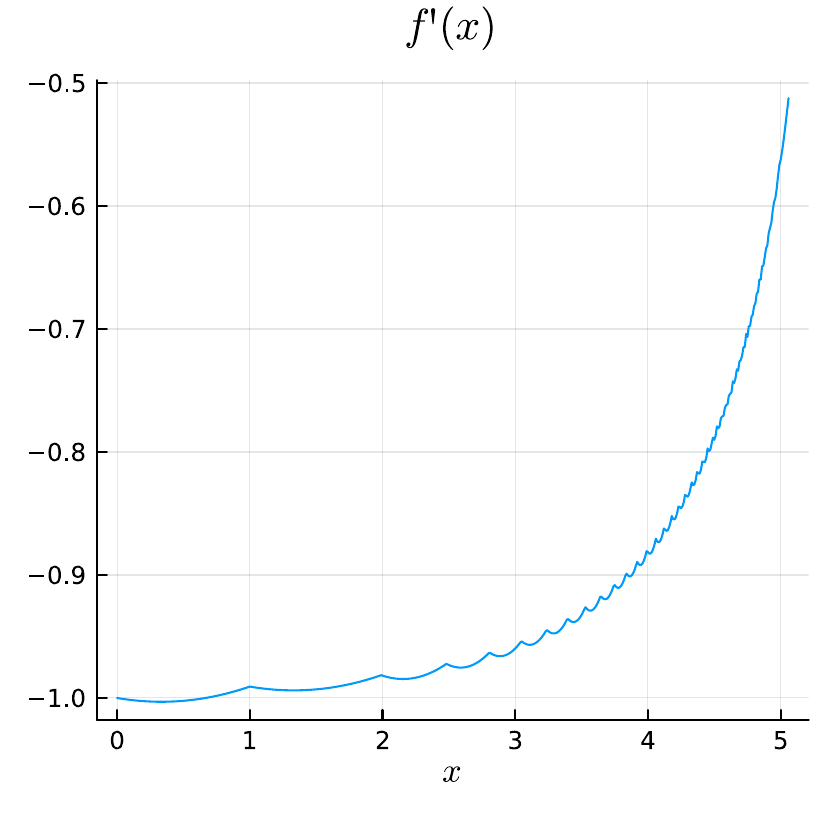}
    \caption{%
      \(p = 1\) and \(c=1\).
    }
  \end{subfigure}
  \caption{%
    \label{fig:slow-convergence}
    Form of \(f\) (top) and \(f'\) (bottom) over \([x_0, x_{k_{\epsilon}}]\) using \(\epsilon = 0.5\) and \(c = 1\) for three values of \(p\) when \(\|B_k\| = O(|\mathcal{S}_{k-1}|^p)\).
  }
\end{figure}

\Cref{fig:slow-convergence} shows plots of the counterexample in \Cref{thm:sharpness}, depicting both the objective function \(f\) and its gradient \(f'\), for \(p \in \{0, 0.5, 1\} \), \(c = 1\) and \(\epsilon = 10^{-2}\) over the interval \([x_0, x_{k_{\epsilon}}]\).

\subsection{Complexity on convex objectives}%
\label{sec:convex-S}

In the following, we derive complexity results of \Cref{alg:TR-GEN} for both convex and strongly convex functions in the presence of unbounded Hessian model.
Note that, as explained to \Cref{sec:S}, the BFGS update satisfies the bound in \Cref{asm:unb} in this case.

As expected, thanks to convexity, improvements in the obtained complexity bounds can be achieved.
In the particular case where the model Hessian is bounded and the objective function is convex or strongly convex, we show that our family yields complexity bounds similar to those established in the state of the art; see \citep{cartis-gould-toint-2012,nesterov-2013,grapiglia-yuan-2016}.

For the convexe case, we will consider the following stationary measure:

\begin{equation}%
  \label{eq:Dk}
  \delta_k \defeq f(x_k) - {f}_{\textup{low}},
\end{equation}
where \({f}_{\textup{low}} \defeq \min_{x \in \R^n} f(x)\), is the global minimum of the convex function \(f\).

Let \(\epsilon > 0\) and \(\widehat{k}_{\epsilon}\) be the first iteration of \Cref{alg:TR-GEN} such that  \(\delta_{\widehat{k}_{\epsilon}} > \epsilon \) and \(\delta_{\widehat{k}_{\epsilon} + 1} \le \epsilon \).
Define
\begin{align*}
  \widehat{\mathcal{S}}(\epsilon) & \defeq \mathcal{S}_{\widehat{k}_{\epsilon}-1} = \{ k \in \mathcal{S} \mid k < \widehat{k}_{\epsilon} \},
  \\
  \widehat{\mathcal{U}}(\epsilon) & \defeq \mathcal{U}_{\widehat{k}_{\epsilon}-1} = \{k \in \N \mid k \not \in \mathcal{S} \text{ and } k < \widehat{k}_{\epsilon} \}.
\end{align*}

For this section, we assume the following

\begin{problemassumption}%
  \label{asm:convex}
  The objective \(f\) is convex, and there exists \(R \ge 1\) such that,
  \[
    \|x - x^*\| \le R \quad \text{for all } x \text{ such that } f(x) \le f(x_0),
  \]
  where \(x^*\) is a minimizer of \(f\).
\end{problemassumption}

The second part of \Cref{asm:convex}, which addresses the boundedness of the level sets of \(f\), while restrictive, is the one used in \citep{cartis-gould-toint-2012}.
The following property enables us to relate the distance to the global minimum to the gradient norm.

\begin{lemma}[{\citealp[Lemma~\(2.4\)]{cartis-gould-toint-2012}}]%
  \label{lem:convex}
  Assume that \Cref{asm:convex} holds.
  Then, for all \(k \in \N \),
  \begin{equation*}%
    \delta_k \le R \| \nabla f(x_k) \|.
  \end{equation*}
\end{lemma}

Now we state the main complexity theorem, in the convex where the Hessian of the model satisfies \Cref{asm:unb}.
An upper bound of \(O(\epsilon^{-1})\) is also derived in the bounded case.

\begin{theorem}%
  \label{thm:complexity-convex:S}
  Let \Cref{asm:problem-obj,asm:convex,asm:model-error,asm:unb} hold.
  Assume that \Cref{alg:TR-GEN} generates infinitely many successful iterations.
  Let \(\kappa_0\) and \(\kappa_2\) be defined as in \Cref{thm:complexity:S}, \(\mu \) and \(p\) be as in \Cref{asm:unb}, \(R\) be as in \Cref{asm:convex}, and \(a_{\min}\) be as in \Cref{lem:gamma-min}.
  Then, if \(p = 0\),
  \begin{equation}%
    \label{eq:S-eps0-convex}
    |\widehat{\mathcal{S}} (\epsilon)| \leq R^2 \kappa_0 \; \epsilon^{-1} = O(\epsilon^{-1}).
  \end{equation}
  If \(0 < p < 1\),
  \begin{align}%
    \label{eq:S-eps2-convex}
    |\widehat{\mathcal{S}} (\epsilon)| & \leq \left[ (1-p)\;  (R^2 \kappa_0 \; \epsilon^{-1}-\kappa_2) + 1\right]^{1/ (1-p)} \\
                                       & = O\left([(1-p) \epsilon^{-1}]^{1/ (1-p)}\right). \nonumber
  \end{align}
  Otherwise, if \(p = 1\),
  \begin{equation}%
    \label{eq:S-exp-convex}
    |\widehat{\mathcal{S}} (\epsilon)| \leq \exp (R^2 \kappa_0 \; \epsilon^{-1}-\kappa_2).
  \end{equation}
\end{theorem}

\begin{proof}
  If \(\widehat{\mathcal{S}}(\epsilon) = \emptyset\), then the theorem holds trivially as \(|\widehat{\mathcal{S}}(\epsilon)| = 0\).

  Otherwise, consider some \(\ell \in \widehat{\mathcal{S}}(\epsilon)\).
  By applying \Cref{lem:sufficient-decrease}, in conjunction with \Cref{asm:unb} and the monotonicity of \(\fpsib\), we obtain
  \begin{align*}
    f(x_{\ell}) - f(x_{\ell} + s_{\ell}) & \geq \eta_1 \left(m_{\ell}(0) - m_{\ell}(s_{\ell})\right)                                                                                                                 \\
                                         & \geq \eta_1 \kdc a_{\min} \frac{ \|\nabla f(x_\ell)\|^{1 + \alpha} \min_{j = 0, \ldots, \ell} \|\nabla f(x_j)\|^{1 - \alpha}}{\psib{\max_{j = 0, \ldots,\ell} \|B_j\|}}   \\
                                         & \geq \eta_1 \kdc a_{\min} \frac{\|\nabla f(x_\ell)\|^{1 + \alpha} \min_{j = 0, \ldots, \ell} \|\nabla f(x_j)\|^{1 - \alpha}}{\psib{\mu(1 + |\mathcal{S}_{\ell - 1}|^p)}}.
  \end{align*}
  Furthermore, due to the decreasing property of \(\delta_k\) and the application of \Cref{lem:convex}, we conclude that for all \(j \in \{0, \ldots, \ell\}\),
  \begin{align*}
    \delta_j \le R \|\nabla f(x_j)\| & \implies \min_{j = 0, \ldots,\ell} \delta_j \le R \min_{j = 0, \ldots,\ell} \|\nabla f(x_j)\| \\
                                     & \implies \delta_\ell \le R \min_{j = 0, \ldots,\ell} \|\nabla f(x_j)\|.
  \end{align*}
  Thus, recalling~\eqref{eq:Dk} and the fact that \(\alpha \in [-1,1]\), we obtain
  \begin{eqnarray*}
    \delta_\ell - \delta_{\ell+1} \geq \frac{\eta_1 \kdc a_{\min}}{R^2} \frac{\delta_\ell^2}{\psib{\mu(1 + |\mathcal{S}_{\ell - 1}|^p)}},
  \end{eqnarray*}
  which, by dividing by \(\delta_{\ell+1}\), implies
  \begin{align*}
    \frac{1}{\delta_{\ell+1}} - \frac{1}{\delta_\ell} & \geq \frac{\eta_1 \kdc a_{\min}}{R^2} \frac{1}{\psib{\mu(1 + |\mathcal{S}_{\ell - 1}|^p)}} \frac{\delta_\ell}{\delta_{\ell+1}} \\
                                                      & \geq \frac{\eta_1 \kdc a_{\min}}{R^2} \frac{1}{\psib{\mu(1 + |\mathcal{S}_{\ell - 1}|^p)}}.
  \end{align*}
  By employing the same telescoping argument as in the proof of \Cref{thm:complexity:S}, and utilizing the facts that \(\delta_0 > 0\) and \(\delta_{\widehat{k}_{\epsilon}} > \epsilon\), we obtain
  \begin{equation}
    \label{eq:sum-series-convex}
    \epsilon^{-1} \geq \frac{1}{\delta_{\widehat{k}_{\epsilon}}} \geq \frac{1}{\delta_{\widehat{k}_{\epsilon}}} - \frac{1}{\delta_0} \geq  \frac{\eta_1 \kdc a_{\min}}{R^2} \sum_{k = 0}^{ |\widehat{\mathcal{S}}(\epsilon)|-1} \frac{1}{\psib{\mu(1 + k^p)}}
  \end{equation}
  Hence, if \(p = 0\), we have
  \begin{equation*}
    \frac{R^2 \epsilon^{-1}}{\eta_1 \kdc a_{\min}} \geq \frac{1}{\psib{2\mu}} |\widehat{\mathcal{S}}(\epsilon)|,
  \end{equation*}
  which gives~\eqref{eq:S-eps0-convex}.

  Otherwise, if \(0 < p \le 1\), \Cref{lem:sum-integral} gives
  \begin{equation}%
    \label{eq:sum-integral-convex}
    \sum_{k = 0}^{ |\widehat{\mathcal{S}}(\epsilon)|-1} \frac{1}{\psib{\mu(1 + k^p)}} \ge \frac{1}{\psib{\mu}} + \frac{1}{\psib{2\mu}}\int_{1}^{|\widehat{\mathcal{S}}(\epsilon)|} \frac{1}{t^p} \mathrm{d} t.
  \end{equation}
  We combine~\eqref{eq:sum-series-convex} with~\eqref{eq:sum-integral-convex} and obtain
  \begin{equation*}
    \frac{R^2 \epsilon^{-1}}{\eta_1 \kdc a_{\min}} \geq \frac{1}{\psib{\mu}} + \frac{1}{\psib{2\mu}}\int_{1}^{|\widehat{\mathcal{S}}(\epsilon)|} \frac{1}{t^p} \mathrm{d} t.
  \end{equation*}
  Two cases follow, the first one, \(0 < p < 1\), gives
  \begin{equation*}
    \frac{R^2 \epsilon^{-1}}{\eta_1 \kdc a_{\min}} \geq   \frac{1}{\psib{\mu}} + \frac{1}{\psib{2\mu}} \left(\frac{{|\widehat{\mathcal{S}}(\epsilon)|}^{1-p} - 1}{1-p} \right),
  \end{equation*}
  which provides~\eqref{eq:S-eps2-convex}.
  The second case, \(p = 1\), gives
  \begin{equation*}
    \frac{R^2 \epsilon^{-1}}{\eta_1 \kdc a_{\min}} \geq   \frac{1}{\psib{\mu}} + \frac{1}{\psib{2\mu}} \log\left(|\widehat{\mathcal{S}}(\epsilon)|  \right) ,
  \end{equation*}
  which establishes~\eqref{eq:S-exp-convex}.
  \qed
\end{proof}
Note that, using \Cref{thm:complexity:U}, the bound on the cardinality of \(\widehat{\mathcal{S}}(\epsilon)\) suffices to derive a similar bound on \(\widehat{k}_{\epsilon}\).
When the objective is strongly convex, under \Cref{asm:unb}, an improved iteration complexity of at most \(O([(1-p)\log(\epsilon^{-1})]^{1/(1-p)})\) for \(0 \leq p < 1\) and \(O(\epsilon^{-\kappa})\) for \(p = 1\), where \(\kappa > 0 \), can be established for \Cref{alg:TR-GEN}, as we show next.

Let us assume that the following holds, i.e., that the objective is strongly convex.

\begin{problemassumption}%
  \label{asm:str-convex}
  There exists \(\mu_c > 0\) such that,
  \[
    f(y) \ge f(x) + \nabla f(x)^\top (y - x) + \frac{\mu_c}{2} \|y - x\|^2 \quad \text{for all } x, y \in \R^n.
  \]
\end{problemassumption}

In this setting, we get a tighter relation between the distance to the global minimum and the gradient norm, that will enable us to get a better complexity bound.

\begin{lemma}[{\citep[Lemma~\(2.6\)]{cartis-gould-toint-2012}}]%
  Assume that \Cref{asm:str-convex} holds.
  Then, for all \(k \in \N\),
  \begin{equation*}%
    \delta_k \le \frac{1}{2\mu_c} \| \nabla f(x_k) \|^2.
  \end{equation*}
\end{lemma}

Now we state the main complexity theorem, in the rather strongly convex case where the Hessian of the model satisfies \Cref{asm:unb}.
An upper bound of \(O(\log(\epsilon^{-1}))\) is also recovered in the bounded case.

\begin{theorem}%
  \label{thm:complexity-convex-str:S}
  Let \Cref{asm:problem-obj,asm:model-error,asm:unb,asm:str-convex} hold.
  Assume that \Cref{alg:TR-GEN} generates infinitely many successful iterations.
  Let \(\kappa_0\) and \(\kappa_2\) be defined as in \Cref{thm:complexity:S}, \(\mu \) and \(p\) be as in \Cref{asm:unb}, \(\mu_c\) be as in \Cref{asm:str-convex}, and \(a_{\min}\) be as in \Cref{lem:gamma-min}.
  Then, if \(p = 0\)
  \begin{equation}
    \label{eq:S-eps0-str-convex}
    |\widehat{\mathcal{S}} (\epsilon)| \leq \frac{\kappa_0}{2 \mu_c} \; \log\left(\frac{\delta_0}{\epsilon}\right) = O\left(\log\left(\epsilon^{-1}\right)\right).
  \end{equation}
  If \(0 < p < 1\),
  \begin{align}
    \label{eq:S-eps2-str-convex}
    |\widehat{\mathcal{S}} (\epsilon)| & \leq \left[ (1-p)\;  \left(\frac{\kappa_0}{2 \mu_c} \; \log\left(\frac{\delta_0}{\epsilon}\right)-\kappa_2\right) + 1\right]^{1/ (1-p)} \\
                                       & = O\left(\left[(1-p)\log (\epsilon^{-1})\right]^{1/ (1-p)}\right). \nonumber
  \end{align}
  Otherwise, if \(p = 1\),
  \begin{equation}
    \label{eq:S-exp-str-convex}
    |\widehat{\mathcal{S}} (\epsilon)| \leq \exp \left(\frac{\kappa_0}{2 \mu_c}\; \log\left(\frac{\delta_0}{\epsilon}\right)-\kappa_2\right) = O\left(\epsilon^{-\frac{\kappa_0}{2 \mu_c}}\right).
  \end{equation}
\end{theorem}

\begin{proof}
  If \(\widehat{\mathcal{S}}(\epsilon) = \emptyset\), then the theorem holds trivially as \(|\widehat{\mathcal{S}}(\epsilon)| = 0\).

  Otherwise, let \(\ell \in \widehat{\mathcal{S}}(\epsilon)\).
  By invoking \Cref{lem:sufficient-decrease}, alongside \Cref{asm:unb} and the monotonic growth of \(\fpsib \), we obtain
  \begin{align*}
    f(x_{\ell}) - f(x_{\ell} + s_{\ell}) & \geq \eta_1 \left(m_{\ell}(0) - m_{\ell}(s_{\ell})\right)                                                                                                                 \\
                                         & \geq \eta_1 \kdc a_{\min} \frac{ \|\nabla f(x_\ell)\|^{1 + \alpha} \min_{j = 0, \ldots, \ell} \|\nabla f(x_j)\|^{1 - \alpha}}{\psib{\max_{j = 0, \ldots,\ell} \|B_j\|}}   \\
                                         & \geq \eta_1 \kdc a_{\min} \frac{\|\nabla f(x_\ell)\|^{1 + \alpha} \min_{j = 0, \ldots, \ell} \|\nabla f(x_j)\|^{1 - \alpha}}{\psib{\mu(1 + |\mathcal{S}_{\ell - 1}|^p)}}.
  \end{align*}
  Furthermore, due to the non-increasing nature of \(\delta_k\), we leverage \Cref{lem:convex} to establish that for all \(j \in \{0, \ldots, \ell\}\),
  \begin{align*}
    \delta_j \le \frac{1}{2\mu_c} \|\nabla f(x_j)\|^2 & \implies \min_{j = 0, \ldots,\ell} \delta_j \le \frac{1}{2\mu_c} \min_{j = 0, \ldots,\ell} \|\nabla f(x_j)\|^2 \\
                                                      & \implies \delta_\ell \le \frac{1}{2\mu_c} \min_{j = 0, \ldots,\ell} \|\nabla f(x_j)\|^2.
  \end{align*}
  Consequently, recalling~\eqref{eq:Dk} and considering that \(\alpha \in [-1,1]\), we deduce
  \begin{eqnarray*}
    \delta_\ell - \delta_{\ell+1} \geq 2\mu_c \eta_1 \kdc a_{\min} \frac{\delta_\ell}{\psib{\mu(1 + |\mathcal{S}_{\ell - 1}|^p)}},
  \end{eqnarray*}
  which directly leads to the bound,
  \begin{equation*}
    \left(1 - \frac{2\mu_c \eta_1 \kdc a_{\min}}{\psib{\mu(1 + |\mathcal{S}_{\ell - 1}|^p)}}\right) \delta_\ell \geq \delta_{\ell+1}.
  \end{equation*}
  Given that for any \(\ell \in \widehat{\mathcal{S}}(\epsilon)\), we have \(\delta_{\ell} \ge \delta_{\ell+1} > \epsilon > 0\), we can safely apply the logarithm function.
  Using the fact that \( \log(1-x) \le -x \) for \(x \in [0,1)\),
  \begin{align*}
    -\frac{2\mu_c \eta_1 \kdc a_{\min}}{\psib{\mu(1 + |\mathcal{S}_{\ell - 1}|^p)}} & \geq \log\left(1 - \frac{2\mu_c \eta_1 \kdc a_{\min}}{\psib{\mu(1 + |\mathcal{S}_{\ell - 1}|^p)}}\right) \\
                                                                                    & \geq \log\left(\delta_{\ell+1}\right) - \log\left(\delta_{\ell}\right).
  \end{align*}
  Proceeding with the same telescoping argument as in the proof of \Cref{thm:complexity:S}, and using the bounds \(\delta_0 > 0\) and \(\delta_{\widehat{k}_{\epsilon}} > \epsilon \), we conclude
  \begin{align}
    \label{eq:sum-series-str-convex}
    \log(\delta_0) - \log(\epsilon) & \geq \log(\delta_0) - \log(\delta_{\widehat{k}_{\epsilon}}) \nonumber                                                \\
                                    & \geq  2\mu_c \eta_1 \kdc a_{\min} \sum_{k = 0}^{ |\widehat{\mathcal{S}}(\epsilon)|-1} \frac{1}{\psib{\mu(1 + k^p)}}.
  \end{align}
  The case where \(p = 0\) leads to
  \begin{equation*}
    \frac{1}{2 \mu_c \eta_1 \kdc a_{\min}} \log\left(\frac{\delta_0}{\epsilon}\right)\geq \frac{1}{\psib{2\mu}} |\widehat{\mathcal{S}}(\epsilon)|,
  \end{equation*}
  which gives~\eqref{eq:S-eps0-str-convex}.

  For \(0 < p \le 1\), applying \Cref{lem:sum-integral} gives
  \begin{equation}%
    \label{eq:sum-integral-str-convex}
    \sum_{k = 0}^{ |\widehat{\mathcal{S}}(\epsilon)|-1} \frac{1}{\psib{\mu(1 + k^p)}} \ge \frac{1}{\psib{\mu}} + \frac{1}{\psib{2\mu}}\int_{1}^{|\widehat{\mathcal{S}}(\epsilon)|} \frac{1}{t^p} \mathrm{d} t.
  \end{equation}
  Combining~\eqref{eq:sum-series-str-convex} with~\eqref{eq:sum-integral-str-convex} leads to
  \begin{equation*}
    \frac{1}{2 \mu_c \eta_1 \kdc a_{\min}} \log\left(\frac{\delta_0}{\epsilon}\right) \geq \frac{1}{\psib{\mu}} + \frac{1}{\psib{2\mu}}\int_{1}^{|\widehat{\mathcal{S}}(\epsilon)|} \frac{1}{t^p} \mathrm{d} t.
  \end{equation*}
  Two cases emerge: for \(0 < p < 1\), we get
  \begin{equation*}
    \frac{1}{2 \mu_c \eta_1 \kdc a_{\min}} \log\left(\frac{\delta_0}{\epsilon}\right) \geq   \frac{1}{\psib{\mu}} + \frac{1}{\psib{2\mu}} \left[\frac{{|\widehat{\mathcal{S}}(\epsilon)|}^{1-p} - 1}{1-p} \right],
  \end{equation*}
  yielding~\eqref{eq:S-eps2-str-convex}.
  When \(p = 1\), we obtain
  \begin{equation*}
    \frac{1}{2\mu_c \eta_1 \kdc a_{\min}} \log\left(\frac{\delta_0}{\epsilon}\right) \geq  \frac{1}{\psib{\mu}} + \frac{1}{\psib{2\mu}} \log\left(|\widehat{\mathcal{S}}(\epsilon)|  \right) ,
  \end{equation*}
  which confirms~\eqref{eq:S-exp-str-convex}.
  \qed
\end{proof}

Several remarks on the complexity results \Cref{thm:complexity-convex:S,thm:complexity-convex-str:S} for (strongly) convex objectives are in order.
This is the first generalization of the complexity result in the unbounded Hessian case that specifically takes advantage of the structure of the function, such as convexity or strong convexity.

Notably, despite \Cref{asm:convex,asm:str-convex}, no convexity assumption is required for \(m_k\), reaffirming that the basic \Cref{alg:TR-GEN} framework is a first-order method.
The only model assumptions needed are \Cref{asm:model-error,asm:unb}.
Nevertheless, assuming boundedness of the level sets of \(f\) in \Cref{asm:convex} may be restrictive in certain cases.
For example, consider \(f : \mathbb{R} \to \mathbb{R}\) with \(f(x) = x^2\) for \(x \ge 0\) and \(f(x) = 0\) for \(x < 0\).
While \(f\) satisfies \Cref{asm:problem-obj} with a Lipschitz continuous gradient and is convex, the set of minimizers is \([0, +\infty)\), which is unbounded and does not satisfy \Cref{asm:convex}.
In future work, we will investigate the possibility of relaxing this assumption.

When \(p = 0\) in \Cref{asm:unb}, our bounds are similar to those established for the classical trust-region framework applied to (strongly) convex functions \citep{grapiglia-yuan-2016}.

Moreover, when \(0 < p \le 1\), \Cref{thm:complexity-convex:S,thm:complexity-convex-str:S} provide a better complexity bound than the one derived by \citep{leconte-orban-2023}, as expected.
Specifically, when \(p = 1\), \Cref{thm:complexity-convex-str:S} shows that in the strongly convex setting, the analysis yields a polynomial complexity bound, which is significantly tighter than the exponential bound established in \Cref{thm:complexity:S} for the general case.
As for the sharpness of the bound in the convex case, we will investigate it in follow-up work, as it might require different interpolation techniques.
Furthermore, it remains unclear whether we achieve the same improved bound when using the gradient metric as the stopping criterion, as developed in \citep{garmanjani-2022} for bounded model Hessians.

\section{Complexity when the bound on \texorpdfstring{\boldmath{\(\|B_k\|\)}}{∥Bₖ∥} depends only on the current iteration}%
\label{sec:k}

We now assume that \(\{B_k\}\) is bounded by a multiple of \(k\).
Among other scenarios, this suggests that in \Cref{alg:TR-GEN}, one is allowed to update the model Hessian at both successful and unsuccessful iterations.
Note the bounded case \(p = 0\) in \Cref{asm:unb} has been studied in \Cref{sec:S}, and we now consider \(p > 0\).

We replace \Cref{asm:unb} with the following assumption.
\begin{modelassumption}%
  \label{asm:unb-k}
  There exists \(\mu > 0\) and \(0 < p \leq 1\) such that
  \begin{equation*}
    \max_{j = 0, \ldots,k} \|B_j\| \le \mu (1 + k^p) \qquad (k \in \N).
  \end{equation*}
\end{modelassumption}

As before, we begin with the case where only a finite number of successful iterations is generated.

\begin{theorem}%
  Let \Cref{asm:problem-obj,asm:model-error,asm:unb-k} be satisfied.
  If \Cref{alg:TR-GEN} only generates finitely many successful iterations, then \(x_k = x^*\) for all sufficiently large \(k\) where \(\|\nabla f(x^*)\| = 0\).
\end{theorem}

\begin{proof}
  Assume by contradiction that there exists \(\nu > 0\) such that \(\|\nabla f(x_k)\| \geq \nu > 0\) for all \(k \in \N\), and let \(k_f\) be the last successful iteration.
  Necessarily, \(x_k = x_{k_f}\) for all \(k \geq k_f\), and hence, \(\{x_k\} \to x^* \defeq x_{k_f}\).
  By \Cref{lem:gamma-min,lem:technical}, for all \(k \in \N\),
  \begin{align*}
    \Delta_k \geq \frac{(\min_{j = 0, \ldots, k} \|\nabla f(x_j)\|)^{1-\alpha} \, a_{\min}}{\max_{j = 0, \ldots,k} [(L + \|B_j\|) (1 + \|B_j\|)^{- \beta}] } & \geq \frac{\nu^{1-\alpha}\, a_{\min}}{ L(1 + \max_{j = 0, \ldots,k}\|B_j\|)^{1- \beta}} \\
                                                                                                                                                             & \geq \frac{\nu^{1-\alpha} a_{\min}}{L(1 + \mu(1 + k^p))^{1 - \beta}}.
  \end{align*}
  On unsuccessful iterations, \Cref{alg:TR-GEN} reduces \(\Delta_k\) by a factor at least \(\gamma_2\).
  Hence, for all \(k \geq k_f\),
  \[
    \Delta_k \leq \gamma_2^{k-k_f} \Delta_{k_f}.
  \]
  We combine the above inequalities, and obtain
  \[
    \frac{\nu^{1-\alpha} a_{\min}}{L(1 + \mu (1 + k^p))^{1 - \beta}} \leq \gamma_2^{k-k_f} \Delta_{k_f},
  \]
  which may be rewritten
  \[
    0 < \frac{\nu^{1-\alpha} a_{\min} \, \gamma_2^{k_f}}{L \Delta_{k_f}} \leq \gamma_2^k (1 + \mu (1 + k^p))^{1 - \beta}.
  \]
  However, the above is a contradiction as the right-hand side goes to zero as \(k \to \infty\).
  Thus, by contradiction, \(\liminf_{k \to \infty} \|\nabla f(x_k)\| =  \|\nabla f(x_{k_f})\| = \|\nabla f(x^*)\| = 0\).
  \qed
\end{proof}

Let \(\tau \in \N_0\) and \(\lambda \in \N \).
Define
\begin{subequations}%
  \label{eq:def-T-W}
  \begin{align}%
    \mathcal{T}^{\tau,\lambda}_k & \defeq \left\{ j = \lambda, \ldots, k \mid j < \tau |\mathcal{S}_j| + \lambda \right\},
    \\
    \label{Wk}
    \mathcal{W}^{\tau,\lambda}_k & \defeq \left\{ j = \lambda, \ldots, k \mid j \geq \tau |\mathcal{S}_j| + \lambda \right\}.
  \end{align}
\end{subequations}

Note that the set \(\mathcal{T}_k^{\tau,\lambda}\) might be empty.
The next technical result will allow us to relate the number of successful iterations to the number of iterations.

\begin{lemma}%
  \label{lem:comparison-series}
  Let \({\{r_j\}}_{j \in \N}\) be a non-decreasing positive real sequence.
  For \(k \ge \lambda\),
  \begin{equation*}
    \tau \sum_{j \in \mathcal{S}_k} \frac{1}{r_j} \ge \sum_{j \in \mathcal{T}^{\tau,\lambda}_k} \frac{1}{r_j} = \sum_{j = \lambda }^k \frac{1}{r_j} - \sum_{j \in\mathcal{W}^{\tau,\lambda}_k } \frac{1}{r_j},
  \end{equation*}
  where \(\mathcal{T}^{\tau,\lambda}_k \) and \(\mathcal{W}^{\tau,\lambda}_k \) are defined in~\eqref{eq:def-T-W}.
\end{lemma}

\begin{proof}
  If \( \mathcal{T}^{\tau, \lambda}_k\) is empty, the result is trivial, using the convention that the sum over an empty set is zero.

  Otherwise, let \(k \ge \lambda\) and \(j \in \mathcal{T}^{\tau,\lambda}_k\), i.e., $\lambda\leq j \leq k$ and \(j < \tau |\mathcal{S}_j| + \lambda\).
  Define \(\mathcal{S}^{\tau}_j\) as the list of elements of \(\mathcal{S}_j\), in ascending order, where each element is repeated \(\tau \) times:
  \begin{equation}%
    \label{eq:def-Skt}
    \mathcal{S}^{\tau}_j \defeq \Big\{ \underbrace{i_1, \ldots, i_1}_{\tau \text{ times}}, \ldots, \underbrace{i_{|\mathcal{S}_j|}, \ldots, i_{|\mathcal{S}_j|}}_{\tau \text{ times}} \Bigm| i_l \in \mathcal{S}_j \text{ for } l \in \{1,\ldots,|\mathcal{S}_j|\} \Big\}.
  \end{equation}
  By construction, \(|\mathcal{S}^{\tau}_j| = \tau |\mathcal{S}_j| \ge j - \lambda + 1\), hence \(\mathcal{S}^{\tau}_j\) must contain at least \(j-\lambda+1\) elements.
  Let \(i^{\tau}_{j-\lambda}\) be the \(j-\lambda+1\)-th element of \(\mathcal{S}^{\tau}_j\).
  In particular, \(i_{j-\lambda}^\tau \in \mathcal{S}_j\), and because each element of \(\mathcal{S}_j\) is less than \(j\), \(i^{\tau}_{j-\lambda} \le j\).
  Because \(\mathcal{S}^{\tau}_j \subseteq \mathcal{S}^{\tau}_{k}\) also holds, \(i^{\tau}_{j-\lambda} \in \mathcal{S}^{\tau}_k\).
  We have just showed that
  \begin{equation}%
    \label{eq:ijt}
    \{ i^\tau_{j-\lambda} \mid j \in \mathcal{T}^{\tau,\lambda}_k \} \subseteq \mathcal{S}^{\tau}_k.
  \end{equation}
  As \(\{r_j\}\) is non-decreasing and \(i_{j-\lambda}^\tau \leq j\), \(r_{i^{\tau}_{j-\lambda}} \le r_j\).
  We sum over \(j \in \mathcal{T}^{\tau,\lambda}_k\):
  \begin{equation}%
    \label{eq:sum-ratio}
    \sum_{j \in \mathcal{T}^{\tau,\lambda}_k} \frac{1}{r_{i^\tau_{j-\lambda}}} \ge \sum_{j \in \mathcal{T}^{\tau,\lambda}_k } \frac{1}{r_j}.
  \end{equation}
  Positivity of \(\{r_j\}_{j \in \N}\) and~\eqref{eq:ijt} then yield
  \begin{equation}%
    \label{eq:sum-ratio-2}
    \sum_{j \in \mathcal{S}^{\tau}_k} \frac{1}{r_j} \ge \sum_{j \in \mathcal{T}^{\tau,\lambda}_k} \frac{1}{r_{i^\tau_{j-\lambda}}}.
  \end{equation}
  Thus, we combine~\eqref{eq:def-T-W},~\eqref{eq:def-Skt},~\eqref{eq:sum-ratio} and~\eqref{eq:sum-ratio-2} to conclude that
  \begin{align*}
    \tau \sum_{j \in \mathcal{S}_k} \frac{1}{r_j} = \sum_{j \in \mathcal{S}^{\tau}_k} \frac{1}{r_j} \ge \sum_{j \in \mathcal{T}^{\tau,\lambda}_k} \frac{1}{r_{i^\tau_{j-\lambda}}} \ge \sum_{j \in \mathcal{T}^{\tau,\lambda}_k} \frac{1}{r_j} = \sum_{j = \lambda}^k \frac{1}{r_j} - \sum_{j \in\mathcal{W}^{\tau,\lambda}_k } \frac{1}{r_j}.
    \tag*{\qed}
  \end{align*}
\end{proof}

The following lemma also plays a key role in deriving our worst-case complexity bound.

\begin{lemma}%
  \label{lem:xi}
  Let \Cref{asm:problem-obj,asm:model-error,asm:unb-k} be satisfied.
  Assume that \(\tau \in \N_0\) is chosen so that \(\gamma_4 \gamma_2^{\tau-1} < 1\) and \(\lambda \in \N \) such that \(\gamma_2^\lambda \le \frac{a_{\min} \epsilon^{1 - \alpha}}{\Delta_0}\).
  Let \(\epsilon > 0\) and \(k_{\epsilon}\) be the first iteration such that \(\|\nabla f(x_{k_{\epsilon}})\| \le \epsilon\).
  Then,
  \begin{equation*}%
    \sum_{k \in \mathcal{W}^{\tau}_{k_\epsilon-1}} \frac{1}{\psib{\mu(1 + k^p)}} \le \xi \,
    \qquad
    \xi \defeq \sum_{k \in \N} \left(\gamma_4 \gamma_2^{\tau-1}\right)^{ k / \tau} < \infty,
  \end{equation*}
  where \(a_{\min }\) is as in \Cref{lem:gamma-min} and $\mathcal{W}^{\tau}_{k_\epsilon-1}$ is defined in~\eqref{Wk}.
\end{lemma}

\begin{proof}
  Let \(k \in \mathcal{W}^{\tau}_{k_\epsilon-1}\).
  By \Cref{asm:unb-k}, the increasing nature of \(\fpsib\) and \(\beta \ge 0\),
  \begin{align*}%
    \frac{a_{\min}\epsilon^{1-\alpha}}{\psib{\mu(1 + k^p)}} & \le \frac{a_{\min}\min_{j = 0, \ldots,k} \|\nabla f(x_j)\|^{1-\alpha} }{\psib{\max_{j = 0, \ldots,k} \|B_j\|}}                                                  \\
                                                            & \le \frac{a_{\min}\min_{j = 0, \ldots,k} \|\nabla f(x_j)\|^{1-\alpha} }{\psib{\max_{j = 0, \ldots,k} \|B_j\|}(1 + \max_{j = 0, \ldots,k} \|B_j\|)^{^{-\beta}}}.
  \end{align*}
  Thus, \Cref{lem:Delta-min} together with~\eqref{eq:L-beta-psi} yield
  \begin{align}%
    \label{eq:Delta-k-min}
    \frac{a_{\min}\epsilon^{1-\alpha}}{\psib{\mu(1 + k^p)}} & \le \frac{a_{\min}\min_{j = 0, \ldots,k} \|\nabla f(x_j)\|^{1-\alpha} }{\psib{\max_{j = 0, \ldots,k} \|B_j\|}(1 + \max_{j = 0, \ldots,k} \|B_j\|)^{^{-\beta}}} \nonumber \\
                                                            & \leq \Delta_k.
  \end{align}
  On the other hand, the update mechanism of \(\Delta_k\) in \Cref{alg:TR-GEN} implies
  \begin{equation*}%
    \Delta_k \leq \gamma_4^{|S_{k-1}|} \gamma_2^{|U_{k-1}|} \Delta_0 \leq \gamma_4^{|S_{k}|} \gamma_2^{k - |S_{k}|} \Delta_0.
  \end{equation*}
  Because \(k \in \mathcal{W}^{\tau}_{k_\epsilon-1}\), \(k \ge \tau |\mathcal{S}_k| + \lambda \), which, together with the fact that \(\gamma_4 \ge 1\), \(0 < \gamma_2 < 1\) leads to
  \begin{align}%
    \label{eq:gamma-ineq}
    \Delta_k \leq \gamma_4^{|S_k|} \gamma_2^{k-|S_k|} \Delta_0 & \leq \gamma_4^{(k - \lambda)/\tau} \gamma_2^{k - (k - \lambda)/\tau} \Delta_0 \nonumber           \\
                                                               & = \left(\gamma_4 \gamma_2^{\tau-1}\right)^{ (k - \lambda) / \tau} \gamma_2^{\lambda} \, \Delta_0.
  \end{align}
  Since \(\gamma_2^\lambda \le \frac{a_{\min} \epsilon^{1 - \alpha}}{\Delta_0}\),~\eqref{eq:gamma-ineq} implies
  \begin{equation}%
    \label{eq:gamma-ineq-lambda}
    \Delta_k \leq \left(\gamma_4 \gamma_2^{\tau-1}\right)^{ (k - \lambda) / \tau} \gamma_2^{\lambda} \, \Delta_0 \leq \left(\gamma_4 \gamma_2^{\tau-1}\right)^{ (k - \lambda) / \tau} a_{\min}  \epsilon^{1 - \alpha}.
  \end{equation}
  By combining~\eqref{eq:Delta-k-min} and~\eqref{eq:gamma-ineq-lambda} and summing over \(k \in \mathcal{W}^{\tau}_{k_\epsilon-1}\),
  \begin{equation*}%
    \sum_{k \in \mathcal{W}^{\tau}_{k_\epsilon-1}} \frac{1}{\psib{\mu(1 + k^p)}} \le \sum_{k \ge \lambda} \left(\gamma_4 \gamma_2^{\tau-1}\right)^{ (k - \lambda) / \tau} \le \sum_{k \in \N} \left(\gamma_4 \gamma_2^{\tau-1}\right)^{ k / \tau}.
    \tag*{\qed}
  \end{equation*}
\end{proof}

Note that it is possible to choose \(\tau \) and \(\lambda \) as required by \Cref{lem:xi}.
In fact, it suffices to pick
\begin{equation}
  \label{eq:leps}
  \tau > \log_{\gamma_2}(\gamma_4^{-1}) + 1 \ge 1
  \quad \textup{and} \quad
  1 \le \lambda \leq 1 + \log_{\gamma_2}\left(\frac{a_{\min} \epsilon^{1 - \alpha}}{\Delta_0}\right).
\end{equation}
The condition on \(\lambda \) is always satisfied since \(\frac{a_{\min} \epsilon^{1 - \alpha}}{\Delta_0} < 1\) by definition of \( a_{\min} \) in \Cref{lem:gamma-min} and \(0 < \gamma_2 < 1\).

In the following, we establish the evaluation complexity of \Cref{alg:TR-GEN} under \Cref{asm:unb-k}.

\subsection{Complexity on nonconvex objectives}%

Given the tools developed in \Cref{sec:k}, we are ready to state our main result on the evaluation complexity of \Cref{alg:TR-GEN} under \Cref{asm:unb-k}, that will enable to directly bound the total number of iterations required to reach a solution with a desired accuracy.

\begin{theorem}%
  \label{thm:complexity-k}
  Let \Cref{asm:problem-obj,asm:model-error,asm:unb-k} be satisfied.
  Assume that \Cref{alg:TR-GEN} generates infinitely many successful iterations.
  Let \(\kappa_1\) be defined as in \Cref{thm:complexity:S}, \(\mu\) and \(p\) be as in \Cref{asm:unb-k}, and \(\tau \) and \(\xi \) be as in \Cref{lem:xi}.
  Let \(\epsilon >0\), and \(k_{\epsilon}\) be the first iteration such that \(\|\nabla f(x_{k_{\epsilon}})\| \le \epsilon\).
  Define
  \begin{equation*}%
    \kappa_{3} \defeq \psib{2\mu} \, \xi > 0
    \quad \mbox{and} \quad
    \kappa_{4} \defeq \log_{\gamma_2}\left(\frac{a_{\min}}{\Delta_0}\right) + 1 > 0.
  \end{equation*}
  If \(0 < p < 1\),
  \begin{align}
    k_\epsilon & \leq {\left[(1-p) \left(\tau \kappa_1 \; \epsilon^{-2} + \kappa_{3} \right) + \left((1-\alpha)\log_{\gamma_2}(\epsilon) + \kappa_{4}\right)^{1-p}\right]}^{1/(1-p)}
    \label{eq:ke-p<1}
    \\
               & = O\left([(1-p) \epsilon^{-2}]^{1/(1-p)} \right). \nonumber
  \end{align}
  If \(p = 1\),
  \begin{align}
    k_\epsilon & \leq \left((1-\alpha)\log_{\gamma_2}(\epsilon) + \kappa_{4}\right) \exp \left(\tau \kappa_1 \; \epsilon^{-2} + \kappa_{3} \right).
    \label{eq:ke-p=1}
  \end{align}
\end{theorem}

\begin{proof}
  Let \(\lambda \) be as in~\eqref{eq:leps}.
  If \(k_\epsilon \le \lambda \), the result holds, as \(\lambda \le \log_{\gamma_2}\left(\frac{a_{\min} \epsilon^{1 - \alpha}}{\Delta_0}\right) + 1\).

  Otherwise, let \(k \in \mathcal{S}(\epsilon)\).
  \Cref{lem:sufficient-decrease}, the increase nature of \(\fpsib\) and \Cref{asm:unb-k} imply
  \begin{align*}
    f(x_k) - f(x_k + s_k) & \geq \eta_1 (m_k(0) - m_k(s_k))                                                                                                                                 \\
                          & \geq \eta_1 \kdc a_{\min} \frac{ \|\nabla f(x_k)\|^{1 + \alpha} \min_{j = 0, \ldots, k} \|\nabla f(x_j)\|^{1 - \alpha}}{\psib{\max_{j = 0, \ldots, k} \|B_j\|}} \\
                          & \geq \eta_1 \kdc a_{\min} \epsilon^{2} \frac{1}{\psib{\mu(1 + k^p)}}.
  \end{align*}
  We sum the above inequality over all \(k \in \mathcal{S}(\epsilon)\), use a telescoping argument, and obtain
  \begin{equation*}
    f(x_0) - {f}_{\textup{low}} \geq \eta_1 \kdc a_{\min} \epsilon^2 \sum_{k \in \mathcal{S}(\epsilon)} \frac{1}{\psib{\mu(1 + k^p)}} .
  \end{equation*}
  The sequence $\{\psib{\mu(1 + k^p)}\}$ is positive and increasing, so \Cref{lem:sum-integral,lem:comparison-series,lem:xi} yield
  \begin{align*}%
    f(x_0) - {f}_{\textup{low}} & \geq \frac{\eta_1 \kdc a_{\min}\epsilon^2}{\tau} (\sum_{k = \lambda}^{k_\epsilon - 1} \frac{1}{\psib{\mu(1 + k^p)}}  - \sum_{k\in \mathcal{W}^{\tau}_{k_\epsilon -1}} \frac{1}{\psib{\mu(1 + k^p)}})
    \\
                                & \geq \frac{\eta_1 \kdc a_{\min}\epsilon^2}{\tau} (\sum_{k = \lambda}^{k_\epsilon - 1} \frac{1}{\psib{\mu(1 + k^p)}} -  \xi ).
  \end{align*}
  By applying \Cref{lem:sum-integral} with \( k_1 = \lambda \ge 1 \) and \( k_2 = k_\epsilon - 1 \), we obtain
  \begin{align*}%
    f(x_0) - {f}_{\textup{low}} & \geq \frac{\eta_1 \kdc a_{\min}\epsilon^2}{\tau} ( \frac{\lambda^p}{\psib{\mu (1 + \lambda^p)}} \int_{\lambda}^{k_\epsilon} \frac{1}{t^p} \, \mathrm{d}t  - \xi )
    \\
                                & \geq \frac{\eta_1 \kdc a_{\min}\epsilon^2}{\tau} ( \frac{1}{\psib{2 \mu}} \int_{\lambda}^{k_\epsilon} \frac{1}{t^p} \, \mathrm{d}t  - \xi ),
  \end{align*}
  where we used that \( \psib{\mu (1 + \lambda^p)} \leq \psib{2 \mu \lambda^p} \leq \lambda^p  \psib{2 \mu} \), since \( \lambda \geq 1 \) and \( \fpsib \) is increasing.

  We distinguish two cases.
  The first one, \(0 < p < 1\), gives
  \begin{equation*}
    f(x_0) - {f}_{\textup{low}} \geq \frac{\eta_1 \kdc a_{\min}\epsilon^2}{\tau} ( \frac{1}{\psib{2 \mu}} \frac{{k_\epsilon}^{1-p} - {\lambda}^{1-p}}{1-p}  - \xi ),
  \end{equation*}
  which combined with~\eqref{eq:leps} yields~\eqref{eq:ke-p<1}.
  The second case, \(p = 1\), gives
  \begin{equation*}
    f(x_0) - {f}_{\textup{low}} \geq \frac{\eta_1 \kdc a_{\min}\epsilon^2}{\tau} ( \frac{1}{\psib{2 \mu}} \log\left(\frac{k_\epsilon}{\lambda}\right)  - \xi ),
  \end{equation*}
  which combined with~\eqref{eq:leps} gives~\eqref{eq:ke-p=1}.
  \qed
\end{proof}

The complexity bounds of \Cref{thm:complexity-k} have the same nature as those of \Cref{thm:complexity:S,thm:complexity:U}; a polynomial bound for $0 < p < 1$ and an exponential bound for $p=1$.
In \Cref{thm:complexity-k}, we bound the total number of iterations to reach an $\epsilon$-first order point; we do no have estimates on the number of successful iterations as we did under \Cref{asm:unb}.

In addition, \Cref{thm:complexity-k} suggests that when we relax \Cref{asm:unb} to \Cref{asm:unb-k}, the complexity bound deteriorates slightly because of the multiplicative factor \(\tau \) of \(\epsilon^{-2}\) in~\eqref{eq:ke-p<1} and~\eqref{eq:ke-p=1}.
Moreover, when \(p = 1\), there is an additional multiplicative factor of order \((1-\alpha)\log(\epsilon^{-1}) + 1\) of the exponential term in~\eqref{eq:ke-p=1}.

In the next section, we show the tightness of the complexity bounds obtained in the case where the model Hessian \(B_k\) depends on \(k\).

\subsubsection{Sharpness of the complexity bound}%

For the case \(p \in [0, 1)\) in \Cref{thm:complexity-k}, one might notice that the leading term in the complexity bound is \(\epsilon^{-2/(1-p)}\), which is the same as the one derived in \Cref{thm:complexity:S}.
Thus, the complexity bound of \Cref{thm:complexity-k} is sharp for \(p \in [0, 1)\) according to \Cref{sec:sharpness}.

However, for the case \(p = 1\), the leading term in the complexity bound is \(\log(\epsilon^{-1}) \exp(\epsilon^{-2})\), which is worse than the one derived in \Cref{thm:complexity:S}.
In the next section, we will show that in this case, the complexity bound of \Cref{thm:complexity-k} is sharp.
Interestingly, to establish this sharpness, we must ensure the existence of an arbitrarily large number of unsuccessful iterations.
This guarantees that the model Hessian \(B_k\) is updated infinitely often during these iterations; otherwise, the analysis would still fall under the previous \Cref{sec:S}.
In other words, unsuccessful iterations must occur to distinguish between the complexity bounds of \Cref{thm:complexity:S} and \Cref{thm:complexity-k}.
Next, we introduce an example function for which the complexity bound related to the case \(p = 1\) is achievable.

Let  \(0 < \gamma_2 < 1\) and \(-1 \le \alpha \le 1\) be parameters as defined in \Cref{alg:TR-GEN}.
Set \(\Delta_0 \defeq 2^{2-\alpha}\).
For simplicity, when an unsuccessful iteration occurs, we consider \(\Delta_{k+1} = \gamma_2 \Delta_k\).
Finally, let \(\epsilon \in (0, \, 1)\) denote the desired accuracy level and \(c > 0\) be a constant.
Our goal is to construct smooth \(f: \R \to \R\) that satisfies \Cref{asm:problem-obj,asm:model-error,asm:unb-k} for \(p = 1\) and for which \Cref{alg:TR-GEN} requires exactly
\begin{equation}%
  \label{eq:def-keps-k}
  k_{\epsilon} \defeq  (\ueps + 1) \left\lfloor \exp(c \epsilon^{-2}) \right\rfloor
\end{equation}
function evaluations to produce \(x_{k_\epsilon}\) with  \(|f'(x_{k_\epsilon})| \leq \epsilon\), where \(u_{\epsilon} = \lfloor(1-\alpha)\log_{\gamma_2}(\epsilon)\rfloor\), and \(\gamma_2\) is as in \Cref{alg:TR-GEN}.
In particular, \(c\) could be chosen as \(\tau \kappa_1\) where \(\tau\) and \(\kappa_1\) are as in \Cref{thm:complexity-k}.

For \(k \in \{0, \ldots, k_\epsilon\}\), we set
\begin{align}
  \label{eq:g_k-k}
  \omega_k \defeq \frac{k_\epsilon - k}{k_\epsilon} \quad \text{and} \quad
  g_k \defeq
  \begin{cases}
    -\epsilon (1 + \omega_0) & \text{for } 0 \le k \le \ueps,        \\
    -\epsilon (1 + \omega_k) & \text{for } \ueps < k \le k_\epsilon.
  \end{cases}
\end{align}
By definition, \(|g_k| > \epsilon\) for all \(k \in \{0, \ldots, k_\epsilon - 1\}\) and \(|g_{k_\epsilon}| = \epsilon\).

Define
\begin{equation}%
  \label{eq:B_k-k}
  B_k \defeq
  \begin{cases}
    1 & \text{for } 0 \le k \le \ueps,        \\
    k & \text{for } \ueps < k \le k_\epsilon,
  \end{cases}
\end{equation}
and
\begin{equation}%
  \label{eq:x_k-k}
  x_0 \defeq 0 \quad \text{and} \quad s_k \defeq -B_k^{-1} g_k > 0 \text{ for all } k = 0, \ldots, k_\epsilon - 1.
\end{equation}
Finally, set
\begin{equation}%
  \label{eq:f_k-k}
  f_k \defeq
  \begin{cases}
    8 \epsilon^2 + 4c         & \text{for } 0 \le k \le \ueps,       \\
    f_{k-1} + g_{k-1} s_{k-1} & \text{for } \ueps< k \le k_\epsilon.
  \end{cases}
\end{equation}

The next lemma establishes properties of \(\{f_k\}\) similar to that in \Cref{lem:fk}.
\begin{lemma}%
  \label{lem:fk-k}
  The sequence \(\{f_k\}\) defined in~\eqref{eq:f_k-k} is non-increasing and
  \begin{equation*}%
    f_k \in [0, \, f_0] \quad \text{for all } k = 0, \ldots, k_\epsilon,
  \end{equation*}
  where \(k_\epsilon\) is defined in~\eqref{eq:def-keps-k}.
\end{lemma}

\begin{proof}
  For \( k \in \{ 0, \ldots, \ueps \} \), the proof is trivial because \( f_k = f_0 \).

  For \( k > \ueps \), we have \( f_{k} - f_{k-1} = g_{k-1} s_{k-1} < 0 \), which implies that \( \{ f_k \} \) is non-increasing.
  It follows that \(f_k \le f_0\) for all \(k = \ueps + 1, \ldots, k_\epsilon\).

  We now show that \(f_k \geq 0\) for all \(k = \ueps + 1, \ldots, k_\epsilon\).
  For \( k = \ueps + 1 \), we have
  \begin{equation*}
    f_0 - f_{\ueps + 1} = f_{\ueps} - f_{\ueps+1} = -g_{\ueps} s_{\ueps} = g_{\ueps} B_{\ueps}^{-1} g_{\ueps} = 4 \epsilon^2,
  \end{equation*}
  hence, \( f_{\ueps + 1} = f_0 - 4 \epsilon^2 = 4 \epsilon^2 + 4 c > 0 \).

  For \( k = \ueps + 2, \ldots, k_\epsilon \),
  \begin{align*}
    f_0 - f_k = f_{u_\epsilon} - f_k = \sum_{i=u_\epsilon}^{k-1} (f_i - f_{i+1}) & = -g_{\ueps} s_{\ueps} - \sum_{i=u_\epsilon + 1}^{k-1} g_i s_i                      \\
                                                                                 & = 4 \epsilon^2 + \sum_{i=1}^{k-1} \epsilon^2 {(1 + \omega_i)}^2 i^{-1}              \\
                                                                                 & = \epsilon^2 \left( 4 + \sum_{i=u_\epsilon}^{k-1} {(1 + \omega_i)}^2 i^{-1}\right).
  \end{align*}
  Because \(1 + \omega_i \le 2\),
  \begin{align*}
    \sum_{i=u_\epsilon + 1}^{k-1} {(1 + \omega_i)}^2 i^{-1} & \leq \sum_{i=u_\epsilon + 1}^{k-1} 4 i^{-1} = 4 \left((\ueps + 1)^{-1} + \sum_{i=\ueps+2}^{k-1} i^{-1}\right)                                            \\
                                                            & \le 4 \left(1 + \sum_{i=\ueps + 2}^{k-1} \int_{i-1}^{i} t^{-1} \, \mathrm{d} t \right) = 4 \left(1 + \int_{\ueps+1}^{k-1} t^{-1} \, \mathrm{d} t\right),
  \end{align*}
  so that
  \[
    f_{0} - f_k \leq 4 \epsilon^2 \left( 2 +  \int_{u_\epsilon + 1}^{k-1} t^{-1} \, \mathrm{d} t\right).
  \]
  Hence,
  \begin{align*}
    f_0 - f_k  \leq 4 \epsilon^2\left( 2 + \log(\frac{k-1}{u_\epsilon + 1})\right) & \leq 4 \epsilon^2\left( 2  + \log(\frac{k_\epsilon}{u_\epsilon + 1})\right) \\
                                                                                   & \leq 4 \epsilon^2\left(2 + c\epsilon^{-2}\right) = f_0,
  \end{align*}
  by definition of \(k_\epsilon\) in~\eqref{eq:def-keps-k}.
  Thus, \(f_0 - f_k \leq f_0\), which implies \(f_k \geq 0\).
  \qed
\end{proof}

The next theorem establishes slow convergence of \Cref{alg:TR-GEN} similar to the one in \Cref{thm:sharpness}.
\begin{theorem}%
  \label{thm:sharpness-k}
  Let \(0 < \epsilon < 1\).
  \Cref{alg:TR-GEN} applied to minimize \(f: \R^n \to \R\) satisfying \Cref{asm:problem-obj,asm:model-error,asm:unb-k} with \(p = 1\) may require as many as \(k_\epsilon\)
  iterations to produce \(x_{k_\epsilon}\) such that  \(\|\nabla f(x_{k_\epsilon})\| \leq \epsilon\), where \(k_\epsilon\) is defined in~\eqref{eq:def-keps-k}.
\end{theorem}

\begin{proof}
  We check that the iterates~\eqref{eq:x_k-k}, the function values~\eqref{eq:f_k-k}, the gradient values~\eqref{eq:g_k-k}, and the step sizes~\eqref{eq:s_k} satisfy the assumptions of \Cref{prop:hermite-interpolation}.
  By construction, for \(k = 0, \ldots, u_\epsilon - 1\), \(|f_{k+1} - f_k - g_k s_k| = |g_k s_k| = B_k |s_k|^2 = |s_k|^2\).
  Whereas, for \(k = u_\epsilon, \ldots, k_\epsilon - 1\), \(|f_{k+1} - f_k - g_k s_k| = 0\).

  On the other hand, for \(k = 0, \ldots, u_\epsilon-1\), \(|g_{k+1} - g_k| = 0 \leq |s_k|\) and for \(k = u_\epsilon\), \(|g_{k+1} - g_k| = \epsilon (u_\epsilon + 1)/k_\epsilon \leq 2 \epsilon \leq |s_k|\).
  For \(k = u_\epsilon + 1, \ldots, k_\epsilon - 1\),
  \begin{equation*}
    |g_{k+1} - g_k| = \epsilon |\omega_k -\omega_{k+1}| = \epsilon k_\epsilon^{-1} \leq B_k^{-1} |g_k| = |s_k|.
  \end{equation*}
  \Cref{lem:fk-k} yields \(f_k \in [0, \, f_0]\) for \(k = 0, \ldots, k_\epsilon\), and~\eqref{eq:g_k-k}--\eqref{eq:x_k-k} imply
  \[
    |g_k| \leq 2\epsilon \leq 2 \quad \text{and} \quad 0 < s_k \leq |g_k| \leq 2, \quad k = 0, \ldots, k_\epsilon.
  \]
  The bounds above pave the way to \Cref{prop:hermite-interpolation}, with \(\kappa_f \defeq \max \left(f_0, 2\right)\).

  The final step is to check consistency with \Cref{alg:TR-GEN}, and ensure that, for \(k = 0, \ldots, \ueps\), only unsuccessful iterations are performed, and, for \(k = \ueps + 1, \ldots, k_\epsilon - 1\), only successful iterations are performed.

  Each iteration \( k = 0, \ldots, \ueps \) is unsuccessful because
  \[
    \rho_k = \frac{f_k - f_{k+1}}{f_k - f_k - g_k s_k - \tfrac{1}{2} B_k s_k^2}  = 0.
  \]
  Thus, \( x_{k+1} = x_k \).
  Each iteration \( k = \ueps + 1, \ldots, k_\epsilon - 1 \) is successful because
  \[
    \rho_k = \frac{f_k - f_{k+1}}{f_k - f_k - g_k s_k - \tfrac{1}{2} B_k s_k^2}  = \frac{-g_k s_k}{\tfrac{1}{2} B_k^{-1} g_k^2}  = \frac{B_k^{-1} g_k^2}{\tfrac{1}{2} B_k^{-1} g_k^2} = 2.
  \]
  Thus, \( x_{k+1} = x_k + s_k \).

  The next step is to show that each \( s_k \) is inside the trust-region.
  Let \(k \in \{0, \ldots, k_\epsilon - 1\}\).
  By construction, \((1 + B_k)^\beta \le 1 + B_k \le 2 B_k\), since \(\beta \le 1\) and \(B_k \ge 1\).
  Thus,
  \begin{align}
    \label{eq:s_k-trust-region}
    |s_k| = |B_k^{-1} g_k| = \frac{|g_k|}{B_k} & \le \frac{2|g_k|^{\alpha} |g_k|^{1-\alpha}}{(1+B_k)^\beta} \nonumber                                                                           \\
                                               & \le \frac{2 |g_k|^\alpha (2)^{1-\alpha} \epsilon^{1-\alpha} }{(1+B_k)^\beta} \le \frac{|g_k|^\alpha}{(1+B_k)^\beta} \gamma_2^{\ueps} \Delta_0,
  \end{align}
  where we used that \(|g_k| \leq 2\epsilon\), \(\Delta_0 = 2^{2-\alpha}\) and \(\epsilon^{1-\alpha} \le \gamma_2^{\ueps}\) by definition of \(\ueps = \lfloor (1 - \alpha) \log_{\gamma_2}(\epsilon) \rfloor\) with \(\gamma_2 < 1\).

  Now, we would like to prove that \(\Delta_k \ge \Delta_0 \gamma_2^{\ueps}\) for all \(k = 0, \ldots, k_\epsilon - 1\).
  If \(k \in \{0, \ldots, \ueps\}\), the iteration is unsuccessful, so
  \[
    \Delta_k = \Delta_0 \gamma_2^{k} \ge \Delta_0 \gamma_2^{\ueps},
  \]
  because \(\gamma_2 < 1\).
  If \(k \in \{\ueps + 1, \ldots, k_\epsilon - 1\}\), the iteration is successful, hence
  \[
    \Delta_k \ge \Delta_{\ueps} = \Delta_0 \gamma_2^{\ueps}.
  \]
  In both cases, we have \(\Delta_k \ge \Delta_0 \gamma_2^{\ueps}\).

  Therefore,~\eqref{eq:s_k-trust-region} yields
  \begin{equation*}
    |s_k| \le \frac{|g_k|^\alpha}{(1+B_k)^\beta} \Delta_k,
  \end{equation*}
  which implies that each \(s_k\) lies inside the trust-region.

  Moreover, \(s_k\) is the minimizer of the quadratic model \(m_k\) at \(x_k\) defined in~\eqref{subprob:TR}.
  Consequently, the sufficient decrease condition~\eqref{eq:suff-decrease} is satisfied.

  We deduce from \Cref{prop:hermite-interpolation} that there exists a continuously differentiable \(f: \R \to \R\) that satisfies \Cref{asm:model-error} by construction and such that, for \(k = 0, \ldots, k_\epsilon \), \(f(x_0) = f_k\) and \(f'(x_k) = g_k\).
  Moreover, \(|g_k| > \epsilon\) for all \(k \in \{0, \ldots, k_\epsilon - 1\}\) and \(|g_{k_\epsilon}| = \epsilon\) by definition.
  Since the range of \(f\) lies in \([-\kappa_f, \, \kappa_f]\), where \(\kappa_f = \max \left(f_0, 2\right)\), \(f\) is bounded below.
  \qed
\end{proof}

\begin{figure}[ht]%
  \centering
  \begin{subfigure}{.32\textwidth}%
    \includegraphics[width=1\linewidth]{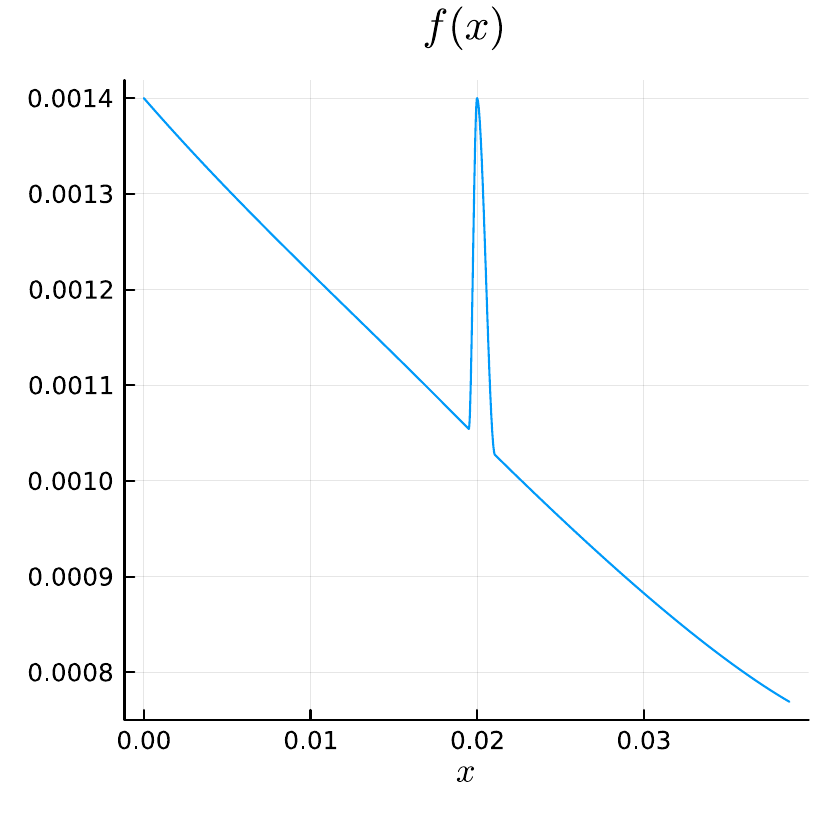}
  \end{subfigure}
  \qquad \qquad \qquad \qquad
  \begin{subfigure}{.32\textwidth}%
    \includegraphics[width=1\linewidth]{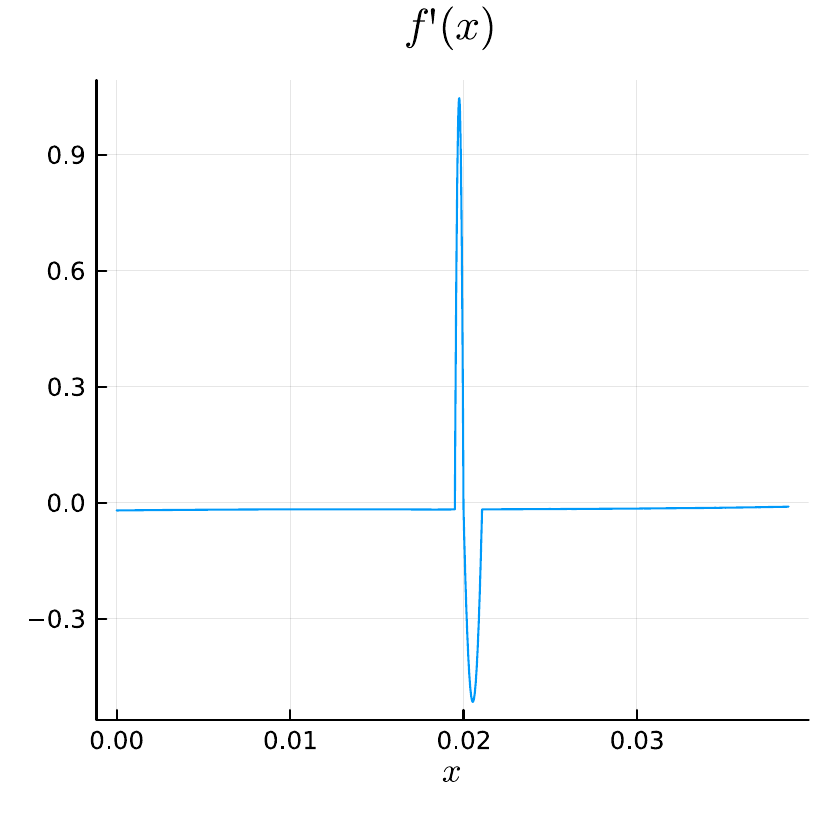}
  \end{subfigure}
  \caption{%
  \label{fig:slow-convergence-k-interpol}
  Form of \(f\) (left) and \(f'\) (right) over \([x_0, x_{k_{\epsilon}}]\) with \(\epsilon = 10^{-2}\), \(c = 1.5 \times 10^{-4}\) and \(B_k\) defined in~\eqref{eq:B_k-k}  (i.e., \(\|B_k\| = O(k)\)).
  }
\end{figure}

\begin{figure}[ht]%
  \centering
  \begin{subfigure}{.32\textwidth}%
    \includegraphics[width=1\linewidth]{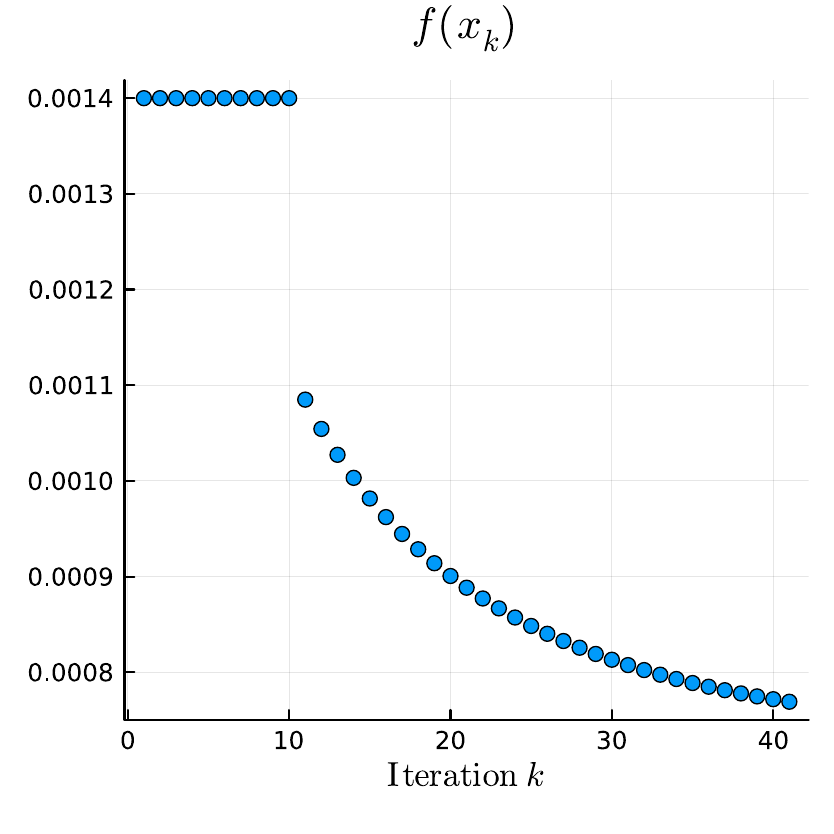}
  \end{subfigure}
  \qquad \qquad \qquad \qquad
  \begin{subfigure}{.32\textwidth}%
    \includegraphics[width=1\linewidth]{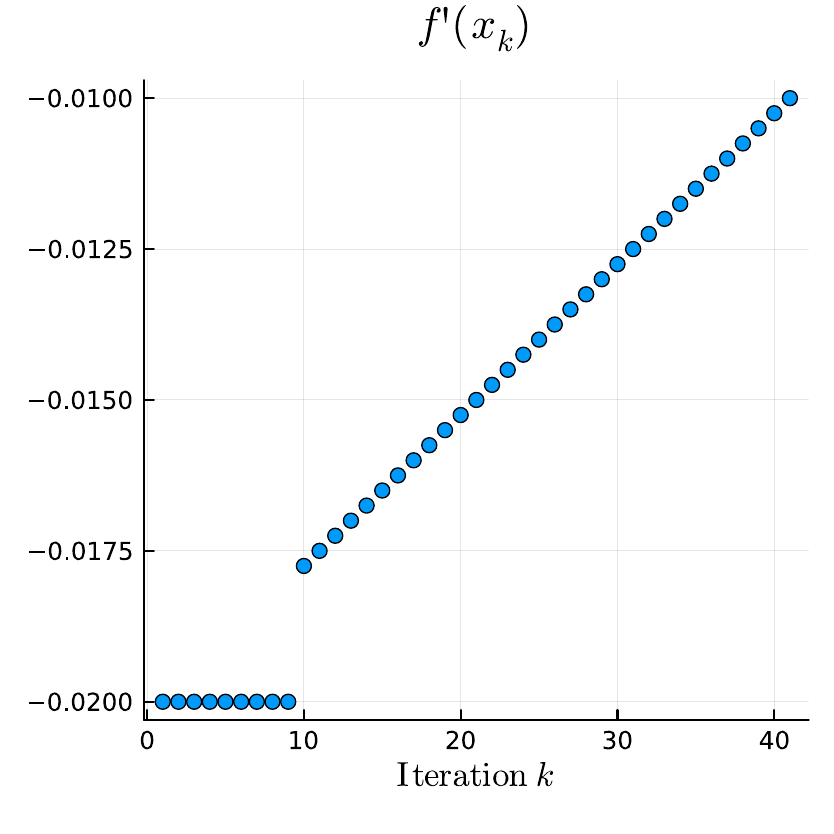}
  \end{subfigure}
  \caption{%
  \label{fig:slow-convergence-k-discret}
  Form of \(f(x_k)\) (left) and \(f'(x_k)\) (right) with \(\epsilon = 10^{-2}\), \(c = 1.5 \times 10^{-4}\) for \(x_k\) defined in~\eqref{eq:x_k-k} and \(B_k\) defined in~\eqref{eq:B_k-k} (i.e., \(\|B_k\| = O(k)\)).
  }
\end{figure}

\Cref{fig:slow-convergence-k-interpol,fig:slow-convergence-k-discret} shows plots of the counterexample in \Cref{thm:sharpness-k}, depicting both the objective function \(f\) and its gradient \(f'\) in both interpolated and discrete formats, with \(c = 1.5 \times 10^{-4}\) and \(\epsilon = 10^{-2}\) over the interval \([x_0, x_{k_{\epsilon}}]\).
We set \(\alpha = 0\) and \(\gamma_2 = e^{-1}\).
These parameters were chosen to make the plots clearer and to avoid generating too many points.
One might notice the existence of several unsuccessful iterations in \Cref{fig:slow-convergence-k-discret}, which is a consequence of the construction of the function \(f\) in \Cref{thm:sharpness-k}.

In the following as in \Cref{sec:convex-S}, we derive complexity bounds for the convex and strongly convex cases under \Cref{asm:unb-k}.

\subsection{Complexity on convex objectives}%

In the beginning, we assume that the objective function is convex and satisfies \Cref{asm:convex}.
We will show that the complexity of \Cref{alg:TR-GEN}, under \Cref{asm:unb-k}, will deteriorate slightly compared to the global complexity bound of \Cref{thm:complexity-convex:S}.
The proof of the following theorem is a direct combination of the proof of \Cref{thm:complexity-k} and the proof of \Cref{thm:complexity-convex:S}, and is thus omitted.

\begin{theorem}%
  \label{thm:complexity-convex:k}
  Let \Cref{asm:problem-obj,asm:convex,asm:model-error,asm:unb-k} be satisfied.
  Assume that \Cref{alg:TR-GEN} generates infinitely many successful iterations.
  Let \(\kappa_0\) be defined as in \Cref{thm:complexity:S}, \(\mu \) and \(p\) be as in \Cref{asm:unb-k}, \(R\) be as in \Cref{asm:convex}, and \(\tau \) be as in \Cref{lem:xi}.
  Moreover recall the constants \(\kappa_{3}\) and \(\kappa_{4}\) from \Cref{thm:complexity-k}.
  Let \(\epsilon > 0\), and \(\widehat{k}_{\epsilon}\) be the first iteration such that \(\delta_{\widehat{k}_{\epsilon}} > \epsilon \) and \(\delta_{\widehat{k}_{\epsilon} + 1} \le \epsilon \).

  Then, if \(0 < p < 1\),
  \begin{align*}%
    \widehat{k}_{\epsilon} & \leq {\left[ (1-p)\;  (\tau R^2 \kappa_0 \; \epsilon^{-1} + \kappa_{3}) + \left( (1-\alpha)\log_{\gamma_2}(\epsilon) + \kappa_{4} \right)^{1-p} \right]}^{1/(1-p)} \\
                           & = O\left([(1-p) \epsilon^{-1}]^{1/(1-p)}\right). \nonumber
  \end{align*}
  Otherwise, if \(p = 1\),
  \begin{equation*}%
    \widehat{k}_{\epsilon} \leq \left( (1-\alpha)\log_{\gamma_2}(\epsilon) + \kappa_{4} \right) \exp \left( \tau R^2 \kappa_0 \; \epsilon^{-1} + \kappa_{3}\right).
  \end{equation*}
\end{theorem}

Now, we assume that the objective function is strongly convex, i.e., satisfies \Cref{asm:str-convex}.
Same as in \Cref{thm:complexity-convex:k}, we will show that the complexity of \Cref{alg:TR-GEN}, under \Cref{asm:unb-k}, will also deteriorate slightly compared to the complexity bound of \Cref{thm:complexity-convex-str:S}.
The proof of the following theorem will be also omitted as it is a direct combination of the proofs of \Cref{thm:complexity-k,thm:complexity-convex-str:S}.

\begin{theorem}%
  \label{thm:complexity-convex-str:k}
  Let \Cref{asm:problem-obj,asm:model-error,asm:unb-k,asm:str-convex} hold.
  Assume that \Cref{alg:TR-GEN} generates infinitely many successful iterations.
  Let \(\kappa_0\) be defined as in \Cref{thm:complexity:S}, \(\mu \) and \(p\) be as in \Cref{asm:unb-k}, \(\mu_c\) be as in \Cref{asm:str-convex}, and \(\tau \) be as in \Cref{lem:xi}.
  Moreover recall the constants \(\kappa_{3}\) and \(\kappa_{4}\) from \Cref{thm:complexity-k}.
  Let \(\epsilon > 0\), and \(\widehat{k}_{\epsilon}\) be the first iteration such that \(\delta_{\widehat{k}_{\epsilon}} > \epsilon \) and \(\delta_{\widehat{k}_{\epsilon} + 1} \le \epsilon \).

  Then, if \(0 < p < 1\),
  \begin{align*}
    \widehat{k}_{\epsilon} & \leq {\left[ (1-p)\;  \left( \frac{ \tau \kappa_0}{2 \mu_c} \; \log\left(\frac{\delta_0}{\epsilon}\right) + \kappa_{3}\right) + \left( (1-\alpha)\log_{\gamma_2}(\epsilon) + \kappa_{4} \right)^{1-p} \right]}^{1/(1-p)}
    \\
                           & = O\left(\left[(1-p)\log (\epsilon^{-1})\right]^{1/ (1-p)}\right).
  \end{align*}
  Otherwise, if \(p = 1\),
  \begin{align*}
    \widehat{k}_{\epsilon} & \leq \left((1-\alpha)\log_{\gamma_2}(\epsilon) + \kappa_{4} \right) \exp \left( \frac{\tau \kappa_0}{2 \mu_c}\; \log\left(\frac{\delta_0}{\epsilon}\right) + \kappa_{3}\right)
    \\
                           & = O\left(\left((1-\alpha)\log_{\gamma_2}(\epsilon) + \kappa_{4}\right) \; \epsilon^{- \frac{\tau \kappa_0}{2 \mu_c}}\right).
  \end{align*}
\end{theorem}

\section{Numerical illustration}%
\label{sec:num}

In order to better understand the connection between worst-case complexity bounds and performance in practice, we implemented \Cref{alg:TR-GEN} in the Julia programming language using the \emph{JuliaSmoothOptimizers} infrastructure \citep{jso}.
We altered the implementation of the \emph{trunk} solver \citep{jso_solvers} to implement the trust-region radius management of \Cref{alg:TR-GEN}.
By default, \emph{trunk} is a standard trust-region method for~\eqref{eq:nlp} in which steps are computed using the truncated conjugate-gradient method \citep{steihaug-1983}.
No other algorithmic changes were performed.
We extracted all unconstrained problems from the \emph{OptimizationProblems} collection \citep{optimization_problems} with at least two variables, which resulted in \(150\) problems with dimension ranging from \(2\) to \(100\).
Larger instances are available by changing the default number of variables in many problems, but our intention here is to obtain preliminary results that may confirm or infirm causality between worst-case complexity and performance.
Typical complexity studies in nonconvex deterministic optimization are, unfortunately, almost never accompanied by numerical experiments.

\begin{figure}[ht]%
  \centering
  \includegraphics[width=.32\linewidth]{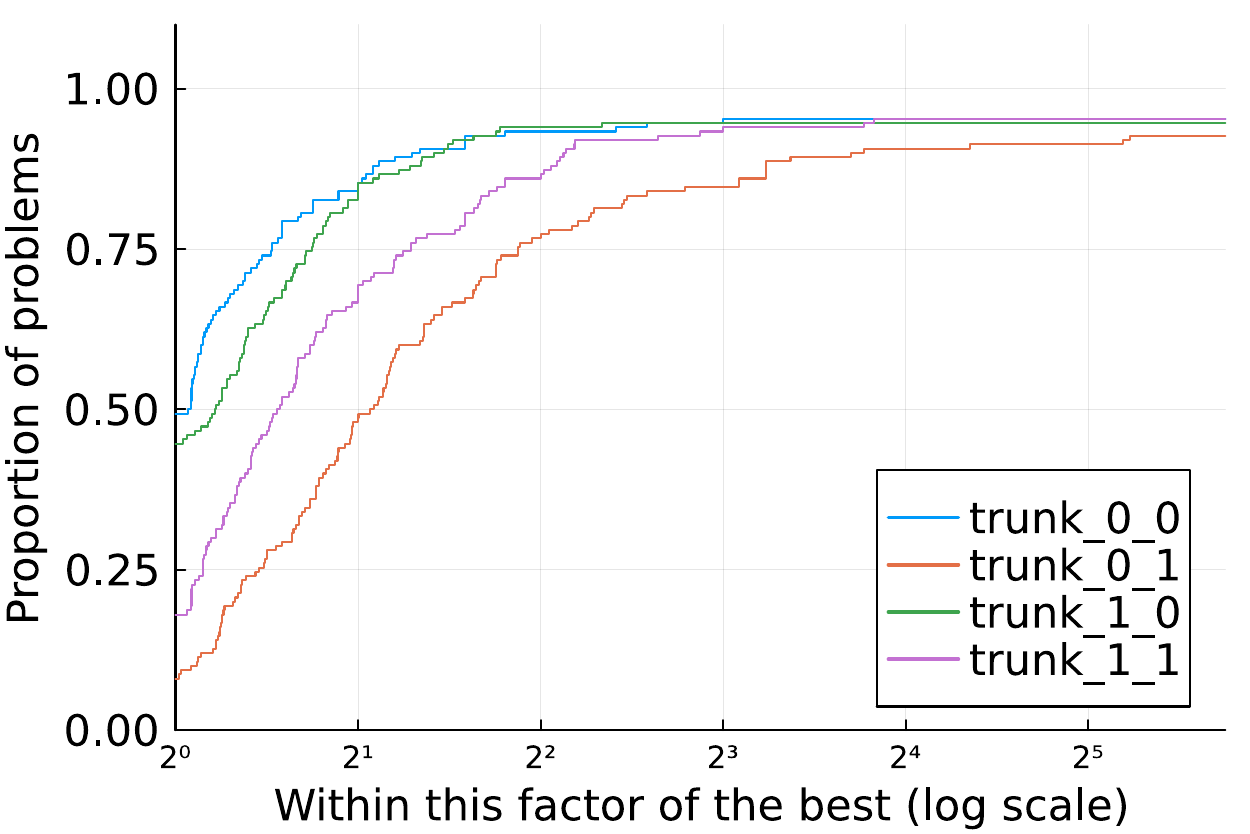}
  \hfill
  \includegraphics[width=.32\linewidth]{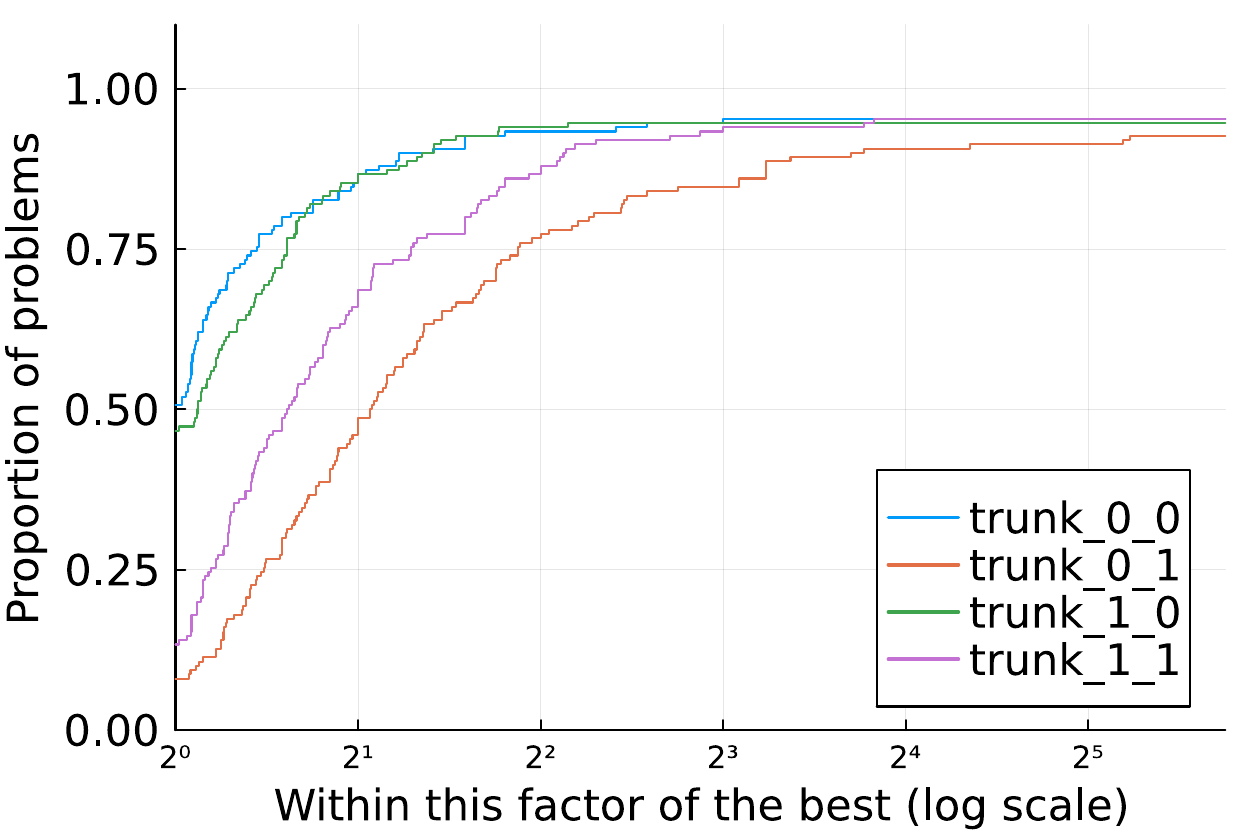}
  \hfill
  \includegraphics[width=.32\linewidth]{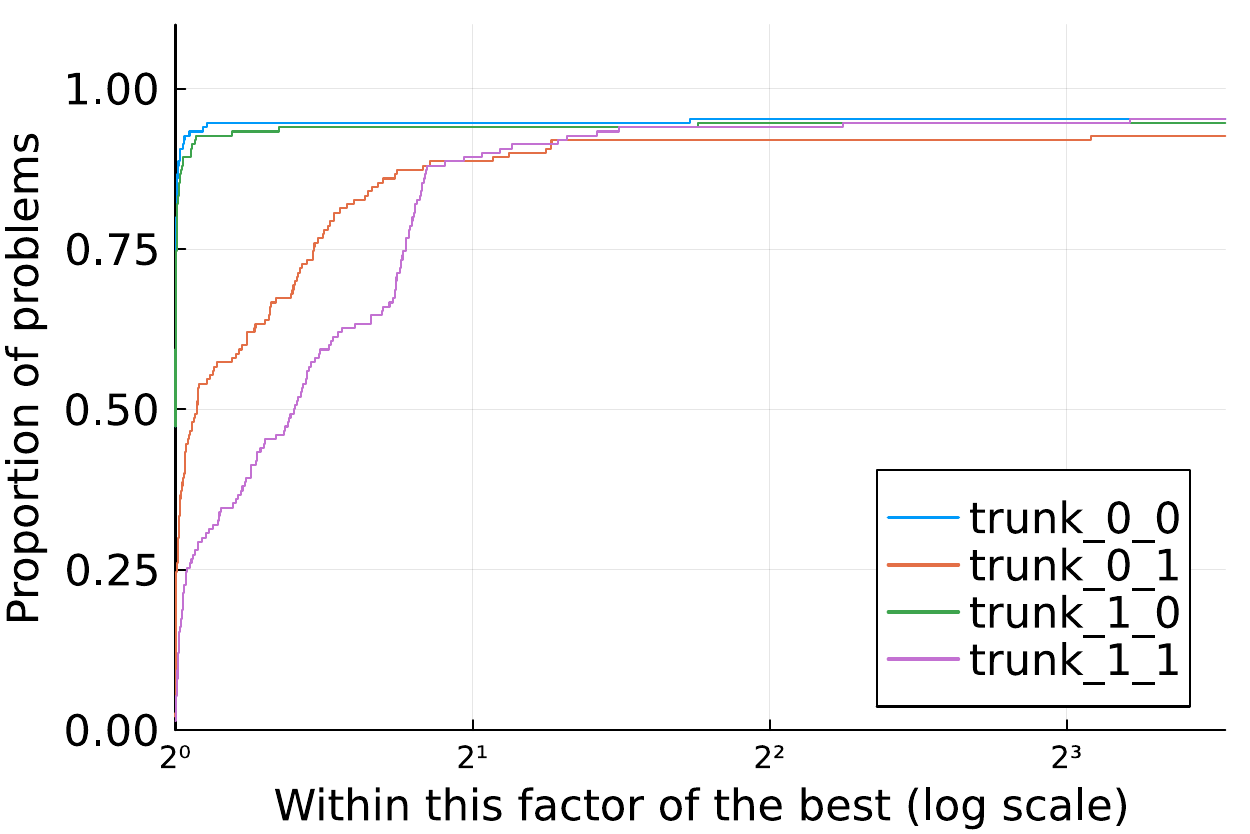}
  \caption{%
    \label{fig:profiles-trunk}
    Performance profiles of \emph{trunk} with exact Hessian in terms of number of objective evaluations (left), gradient evaluations (center), and CPU time (right).
    Method \texttt{trunk\_a\_b} refers to the variant \(\alpha = a\) and \(\beta = b\).
  }
\end{figure}

\Cref{fig:profiles-trunk} reports our results in the form of \citet{dolan-more-2002} performance profiles in terms of total number of evaluations of \(f\) (left), \(\nabla f\) (center), and CPU time (right).
At each iteration, \(B_k\) is the exact Hessian computed via automatic differentiation.
The variants tested are \(\alpha = \beta = 0\), named \texttt{trunk\_0\_0}, \(\alpha = 0\) and \(\beta = 1\), named \texttt{trunk\_0\_1}, \(\alpha = 1\) and \(\beta = 0\), named \texttt{trunk\_1\_0}, and \(\alpha = \beta = 1\), named \texttt{trunk\_1\_1}.

Not all variants are equally efficient or robust.
The standard implementation \(\alpha = \beta = 0\) is both the fastest and most robust with failures on \(7\) problems.
The method using the radius proposed by \citet{fan-yuan-2001}, \(\alpha = 1\) and \(\beta = 0\), is a close second in terms of efficiency, but fails on \(8\) problems.
The two other variants are visibly less efficient.
The variant \texttt{trunk\_1\_1} ranks third in terms of number of evaluations, but last in terms of CPU time (though the gap with the next method is small), and fails on \(7\) problems.
Finally, \texttt{trunk\_0\_1} ranks last in terms of number of evaluations, and third in terms of CPU time, and fails on \(11\) problems.
In terms of CPU time, the two variants that use \(\beta = 0\) are faster because they do not need to compute \(\|B_k\|\) at each iteration.

\begin{figure}[ht]%
  \centering
  \includegraphics[width=.32\linewidth]{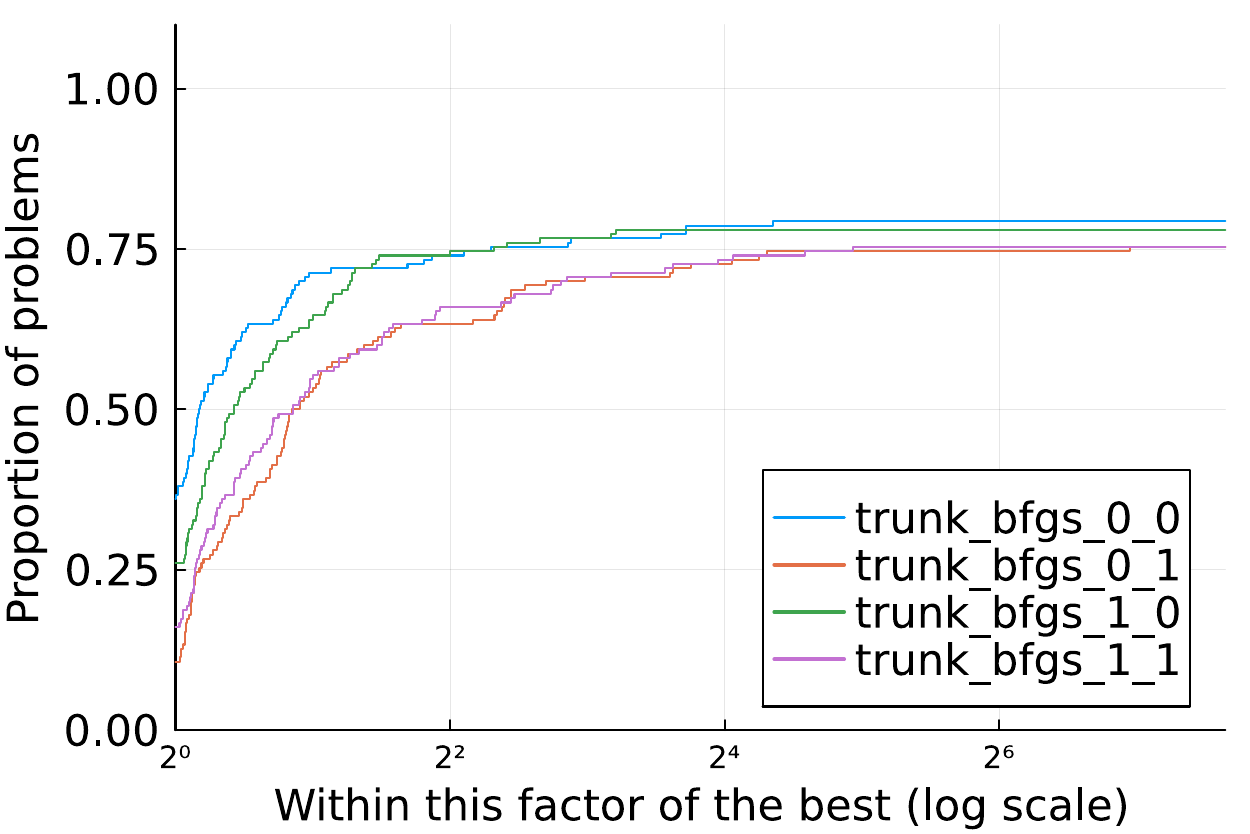}
  \hfill
  \includegraphics[width=.32\linewidth]{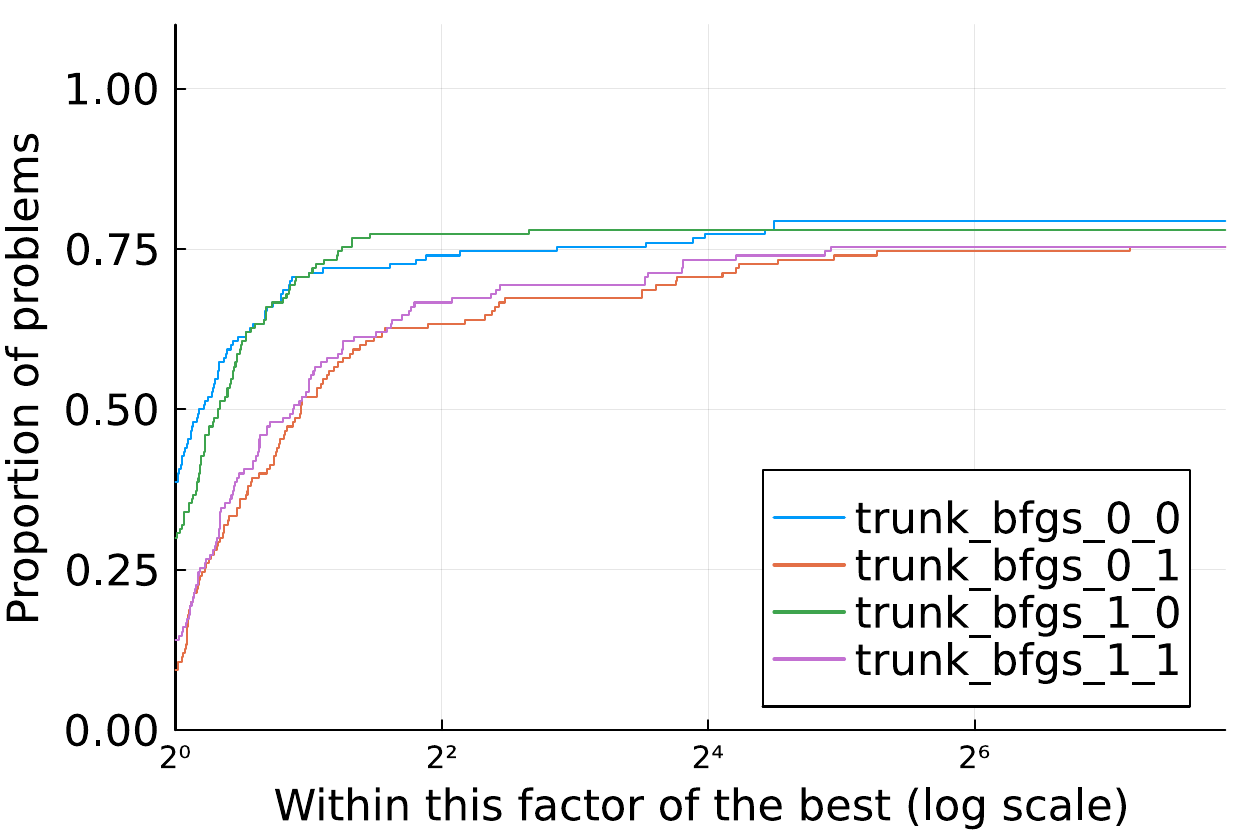}
  \hfill
  \includegraphics[width=.32\linewidth]{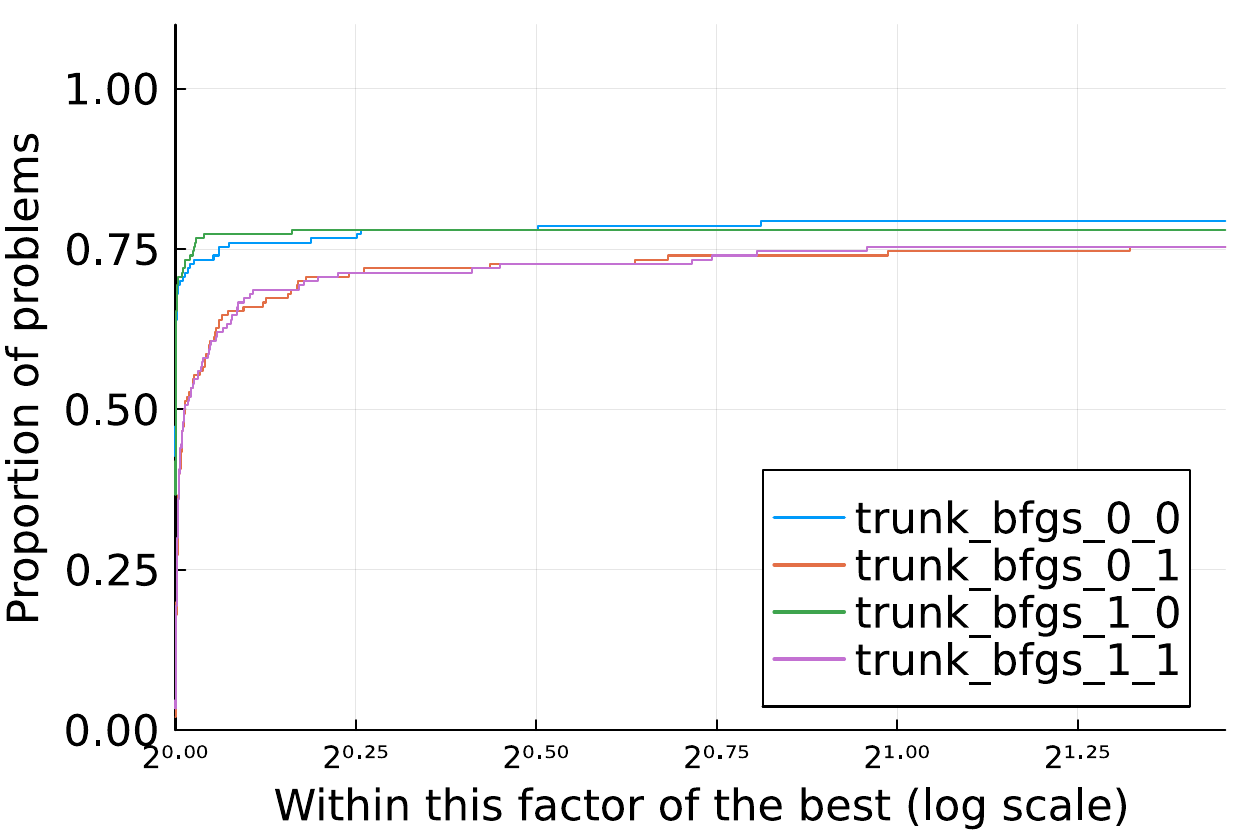}
  \caption{%
    \label{fig:profiles-bfgs}
    Performance profiles of \emph{trunk} with LBFGS Hessian approximation in terms of number of objective evaluations (left), gradient evaluations (center), and CPU time (right).
    Method \texttt{trunk\_bfgs\_a\_b} refers to the variant \(\alpha = a\) and \(\beta = b\).
  }
\end{figure}

\begin{figure}[ht]%
  \centering
  \includegraphics[width=.32\linewidth]{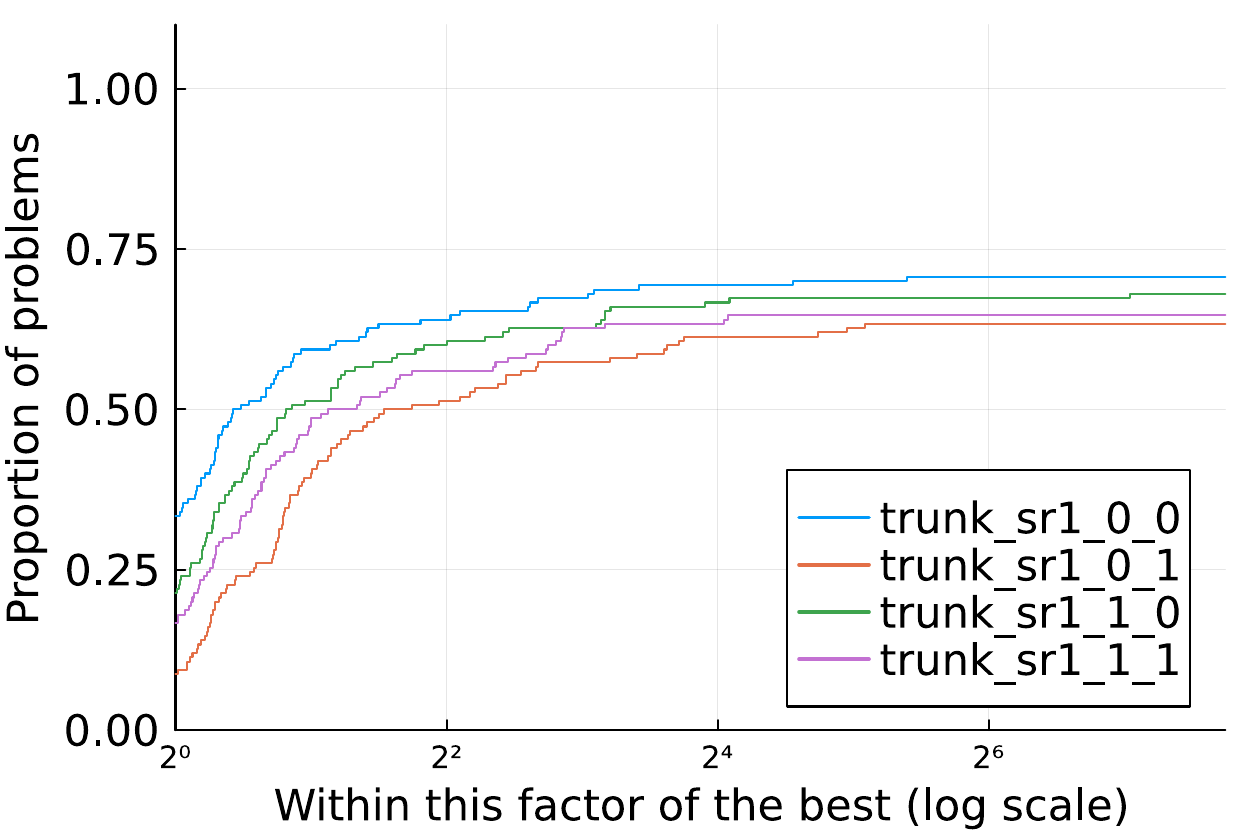}
  \hfill
  \includegraphics[width=.32\linewidth]{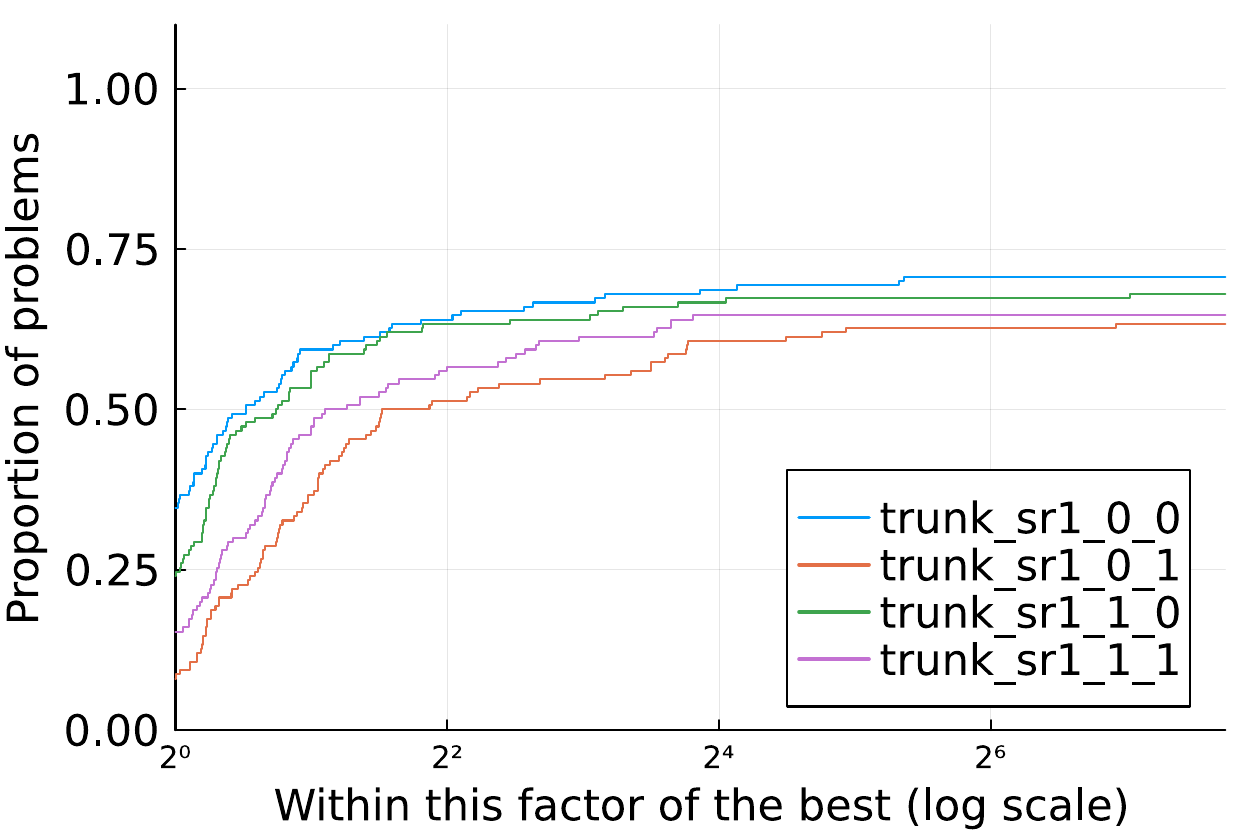}
  \hfill
  \includegraphics[width=.32\linewidth]{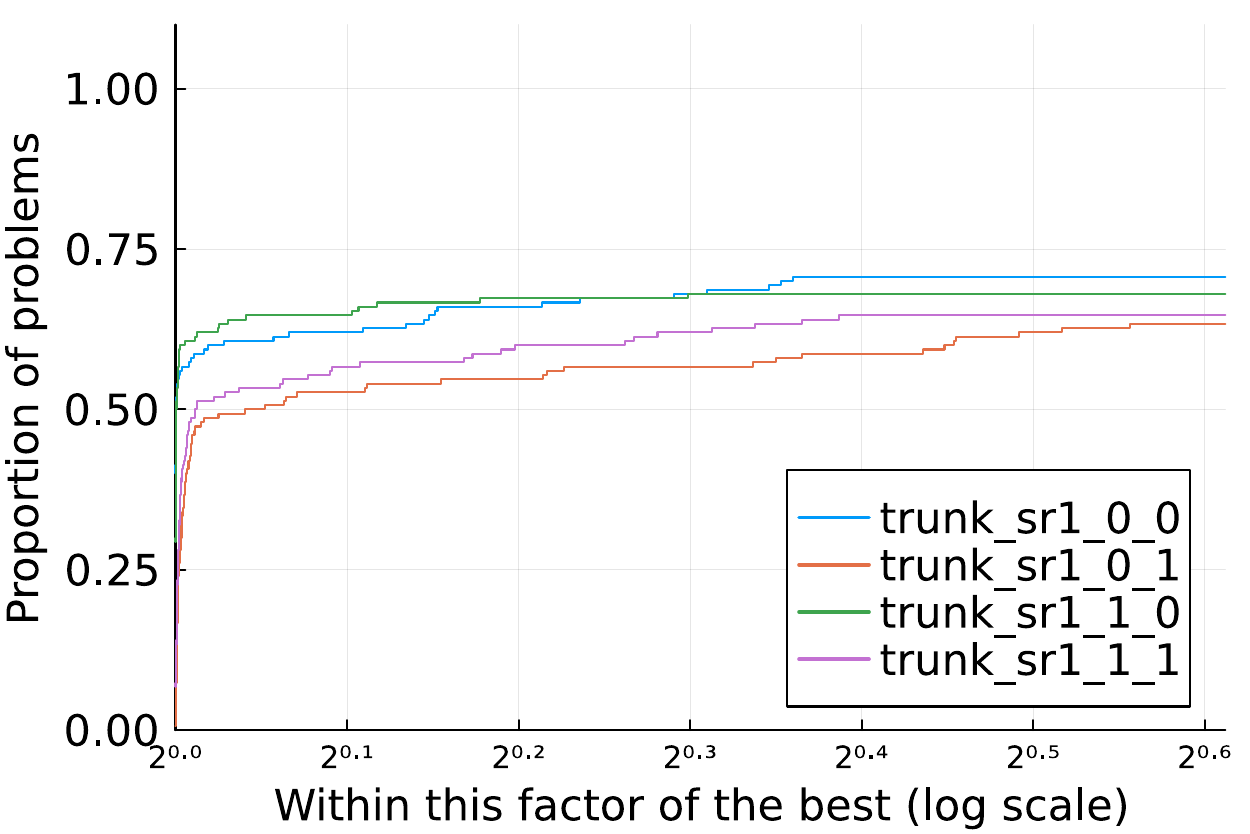}
  \caption{%
    \label{fig:profiles-sr1}
    Performance profiles of \emph{trunk} with LSR1 Hessian approximation in terms of number of objective evaluations (left), gradient evaluations (center), and CPU time (right).
    Method \texttt{trunk\_sr1\_a\_b} refers to the variant \(\alpha = a\) and \(\beta = b\).
  }
\end{figure}

Because our analysis allows for Hessian approximations, we ran the same experiment using model Hessians defined by limited-memory BFGS or SR1 approximations with memory \(5\).
Although these approximations remain bounded under mild assumptions (see \Cref{sec:S}), the objective of this section is to investigate whether a connection exists between the worst-case complexity bounds and the practical performance of \Cref{alg:TR-GEN} for various values of \(\alpha\) and \(\beta\).
The results are reported in \Cref{fig:profiles-bfgs,fig:profiles-sr1}.
The variants are color-coded as in \Cref{fig:profiles-trunk}, and the trend is the same.

\Cref{fig:profiles-trunk,fig:profiles-bfgs,fig:profiles-sr1} suggest no clear relationship between worst-case complexity and performance in practice, at least on this preliminary, and admittedly restricted, test set.
To the authors, this situation is reminiscent of comparisons between the standard trust-region implementation \texttt{trunk\_0\_0} and adaptive cubic regularization (ARC) methods, which became popular due to their favorable \(O(\epsilon^{-3/2})\) worst-case evaluation complexity.
ARC methods are also covered by the analysis of \citet{cartis-gould-toint-2022}, and those authors are responsible for much of their development.
In the experiments conducted by \citet{dussault-migot-orban-2024}, however, the performance of ARC methods is not sufficiently compelling to become a method of choice in practice.
Future developments may, of course, have us revise that opinion.

\section{Discussion}%
\label{sec:discussion}

We extended the complexity analysis of trust-region methods to handle potentially unbounded model Hessians in unconstrained optimization.
Unlike traditional complexity analyses that assume uniformly bounded model Hessians, our study covers practical cases, including quasi-Newton updates such as PSB, BFGS, and SR1.
We analyzed two regimes of the model Hessian growth: linear in the number of successful iterations and in the total number of iterations.
In both cases, we derived sharp bounds for the number of iterations required to reach an \(\epsilon\)-stationary point.
Among other contributions, our results address the intuition of \citet{powell-2010} regarding complexity for several quasi-Newton approximations.
Although he made the assumption \(\|B_k\| = O(k)\), his analysis only considers \(B_k\) updates on successful iterations, which suggests to us that he really meant to assume that \(\|B_k\| = O(|\mathcal{S}_{k-1}|)\).
\citeauthor{powell-2010}'s intuition is that the complexity would be of double exponential order, which is significantly worse than our sharp exponential bound.

Ongoing work aims to extend our analysis to cover the adaptive regularization with cubics (ARC) framework \citep{cartis-gould-toint-2011}, and investigate whether full quasi-Newton approximations may indeed grow linearly with the number of (successful) iterations.
We will also investigate the existence of a minimization problem that shows that the complexity bounds derived in the (strongly) convex case are sharp.

\subsubsection*{Acknowledgements}

We express our sincere gratitude to Coralia Cartis for insightful discussions that improved our arguments, and in particular, for bringing \citep{powell-2010} to our attention.

\small
\bibliographystyle{abbrvnat}
\bibliography{abbrv,unbounded}

\begin{thebibliography}{30}
\providecommand{\natexlab}[1]{#1}
\providecommand{\url}[1]{\texttt{#1}}
\expandafter\ifx\csname urlstyle\endcsname\relax
  \providecommand{\doi}[1]{doi: #1}\else
  \providecommand{\doi}{doi: \begingroup \urlstyle{rm}\Url}\fi

\bibitem[Aravkin et~al.(2021)Aravkin, Baraldi, and
  Orban]{aravkin-baraldi-orban-2021}
A.~Y. Aravkin, R.~Baraldi, and D.~Orban.
\newblock \href{https://arxiv.org/abs/2103.15993v1}{A proximal quasi-{N}ewton
  trust-region method for nonsmooth regularized optimization}.
\newblock Preliminary Report arXiv:2103.15993v1, 2021.

\bibitem[Burdakov et~al.(2017)Burdakov, Gong, Zirkin, and
  Yuan]{burdakov-gong-zirkin-yuan-2017}
O.~Burdakov, L.~Gong, S.~Zirkin, and Y.-X. Yuan.
\newblock \doilink{10.1007/s12532-016-0109-7}{On efficiently combining
  limited-memory and trust-region techniques}.
\newblock \emph{Math. Program. Comp.}, 9:\penalty0 101--134, 2017.

\bibitem[Cartis et~al.(2011)Cartis, Gould, and Toint]{cartis-gould-toint-2011}
C.~Cartis, N.~I.~M. Gould, and {\relax Ph}.~L. Toint.
\newblock \doilink{10.1007/s10107-009-0286-5}{Adaptive cubic regularisation
  methods for unconstrained optimization. {P}art {I}: motivation, convergence
  and numerical results}.
\newblock \emph{Math. Program.}, 127\penalty0 (2):\penalty0 245--295, 2011.

\bibitem[Cartis et~al.(2012)Cartis, Gould, and Toint]{cartis-gould-toint-2012}
C.~Cartis, N.~I.~M. Gould, and {\relax Ph}.~L. Toint.
\newblock \doilink{10.1080/10556788.2011.602076}{Evaluation complexity of
  adaptive cubic regularization methods for convex unconstrained optimization}.
\newblock \emph{Optim. Method Softw.}, 27\penalty0 (2):\penalty0 197--219,
  2012.

\bibitem[Cartis et~al.(2022)Cartis, Gould, and Toint]{cartis-gould-toint-2022}
C.~Cartis, N.~I.~M. Gould, and {\relax Ph}.~L. Toint.
\newblock \doilink{10.1137/1.9781611976991}{\emph{Evaluation Complexity of
  algorithms for nonconvex optimization}}.
\newblock Number~30 in MOS-SIAM Series on Optimization. SIAM, Philadelphia,
  USA, 2022.

\bibitem[Conn et~al.(2000)Conn, Gould, and Toint]{conn-gould-toint-2000}
A.~R. Conn, N.~I.~M. Gould, and {\relax Ph}.~L. Toint.
\newblock \doilink{10.1137/1.9780898719857}{\emph{Trust-Region Methods}}.
\newblock Number~1 in MPS-SIAM Series on Optimization. SIAM, Philadelphia, USA,
  2000.

\bibitem[Curtis et~al.(2018)Curtis, Lubberts, and
  Robinson]{curtis-lubberts-robinson-2018}
F.~E. Curtis, Z.~Lubberts, and D.~P. Robinson.
\newblock \doilink{doi.org/10.1007/s11590-018-1286-2}{Concise complexity
  analyses for trust region methods}.
\newblock \emph{Optim. Lett.}, 12:\penalty0 1713--1724, 2018.

\bibitem[Dennis and Mor\'{e}(1977)]{dennis-more-1977}
J.~E. Dennis, Jr. and J.~J. Mor\'{e}.
\newblock \doilink{10.1137/1019005}{Quasi-{N}ewton methods, motivation and
  theory}.
\newblock \emph{SIAM Rev.}, 19\penalty0 (1):\penalty0 46--89, 1977.

\bibitem[Dolan and Mor\'e(2002)]{dolan-more-2002}
E.~D. Dolan and J.~J. Mor\'e.
\newblock \doilink{10.1007/s101070100}{Benchmarking optimization software with
  performance profiles}.
\newblock \emph{Math. Program., Series~B}, 29\penalty0 (2):\penalty0 201--213,
  2002.

\bibitem[Dussault et~al.(2024)Dussault, Migot, and
  Orban]{dussault-migot-orban-2024}
J.-P. Dussault, T.~Migot, and D.~Orban.
\newblock \doilink{10.1007/s10107-023-02007-6}{Scalable adaptive cubic
  regularization methods}.
\newblock \emph{Math. Program.}, 207:\penalty0 191--225, 2024.

\bibitem[Fan and Yuan(2001)]{fan-yuan-2001}
J.~Fan and Y.~Yuan.
\newblock \href{https://ftp.cc.ac.cn/pub/home/yyx/papers/p014.pdf}{A new trust
  region algorithm with trust region radius converging to zero}.
\newblock In \emph{Proc. 5th Int. Conf. Optimization: Techniques and
  Applications}, volume~4, pages 786--794. World Scientific Publishing, 12
  2001.
\newblock ISBN 962-86475-4-7.

\bibitem[Fletcher(1972)]{fletcher-1972}
R.~Fletcher.
\newblock \doilink{10.1007/BF01584540}{An algorithm for solving linearly
  constrained optimization problems}.
\newblock \emph{Math. Program.}, \penalty0 (2):\penalty0 133--165, 1972.

\bibitem[Garmanjani(2022)]{garmanjani-2022}
R.~Garmanjani.
\newblock \doilink{10.1080/02331934.2020.1830088}{A note on the worst-case
  complexity of nonlinear stepsize control methods for convex smooth
  unconstrained optimization}.
\newblock \emph{Optimization}, 71\penalty0 (6):\penalty0 1709--1719, 2022.

\bibitem[Grapiglia et~al.(2014)Grapiglia, Yuan, and Yuan]{grapiglia-2014}
G.~Grapiglia, J.~Yuan, and Y.~Yuan.
\newblock \doilink{10.1007/s10107-014-0794-9}{On the convergence and worst-case
  complexity of trust-region and regularization methods for unconstrained
  optimization}.
\newblock \emph{Math. Program.}, 152:\penalty0 491--520, 07 2014.

\bibitem[Grapiglia et~al.(2016)Grapiglia, Yuan, and Yuan]{grapiglia-yuan-2016}
G.~N. Grapiglia, J.~Yuan, and Y.~Yuan.
\newblock \doilink{10.1080/10556788.2015.1130129}{On the worst-case complexity
  of nonlinear stepsize control algorithms for convex unconstrained
  optimization}.
\newblock \emph{Optim. Method Softw.}, 31\penalty0 (3):\penalty0 591--604,
  2016.

\bibitem[Kanzow and Steck(2023)]{kanzow-steck-2023}
C.~Kanzow and D.~Steck.
\newblock \doilink{10.1007/s12532-023-00238-4}{Regularization of limited memory
  quasi-{N}ewton methods for large-scale nonconvex minimization}.
\newblock \emph{Math. Program. Comp.}, 15\penalty0 (3):\penalty0 417--444,
  2023.

\bibitem[Leconte and Orban(2024)]{leconte-orban-2023}
G.~Leconte and D.~Orban.
\newblock \doilink{10.13140/RG.2.2.22451.40486}{Complexity of trust-region
  methods with unbounded {H}essian approximations for smooth and nonsmooth
  optimization}.
\newblock Cahier du GERAD G-2023-65, GERAD, Montr\'eal, QC, Canada, 2024.

\bibitem[Lu(1996)]{lu-1996}
X.~Lu.
\newblock \emph{A Study of the Limited Memory SR1 Method in Practice}.
\newblock PhD thesis, University of Colorado, 1996.
\newblock \url{https://books.google.ca/books?id=nkyZNwAACAAJ}.

\bibitem[Migot et~al.(2024{\natexlab{a}})Migot, Orban, Soares~Siqueira, and
  contributors]{jso}
T.~Migot, D.~Orban, A.~Soares~Siqueira, and contributors.
\newblock \doilink{10.5281/zenodo.2655082}{{The JuliaSmoothOptimizers Ecosystem
  for Numerical Linear Algebra and Optimization in Julia}}, July
  2024{\natexlab{a}}.

\bibitem[Migot et~al.(2024{\natexlab{b}})Migot, Orban, Soares~Siqueira, and
  contributors]{jso_solvers}
T.~Migot, D.~Orban, A.~Soares~Siqueira, and contributors.
\newblock
  \href{https://github.com/JuliaSmoothOptimizers/JSOSolvers.jl}{{JSOSolvers.jl:
  JuliaSmoothOptimizers optimization solvers}}, July 2024{\natexlab{b}}.

\bibitem[Migot et~al.(2024{\natexlab{c}})Migot, Orban, Soares~Siqueira, and
  contributors]{optimization_problems}
T.~Migot, D.~Orban, A.~Soares~Siqueira, and contributors.
\newblock
  \href{https://github.com/JuliaSmoothOptimizers/OptimizationProblems.jl}{{OptimizationProblems.jl:
  A collection of optimization problems in Julia}}, July 2024{\natexlab{c}}.

\bibitem[Nesterov(2013)]{nesterov-2013}
Y.~Nesterov.
\newblock \doilink{10.1007/978-1-4419-8853-9}{\emph{Introductory lectures on
  convex optimization: A basic course}}, volume~87.
\newblock Springer Verlag, Berlin, 2013.

\bibitem[Nocedal and Wright(1999)]{nocedal-wright-1999}
J.~Nocedal and S.~J. Wright.
\newblock \doilink{10.1007/b98874}{\emph{Numerical optimization}}.
\newblock Springer Verlag, Berlin, 2nd edition, 1999.

\bibitem[Powell(1970)]{powell-1970}
M.~J.~D. Powell.
\newblock \doilink{10.1016/B978-0-12-597050-1.50006-3}{A new algorithm for
  unconstrained optimization}.
\newblock In J.~B. Rosen, O.~L. Mangasarian, and K.~Ritter, editors,
  \emph{Nonlinear Programming}, pages 31--65. Academic Press, 1970.

\bibitem[Powell(1975)]{powell-1975}
M.~J.~D. Powell.
\newblock \doilink{10.1016/B978-0-12-468650-2.50005-5}{Convergence properties
  of a class of minimization algorithms}.
\newblock In O.~L. Mangasarian, R.~R. Meyer, and S.~M. Robinson, editors,
  \emph{Nonlinear Programming 2}, pages 1--27. Academic Press, 1975.

\bibitem[Powell(1984)]{powell-1984}
M.~J.~D. Powell.
\newblock \doilink{10.1007/BF02591998}{On the global convergence of trust
  region algorithms for unconstrained minimization}.
\newblock \emph{Math. Program.}, \penalty0 (29):\penalty0 297--303, 1984.

\bibitem[Powell(2010)]{powell-2010}
M.~J.~D. Powell.
\newblock \doilink{10.1093/imanum/drp021}{On the convergence of a wide range of
  trust region methods for unconstrained optimization}.
\newblock \emph{IMA J. Numer. Anal.}, 30\penalty0 (1):\penalty0 289--301, 2010.

\bibitem[Steihaug(1983)]{steihaug-1983}
T.~Steihaug.
\newblock \doilink{10.1137/0720042}{The conjugate gradient method and trust
  regions in large scale optimization}.
\newblock \emph{SIAM J. Numer. Anal.}, 20\penalty0 (3):\penalty0 626--637,
  1983.

\bibitem[Toint(1988)]{toint-1988}
{\relax Ph}.~Toint.
\newblock \doilink{10.1093/imanum/8.2.231}{{Global Convergence of a Class of
  Trust-Region Methods for Nonconvex Minimization in Hilbert Space}}.
\newblock \emph{IMA J. Numer. Anal.}, 8\penalty0 (2):\penalty0 231--252, 1988.

\bibitem[Toint(2011)]{toint-2011}
{\relax Ph}.~Toint.
\newblock \doilink{10.1080/10556788.2011.610458}{Nonlinear stepsize control,
  trust regions and regularizations for unconstrained optimization}.
\newblock \emph{Optim. Method Softw.}, 28:\penalty0 1--14, 01 2011.

\end{thebibliography}
\normalsize

\appendix
\clearpage

\end{document}